\documentclass[12pt,twoside]{amsart}

\def\ECH{\operatorname{ECH}}
\def\SWF{\operatorname{SWF}}
\def\PSS{\scriptscriptstyle \operatorname{PSS}}
\def\cyl{\operatorname{cyl}}
\def\good{\operatorname{good}}
\def\cuplength{\operatorname{cl}}
                       
\usepackage{amssymb,latexsym,bbm,graphicx,epsfig,epic,eepic,oldgerm,psfrag,hhline}
\usepackage{a4wide}
\usepackage{amscd}
\usepackage{color}
\usepackage{enumerate}

\usepackage{xypic} 
\usepackage{xfrac}
\input xy
 
\xyoption{all}

\theoremstyle{plain}
\newtheorem{theorem}{Theorem}[section]

\newtheorem{corollary}[theorem]{Corollary}
\newtheorem*{remark*}{Remark}
\newtheorem*{remarks*}{Remarks}
\newtheorem{remark}[theorem]{Remark}

\newtheorem{example}[theorem]{Example}
\newtheorem{examples}[theorem]{Examples}
\newtheorem*{example*}{Example}
\newtheorem*{examples*}{Examples}

\newtheorem*{definition*}{Definition}

\newtheorem*{questions*}{Questions}

\newtheorem{assumption}[theorem]{Assumption}

\newtheorem{Arnold.Ham}[theorem]{Hamiltonian Arnol'd conjecture}

\numberwithin{figure}{section}
\numberwithin{equation}{section}

\newcommand{\proofend}{\hspace*{\fill} $\Box$\\}
\newcommand{\diam}{\hspace*{\fill} $\Diamond$}

\def\1{\:\!}
\def\2{\;\!}
\def\s{\smallskip}
\def\m{\medskip}
\def\eps{\varepsilon}
\def\ss{\operatorname{s}}
\def\u{\operatorname{u}}

\def\im{\operatorname {im}}

\def\pt{\operatorname{pt}}
\def\Vol{\operatorname {Vol}\:\!}

\def\Diffc0{\operatorname{Diff^c_0}}

\def\Sympc0{\operatorname{Symp^c_0}}

\def\Ham{\operatorname{Ham}}

\def\rank{\operatorname{rank}}

\def\multispec{\operatorname{multispec}}

\def\ind{\operatorname{ind}}
\def\Crit{\operatorname{Crit}}

\def\can{\operatorname{can}}
\def\st{\operatorname{st}}

\def\top{\operatorname{top}}

\def\vol{\operatorname{vol}}

\def\H{\operatorname{H}}

\def\SL{\operatorname{SL}}

\def\aa{\boldsymbol{a}}
\def\kk{\boldsymbol{k}}

\def\odd{\operatorname{odd}}
\def\even{\operatorname{even}}

\def\sCZ{ {\scriptscriptstyle\operatorname{CZ}}}
\def\CM{\operatorname{CM}}
\def\HM{\operatorname{HM}}
\def\CF{\operatorname{CF}}
\def\HF{\operatorname{HF}}
\def\HW{\operatorname{HW}}
\def\HC{\operatorname{HC}}
\def\SH{\operatorname{SH}}

\def\RFH{\operatorname{RFH}}

\def\ECC{\operatorname{ECC}} 
\def\ECH{\operatorname{ECH}}  

\def\ga{\alpha}
\def\gb{\beta}
\def\gg{\gamma}

\def\gve{\varepsilon}
\def\gf{\varphi}

\def\go{\omega}
\def\gs{\sigma}

\def\B{\operatorname{B}}

\def\E{\operatorname{E}}

\def\Z{\operatorname{Z}}

\def\ca{{\mathcal A}}
\def\cb{{\mathcal B}}

\def\ce{{\mathcal E}}

\def\cm{{\mathcal M}}

\def\cp{{\mathcal P}}

\def\cs{{\mathcal S}}

\def\CC{\mathbbm{C}}
\def\DD{\mathbbm{D}}

\def\NN{\mathbbm{N}}

\def\QQ{\mathbbm{Q}}
\def\RR{\mathbbm{R}}

\def\TT{\mathbbm{T}}
\def\ZZ{\mathbbm{Z}}

\def\R{\operatorname{\mathbbm{R}}}
\def\RP{\operatorname{\mathbbm{R}P}}
\def\CP{\operatorname{\mathbbm{C}P}}

\def\SU{\operatorname{SU}}

\def\pp{\partial}

\def\ni{\noindent}
\def\b{\bigskip}
\def\m{\medskip}

\def\id{\mbox{id}}

\def\proof{\noindent {\it Proof. \;}}

\begin{document}

\title[Floer homologies]
{Floer homologies, with applications}

\author{Alberto Abbondandolo}
\thanks{This work is part of AA's activities within the DFG collaboration scheme SFB CRC 191 ``Symplectic structures in geometry, algebra and dynamics''.}
\address{Alberto Abbondandolo,
  Fakult\"at f\"ur Mathematik, 
  Ruhr-Universit\"at Bochum}
\email{alberto.abbondandolo@rub.de}

\author{Felix Schlenk}  
\thanks{FS partially supported by SNF grant 200020-144432/1.}
\address{Felix Schlenk,
Institut de Math\'ematiques,
Universit\'e de Neuch\^atel}
\email{schlenk@unine.ch}

\keywords{Floer homology}

\date{\today}
\thanks{2000 {\it Mathematics Subject Classification.}
Primary 53D40, Secondary~37J45, 53D35}

\begin{abstract}
Floer invented his theory in the mid eighties in order to prove the Arnol'd conjectures
on the number of fixed point of Hamiltonian diffeomorphisms and Lagrangian intersections. 
Over the last thirty years, many versions of Floer homology have been constructed. 
In symplectic and contact dynamics and geometry they have become a principal tool,
with applications that go far beyond the Arnol'd conjectures:
The proof of the Conley conjecture and of many instances of the Weinstein conjecture, 
rigidity results on Lagrangian submanifolds and on the group of symplectomorphisms, 
lower bounds for the topological entropy of Reeb flows and
obstructions to symplectic embeddings are just some of the applications of Floer's seminal ideas.
Other Floer homologies are of topological nature.
Among their applications are Property~P for knots
and the construction of compact topological manifolds of dimension greater than five 
that are not triangulisable.

This is by no means a comprehensive survey on the presently known Floer homologies and their applications. 
Such a survey would take several hundred pages. 
We just describe some of the most classical versions and applications, 
together with the results that we know or like best.
The text is written for non-specialists and the focus is on ideas rather than generality. 
Two intermediate sections recall basic notions and concepts from symplectic dynamics and geometry.
\end{abstract}

\maketitle
\tableofcontents


\section{Introduction}  \label{s:intro}

In classical Morse theory, 
the structure of a smooth manifold~$M$ is related to the critical
points of any Morse function~$f$ on~$M$. 
The manifold~$M$ can be successively reconstructed 
by attaching for each critical point~$x$ a handle whose dimension is the Morse index of~$x$.
Another approach to Morse theory is Morse homology. 
Here, one chooses a Riemannian metric on~$M$ and recovers the homology of~$M$ 
by looking at the critical points of a Morse function~$f$ and at the flow lines
of~$-\nabla f$ connecting critical points of index difference one. 

While the first approach can be used in several infinite dimensional situations, 
such as the energy functional on the loop space, it fails when the functional 
in question has only critical points of infinite index. 
In some important situations however, after suitable renormalisation of the index,  
the index {\it difference}\/ between any two critical points 
is still finite. One can then try to build a Morse homology. 
Such a homology is called a Floer homology.

Floer built such homologies in two situations: 
For the action functional of classical mechanics on the space of curves on a symplectic manifold 
and for the Chern--Simons functional on the space of $\SU (2)$ connections 
on a 3-dimensional oriented integral homology sphere~$Y$.
The first theory is called Hamiltonian or Lagrangian Floer homology - 
depending on whether one considers closed or open curves 
- and proves at once important cases of Arnol'd's conjectures on the number 
of fixed points of Hamiltonian diffeomorphisms and of Lagrangian intersection points.
The second one, called instanton homology, refines a topological invariant:
Its Euler characteristic is twice the Casson invariant of~$Y$, 
namely the number of conjugacy classes of representations $\pi_1(Y) \to \SU (2)$.
By now, many other Floer homologies have been constructed, 
both of symplectic and of topological nature,
and they have many more applications to dynamics and geometry.

Our survey starts with a discussion of Morse homology on a finite dimensional compact manifold,
since this homology is the model for every Floer homology.
All features of Morse homology have a counterpart in every Floer homology, 
and the first thing one should do if one wants to check if a Floer homology has a certain property 
is to see if this property holds for Morse homology.
In Section~\ref{s:Ham} we then describe in some detail the arguably easiest version of Floer homology, 
namely Hamiltonian Floer homology for symplectically aspherical closed symplectic manifolds.
Many aspects and technical issues of other Floer homologies are analogous, so that in our subsequent discussion of
other Floer homologies we can give less details on the construction and
focus on their applications.

Symplectic homology (Section~\ref{s:symphom}) is an invariant of suitable subsets of a symplectic manifold.
It leads to a proof of the Weinstein conjecture on the existence of a closed Hamiltonian orbit on an energy surface, 
and it provides obstructions to the existence of Lagrangian submanifolds
and of certain symplectic embeddings.
Lagrangian Floer homology (Section~\ref{s:Lag}) gives lower bounds on the number of intersections of certain pairs of
Lagrangian submanifolds and can be used to prove that for many contact manifolds all Reeb flows must have 
positive topological entropy.
Finally, in Section~\ref{s:contact} we describe two Floer homologies for contact manifolds.
Contact homology can be used to distinguish contact structures on compact contact manifolds.
Embedded contact homology, that is defined for three-dimensional compact contact manifolds,
can be used to prove the Weinstein conjecture in dimension three, and it gives rise to numerical invariants 
that yield very fine symplectic embedding obstructions in dimension four 
and imply the $C^\infty$~closing lemma for Hamiltonian diffeomorphisms on surfaces.

These are all ``symplectic'' Floer homologies and applications to symplectic and contact
dynamics and geometry.
Other Floer homologies have applications to the structure of manifolds, 
like Property~P for knots, 
the question which compact 3-manifolds can be obtained by Dehn surgery of~$S^3$ on which knot,
and the existence of topological manifolds of dimension~$\geqslant 5$ that admit no triangulation. 
Given Manolescu's excellent survey \cite{Ma15}, we just state a few of these results in Section~\ref{s:top}.

%

\m \ni
{\bf Acknowledgment.}
AA cordially thanks the hospital at Visp for the excellent working conditions 
during his ski vacation.
We are grateful to Marcelo Alves, Lucas Dahinden, Carsten Haug, and Pedram Safaee for interesting discussions. 

\section{From Morse theory to Morse homology} \label{s:Morse}

\subsection{Classical Morse Theory} \label{s:classical}
Consider a smooth closed (that is, compact and without boundary) manifold~$M$ of dimension~$d$. Given a smooth function $f \colon M \to \RR$ let 
\[
\Crit f = \{ x \in M \mid df (x) =0\}
\]
be the set of its critical points. The function~$f$ is called {\em Morse}\/ if every $x$ 
in~$\Crit f$ is non-degenerate, meaning that the second differential of~$f$ at~$x$ is a non-degenerate bilinear form. By the Morse Lemma, this is equivalent to say that there is 
a chart~$\varphi$ near~$x$ such that
\begin{equation} \label{e:Morse}
f \circ \varphi^{-1} (y_1, \dots, y_d) \,=\, f(x) - ( y_1^2 + \dots + y_k^2) + ( y_{k+1}^2 + \dots + y_d^2) .
\end{equation}
The number $k$ in this expression does not depend on the chart and coincides with the number of negative eigenvalues of the self adjoint endomorphism of~$T_x M$ representing the second differential of~$f$ at~$x$ with respect to some inner product. 
It is called the {\em Morse index}\/ of~$f$ at~$x$, and we denote it by~$\ind_f(x)$. 
The normal form~\eqref{e:Morse} (or the inverse mapping theorem applied to~$df$) shows that 
non-degenerate critical points are isolated. 
In particular, the critical set~$\Crit f$ of a Morse function on a closed manifold is finite.

Denote by $c_k(f)$ the number of critical points of~$f$ with Morse index~$k$ and by~$b_k(M)$ 
the $k$-th Betti number of~$M$, that is, the rank of the $k$-th homology group~$\H_k(M)$. 
Here $\H_*$ denotes singular homology with integer coefficients. 
Morse theory relates the critical set of~$f$ to the topology of~$M$ by the following relationship between the non-negative integers $c_k(f)$ and~$b_k(M)$:

\begin{theorem} \label{t:Morse}
For any Morse function $f$ on the closed manifold~$M$ there exists a polynomial~$Q$ with 
non-negative integer coefficients such that
\begin{equation}
\label{morse}
\sum_{k=0}^d c_k(f)\, z^k \,=\, \sum_{k=0}^d b_k (M)\, z^k + (1+z)\, Q(z).
\end{equation}
In particular, $c_k(f) \geqslant b_k(M)$ for every $k$.
\end{theorem}

This theorem is proved by analyzing the change of 
topology of the sublevel sets 
$\{ x \in M \mid f(x) \leqslant a\}$ when $a$ crosses a critical level 
and by using elementary properties of singular homology.

An alternative equivalent way to state the above theorem is as follows. Let 
\[
\CM_k(f) \cong  \ZZ^{c_k(f)}
\]
be the free Abelian group generated by the critical points of~$f$ with Morse index~$k$ 
and set $C_{-1}(f):= \{0\}$. Then there is a sequence of homomorphisms
\[
\partial_k \colon \CM_k(f) \rightarrow \CM_{k-1}(f), \qquad k\geqslant 0,
\]
such that $\partial_k \circ \partial_{k+1}=0$ and
\[
\frac{\mathrm{ker} \,\partial_k}{\mathrm{im}\, \partial_{k+1}} \cong \H_k(M)
\]
for every $k \geqslant 0$.
In other words, the Abelian groups $\CM_k(f)$ are the spaces of a chain complex whose homology gives the singular homology of~$M$. Going from one formulation to the other is a simple exercise in homological algebra, but the latter formulation looks intriguing and raises the following interesting questions: What is the meaning of the boundary operators~$\partial_k$?  
Do they have a direct definition?

\subsection{The Morse complex}  \label{s:morse}
It turns out that it is possible to define ``natural'' boundary operators~$\partial_k$ satisfying the above requirements by choosing a {\em Morse--Smale negative gradient flow}\/ for~$f$, 
as we now explain. Fix a Riemannian metric~$g$ on~$M$. The gradient~$\nabla f$ of~$f$ with respect to~$g$ is the vector field on~$M$ which is $g$-dual to~$df$:
\[
g_x (\nabla f(x), v )  \,=\,  df(x)[v], \qquad x  \in M, \; v \in T_x M .
\]
The negative gradient flow of $f$ is the flow~$\phi^s$ of the vector field $-\nabla f$; its orbits are the solutions of the ordinary differential equation
\[
\dot u(s) \,=\, - \nabla f(u(s)), \quad s \in \RR . 
\]
The critical points of $f$ are the stationary points of this flow, 
and $f$ strictly decreases 
on every non-stationary orbit thanks to the identity
\[
\frac{d}{ds} f(u(s)) \,=\, - \| \nabla f(u(s)) \|_{u(s)}^2,
\]
which holds for every orbit $u(s)=\phi^s(x)$. Here $\|\cdot\|_x$ denotes the norm induced by the scalar product~$g_x$ on~$T_x M$. The Morse condition guarantees that all the stationary points of the negative gradient flow are hyperbolic, meaning that the linearisation of the vector field at them has no purely imaginary eigenvalues, and this in turn implies that each critical point~$x$ 
of~$f$ has a well-defined stable and unstable manifold:
\[
W^{\ss}(x) := \left\{ z \in M \mid \lim_{s \to +\infty} \phi^s(z) = x \right\}, 
\qquad
W^{\u}(x) := \left\{ z \in M \mid \lim_{s \to -\infty} \phi^s(z) = x \right\} .
\]
These are smoothly embedded copies of $\R^{d-k}$ and $\R^k$, respectively, where $k$ is the Morse index of~$x$. The flow of $-\nabla f$ is said to satisfy the {\em Morse--Smale condition}\/ when the stable and unstable manifolds of any two critical points intersect transversally. 
Recall that two submanifolds $S_1, S_2$ of~$M$ are said to intersect transversally
if $T_z S_1 + T_z S_2 = T_z M$ for every $z \in S_1 \cap S_2$.
In this case, $S_1 \cap S_2$ is either empty or is a submanifold of dimension 
$\dim S_1 + \dim S_2 -d$. Since the non-empty intersection of the unstable and the stable manifold of two distinct critical points is at least 1-dimensional (because it contains non constant orbits), the Morse--Smale condition implies that $W^{\u}(x) \cap W^{\ss}(y)$ 
is empty when $\ind_f(x) \leqslant \ind_f(y)$ and $x \neq y$.

It can be proved that for a generic choice of the Riemannian metric~$g$ the flow of~$-\nabla f$ is Morse--Smale. A Morse--Smale negative gradient flow for~$f$ induces a nice {\em cellular filtration}\/ of~$M$, whose associated boundary operators are the homomorphisms~$\partial_k$ 
we are looking for. Recall that a (finite) cellular filtration of a topological space~$M$ is a family of open subsets 
\[
\emptyset = U_{-1} \subset U_0 \subset \dots \subset U_{d-1} \subset U_d = M
\]
such that for every $k$ the homology of the pair $(U_k,U_{k-1})$ is concentrated in degree~$k$. In this case, the groups
\[
W_k \,:=\, \H_k(U_k,U_{k-1})
\]
are the spaces of a chain complex
\[
\partial_k \colon W_k \rightarrow W_{k-1}
\]
where the homomorphism $\partial_k$ is the composition of the singular boundary operator
\[
\H_k(U_k,U_{k-1}) \rightarrow \H_{k-1}(U_{k-1})
\]
with the homomorphism
\[
\H_{k-1}(U_{k-1}) \rightarrow \H_{k-1}(U_{k-1},U_{k-2})
\]
induced by inclusion. It is a standard fact in homology theory that 
$\partial_k \circ \partial_{k+1}=0$ and that the homology of the chain complex 
$\{W_k,\partial_k\}$ is isomorphic to the singular homology of~$M$.

The Morse--Smale negative gradient flow $\phi^s$ of $(f,g)$ induces a
cellular filtration~$\{U_k\}$ of~$M$, in which $U_k$ is the forward evolution by the flow of 
the union of suitable small neighbourhoods of all critical points of~$f$ with Morse index not exceeding~$k$. One easily sees that $U_k$ is obtained from~$U_{k-1}$ by attaching the $k$-handles which are given by thickening the unstable manifolds of the critical points of index~$k$, 
and therefore $\H_k (U_k,U_{k-1})$ has one free generator for every such handle and is hence isomorphic to~$\CM_k(f)$. Checking these properties uses the fact that, by the Morse--Smale condition, no negative gradient flow lines go from a critical point of index~$k$ to another one of index~$k$ or larger. The boundary operators of this cellular filtration can now be seen as homomorphisms
\[
\partial_k \colon \CM_k(f) \rightarrow \CM_{k-1}(f),
\]
and the chain complex $\{\CM_k(f),\partial_k\}$ is called the {\em Morse complex}\/ 
of the pair~$(f,g)$. Indeed, the groups~$\CM_k(f)$ do not depend on the generic metric~$g$, 
but the boundary operators~$\partial_k$ do. The homology of this chain complex is called 
the {\em Morse homology}\/ of~$f$ and is denoted by
\[
\HM_k(f) := \frac{\ker \,\partial_k}{\im \, \partial_{k+1}}.
\] 
The fact that this homology is isomorphic to the singular homology of~$M$ gives an alternative proof of Theorem~\ref{t:Morse}. This fact also shows that the Morse homology is independent of the metric~$g$, which is the reason why the metric does not appear in the notation. 
It is actually also independent of the Morse function~$f$, but it is customary to keep~$f$ 
in the notation, for instance because in the case of a non-compact finite dimensional 
manifold~$M$ and of a Morse function~$f$ on~$M$ satisfying suitable asymptotic conditions, 
an analogous construction would give a Morse homology which depends on the asymptotic behaviour 
of~$f$. A better reason is given at the end of this section.

\subsection{A closer look at the boundary operator}  \label{s:closer}
Although the definition that we have given is topological, the boundary operator~$\partial_k$ 
can also be read directly from the dynamics of the negative gradient flow of~$f$ by looking at the intersections of stable and unstable manifolds, in a way which we now explain. 

\begin{figure}[h]   
 \begin{center}
  \psfrag{x}{$x$}
  \psfrag{x'}{$x'$}
  \psfrag{y}{$y$}
  \psfrag{y'}{$y'$}
  \leavevmode\includegraphics{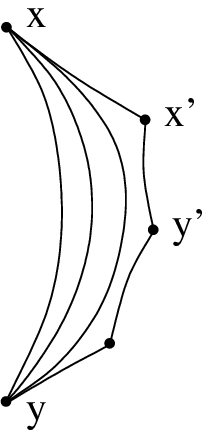}
 \end{center}
 \caption{Breaking}  \label{fig.breaking.Morse}
\end{figure}

We start by fixing an orientation of the unstable manifold of each critical point. 
This induces a co-orientation (that is, an orientation of the normal bundle) 
of the stable manifold of each critical point.
By the Morse--Smale condition, for every pair of distinct critical points~$x$ and~$y$ the set
\[
\widehat{\mathcal{M}}_{f,g}(x,y) := W^{\u}(x) \cap W^{\ss}(y),
\]
when non-empty, is a smooth submanifold of dimension $\ind_f(x) - \ind_f(y)$. Being the transverse intersection of an oriented and a co-oriented submanifold, this submanifold is oriented. Moreover, it
is invariant under the action of~$\RR$ given by the negative gradient flow, 
and we denote by $\mathcal{M}_{f,g}(x,y)$ its quotient by this action. 
Therefore, $\mathcal{M}_{f,g}(x,y)$ is the set of unparametrized negative gradient flow lines 
going from~$x$ to~$y$. It has the structure of an oriented manifold of dimension
\[
\dim \mathcal{M}_{f,g}(x,y) = \ind_f(x) - \ind_f(y) - 1.
\]
In general, the set $\widehat{\mathcal{M}}_{f,g}(x,y)$ is not closed in $M \setminus \Crit f$, 
and its closure in this space contains
points which belong to the unstable manifold of some critical point~$x'$ with $\ind_f(x') \leqslant \ind_f(x)$ (with equality iff $x'=x$) and to the stable manifold of some critical point~$y'$ with $\ind_f(y') \geqslant \ind_f(y)$ (with equality iff $y'=y$). 
It follows that the manifolds $\mathcal{M}_{f,g}(x,y)$ are in general not compact, 
but can be compactified by attaching suitable components of $\mathcal{M}_{f,g}(x',y')$ 
for some critical points~$x'$ and~$y'$ with 
\[
\ind_f(y) \leqslant \ind_f(y') < \ind_f(x') \leqslant \ind_f(x),
\]
where the first (resp.\ last) inequality is an equality iff $y'=y$ (resp.\ $x'=x$), 
see Figure~\ref{fig.breaking.Morse}. 

In the special case 
\[
\ind_f(x) - \ind_f(y)=1,
\]
there are no critical points with intermediate Morse index which can appear in the compactification, 
and hence $\mathcal{M}_{f,g}(x,y)$ is compact. Being a zero-dimensional manifold, it is a finite set, 
and being oriented, each of its finitely many elements is endowed with a coefficient~$\pm 1$. 
The algebraic sum of these coefficients defines an integer~$\nu(x,y)$. It can be proved that 
the boundary operator~$\partial_k$ which is associated to the cellular filtration induced by the negative gradient flow of~$f$ has the form
\begin{equation}
\label{bdry}
\partial_k \, x = \sum_{\substack{y\in \mathrm{Crit}\, f \\ \ind_f(y)=k-1}} \nu(x,y)\, y \qquad \forall x\in \Crit f, \; \ind_f(x)=k.
\end{equation}
This neat formula has some far reaching consequences, which are the subject of this survey. 
But before discussing them, let us look at this formula in action in some simple examples. 
In order to avoid having to consider orientations, we replace the coefficient ring~$\ZZ$ 
by the field~$\ZZ_2$:
This means that $\CM_k(f)$ is the $\ZZ_2$-vector space $\ZZ_2^{c_k(f)}$ and, when $\ind_f(x)-\ind_f(y)=1$, 
the number $\nu(x,y) \in \ZZ_2$ is just the parity of the finite set~$\mathcal{M}_{f,g}(x,y)$. 
In this case, $\HM_k(f)$ is isomorphic to the singular homology of~$M$ with $\ZZ_2$~coefficients.

\begin{example} \label{ex:sphere}
{\rm
We consider three Morse functions on the sphere~$S^2$, 
namely the height function on the round sphere, on the smooth heart, and on the duckling,
as in Figure~\ref{fig.ente}.
The Riemannian metric~$g$ is the one induced from the ambient Euclidean space.
Critical points of index two (local maxima), of index one (saddles), and of index zero (local minima)
are denoted by $a_i$, $b_i$, and $c_i$, respectively. The corresponding negative gradient flows 
are Morse--Smale.
\begin{figure}[h]  
 \begin{center}
  \psfrag{fi}{$f_i$}
  \psfrag{i}{$\mbox{(1) round sphere}$}
  \psfrag{ii}{$\mbox{(2) smooth heart}$}
  \psfrag{iii}{$\mbox{(3) duckling}$}
  \psfrag{a2}{$a_2$}
  \psfrag{a}{$a$}
  \psfrag{b}{$b$}
  \psfrag{c}{$c$}
  \psfrag{a1}{$a_1$}
  \psfrag{a2}{$a_2$}
  \psfrag{b1}{$b_2$}
  \psfrag{b2}{$b_1$}
  \psfrag{c1}{$c_1$}
  \psfrag{c2}{$c_2$}
  \leavevmode\includegraphics{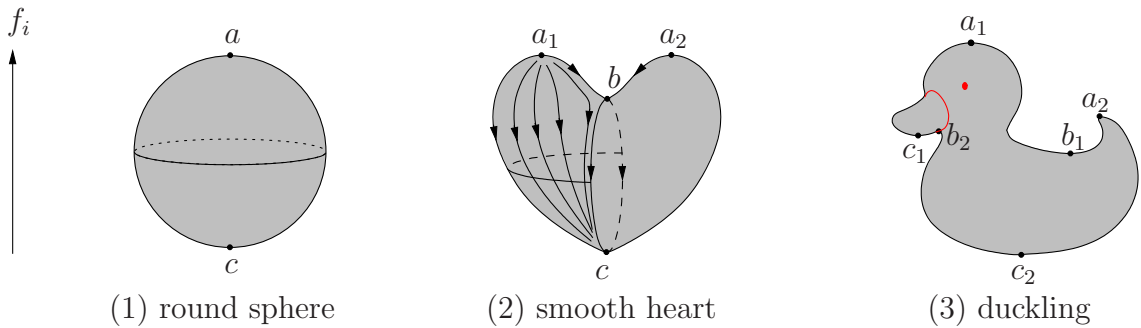}
 \end{center}
 \caption{Three Morse--Smale gradient flows on the 2-sphere}
 \label{fig.ente}
\end{figure}
%

\ni
For the round sphere, the moduli space $\cm_{f_1,g}(a,c)$ can be identified with the 
equatorial circle, consisting of the set of meridians. In this case, the Morse complex has one generator in degree zero and one generator in degree two. Since there are no pairs of critical points with index difference~1, the boundary operator vanishes identically, and hence the Morse homology is isomorphic to the Morse complex. From this, we recover the familiar fact that the singular homology of the 2-sphere is concentrated in degrees~0 and~2, in both of which it is 1-dimensional.

For the smooth heart, $\cm_{f_2,g}(a_1,b)$ and $\cm_{f_2,g}(a_2,b)$ are singletons, while
$\cm_{f_2,g}(b,c)$ consists of two points;
moreover, $\cm_{f_2,g}(a_1,c)$ and $\cm_{f_2,g}(a_2,c)$ are open intervals. The boundary operator is given by
\[
\pp a_1 = \pp a_2 = b, \qquad \pp b =  c + c = 0, \qquad \pp c =0,
\]
and we find 
\[
\HM_2(f_2) = \ZZ_2 \langle a_1+a_2 \rangle \cong \ZZ_2, \quad
\HM_1(f_2) = \ZZ_2 \langle b \rangle / \ZZ_2 \langle b \rangle = 0, \quad
\HM_0(f_2) = \ZZ_2 \langle b \rangle \cong \ZZ_2 
\]
which again agrees with the singular homology of the 2-sphere with $\ZZ_2$~coefficients. 
If we were using integer coefficients, a closer look at the orientations would show that 
$\pp b = c - c =0$.

For the duckling we have
\[
\pp a_1 = \pp a_2 = b_1, \qquad \pp b_1 =0, \qquad \pp b_2 = c_1+c_2, \qquad \pp c_1=\pp c_2= 0,
\]
which again defines a chain complex whose homology is isomorphic to $\H_*(S^2)$. 
}
\end{example}

\begin{example}
{\rm
Consider the straight torus of revolution as on the left of Figure~\ref{fig.torus}. The height function~$f$ is again a Morse function, having a critical point~$a$ of index~2, two critical points~$b_1$ and~$b_2$ of index~1, and a critical point~$c$ of index~$0$. The gradient flow of~$f$ with respect to the metric induced by the ambient Euclidean space is not Morse--Smale, because the unstable manifold of the upper saddle~$b_1$ overlaps with the stable manifold of the lower one~$b_2$, 
while the Morse--Smale condition forbids the intersection of unstable and stable manifolds of distinct critical points with the same Morse index. However, if we apply a small perturbation, for instance by tilting the torus a bit, both branches of the unstable manifold of~$b_1$ will converge to $c$ and both branches of the stable manifold of $b_2$ will converge to $a$. Then
the boundary operator is given by
\[
\partial a =  b_1 + b_1 + b_2 + b_2= 0, \qquad \partial b_1 = \partial b_2 = c + c = 0, \qquad \pp c =0.
\]
The vanishing of the boundary operator implies that the Morse homology of the function~$f$ is isomorphic to its Morse complex, that is 
\[
\HM_0(f) \cong \ZZ_2, \qquad \HM_1(f) \cong \ZZ_2\oplus \ZZ_2, \qquad \HM_1(f) \cong \ZZ_2,
\]
which is precisely the singular homology of the 2-torus with $\ZZ_2$ coefficients.
}
\begin{figure}[h]   
 \begin{center}
  \psfrag{a}{$a$}
  \psfrag{b1}{$b_1$}
  \psfrag{b2}{$b_2$}
  \psfrag{c}{$c$}
  \psfrag{g1}{$u$}
  \leavevmode\includegraphics{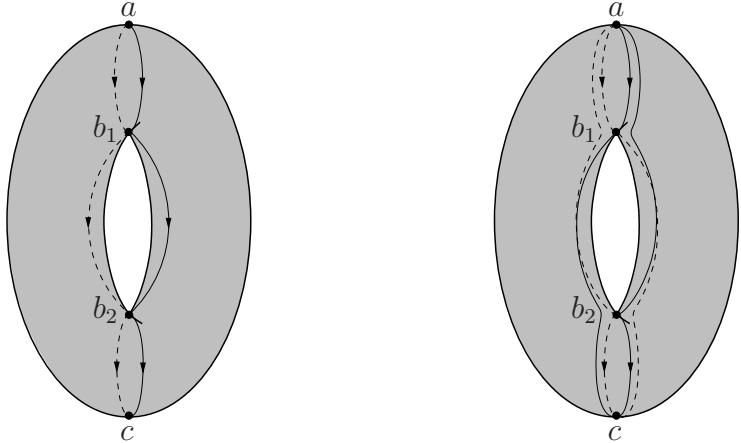}
 \end{center}
 \caption{The torus of revolution, straight and after a perturbation} \label{fig.torus}
\end{figure}
\end{example}

\subsection{Changing the point of view} 
After these examples, we start discussing the new possibilities which are opened by formula~\eqref{bdry}. 
We have deduced this formula starting from the boundary operator in cellular homology, but we could change the point of view and, applying a strategy which is ubiquitous in mathematics, see it as the 
{\em definition}\/ of a homomorphism from $\CM_k(f)$ to $\CM_{k-1}(f)$. Indeed, this formula involves only the data given by a Morse--Smale gradient flow of~$f$ and requires no prior knowledge of singular or cellular homology.

If we adopt this point of view, the first question we should ask ourselves is why the boundary property 
\[
\partial_k \circ \partial_{k+1}=0
\]
holds. Indeed, this is by no means obvious from \eqref{bdry}: Let $x$ and $z$ be critical points of~$f$ 
of Morse index~$k+1$ and~$k-1$, respectively. The coefficient of~$z$ in the expression for 
$\partial_k \circ \partial_{k+1} \,x$ given by a double application of~\eqref{bdry} is the number
\[
\sum_{\substack{y\in \Crit f\\ \ind_f(y)=k}} \nu(x,y) \, \nu(y,z),
\]
which is the algebraic count of the pairs $(u,v)$ in the set
\begin{equation}
\label{set}
\bigcup_{\substack{y\in \Crit f\\ \ind_f(y)=k}} \mathcal{M}_{f,g}(x,y) \times \mathcal{M}_{f,g}(y,z).
\end{equation}
Why is this number zero?

\begin{figure}[h]  
 \begin{center}
  \psfrag{x}{$x$}
  \psfrag{y}{$y$}
  \psfrag{y'}{$y'$}
  \psfrag{z}{$z$}
  \psfrag{u}{$u$}
  \psfrag{u'}{$u'$}
  \psfrag{v}{$v$}
  \psfrag{v'}{$v'$}
  \leavevmode\includegraphics{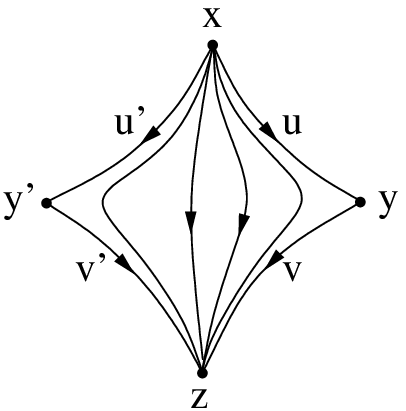}
 \end{center}
 \caption{$\pp \cm_{f,g}(x,z) = \bigcup_y \cm_{f,g}(x,y) \times \cm_{f,g}(y,z)$}
 \label{fig.gluing}
\end{figure}

This follows from a cobordism argument. Indeed, each element of the above set is a pair~$(u,v)$ of consecutive gradient flow lines -- the first from $x$ to some critical point~$y$ and the second from $y$ 
to~$z$ -- and it is possible to prove that there exists a unique 1-parameter family of gradient flow lines 
in~$\mathcal{M}(x,z)$ which converges to the ``broken gradient line'' $u\# v$. 
In other words, the set~\eqref{set} is precisely the set of boundary points of the one-dimensional 
manifold~$\mathcal{M}_{f,g}(x,y)$, see Figure~\ref{fig.gluing} or one of the two halves of the 
smooth heart in Figure~\ref{fig.ente}.
Since compact 1-dimensional manifolds have an even number of 
boundary points, 
this shows that the coefficient of $\partial_k \circ \partial_{k+1} \,x$ is zero when 
we work over~$\ZZ_2$.
Working over the integers requires checking that the two pairs~$(u,v)$ and~$(u',v')$ appearing as boundary points of a component of~$\mathcal{M}_{f,g}(x,z)$ contribute with different signs to the above sum.

This shows that the choice of a Morse function~$f$ and of a Morse--Smale gradient flow for it induces, 
via formula~\eqref{bdry}, the sequence of homology groups~$\HM_k (f)$ of the {\it Morse complex}\/ $\{ \CM_k(f),\partial_k \}$. 
These groups agree with the previously defined Morse homology of~$f$.
It is actually possible to complete this construction to a full homology theory on the category of smooth compact manifolds and smooth maps between them which satisfies the Eilenberg--Steenrod axioms. 
For instance, a smooth map
\[
\varphi \colon M_1 \rightarrow M_2
\]
between closed smooth manifolds induces a homomorphism 
\[
\varphi_* \colon \HM_*(f_1) \rightarrow \HM_*(f_2)
\]
between the Morse homologies of Morse functions $f_1 \colon M_1 \rightarrow \R$ and 
$f_2 \colon M_2 \rightarrow \R$, which is defined by a formula similar to~\eqref{bdry} 
which counts the intersections of $W^{\u}(x;-\nabla f_1)$ with $\varphi^{-1}(W^{\ss}(y;-\nabla f_2))$ 
for critical points $x \in \Crit f_1$ and $y \in \Crit f_2$ of the same Morse index. 
Once this has been done, the fact that the Morse homology of a function~$f$ is isomorphic to the singular homology of its domain~$M$ can also be deduced from the uniqueness of homology theories satisfying the Eilenberg--Steenrod axioms.

Moreover, by looking at the positive (instead of the negative) gradient flow of a Morse function one can define a Morse cohomology, and several constructions in algebraic topology can be read within the Morse homology and cohomology setting. This is the case, for instance, 
for Poincar\'e duality, the K\"unneth formula, 
cup and cap products, Thom isomorphisms, umkehr maps and Steenrod squares, just to mention a few examples. This is maybe not too surprising, since the homomorphism defined in~\eqref{bdry} is the boundary operator of a cellular filtration, and all these algebraic topological constructions can be read in cellular homology. But what we will see in the next sections is that a lot of this can be extended also to an infinite dimensional setting, in which a formula analogous to~\eqref{bdry} will define a boundary operator which in general does not come from a cellular filtration, or any other purely topological theory.

\subsection*{Historical notes and bibliography}
Morse created his theory around 1930, see~\cite{Mor32}.      
The bible for Morse theory is Milnor's book~\cite{Mi63} 
and other textbooks are~\cite{Ma02,Ni11,TaMePi01}.
The idea of Morse homology is essentially contained in Thom's note~\cite{Th49} and in Milnor's book 
on the $h$-cobordism theorem~\cite{Mil65}, and it was then rediscovered several times, 
by Smale~\cite{Sm67}, Witten~\cite{Wi82}, and Floer~\cite{Flo88:Lag}.
A vivid account of this history is given in Bott's~\cite{Bo89}. A concise exposition of Morse homology can be found in~\cite{Web06}. Textbooks on Morse homology are~\cite{AuDa14,BaHu04}, 
and~\cite{Hu02,Schw93} present Morse homology with a perspective to Floer homology. 
Already in the sixties, Palais~\cite{Pal63} and Smale~\cite{Sma64a,Sma64b} extended Morse theory to functionals on infinite dimensional Hilbert manifolds having critical points with finite Morse index. 
In this infinite dimensional theory, the compactness of the domain is replaced by a compactness assumption on the functional which is now known as the Palais--Smale condition. A textbook on this version of infinite dimensional Morse theory is~\cite{Cha93}. The Morse homology approach works nicely also in this setting, 
see~\cite{AbbMaj06}.

\section{Basics on Hamiltonian dynamics and symplectic geometry, I} 
\label{s:basicsI}

Before discussing the first and easiest instance of Floer homology, we need to describe some notions from Hamiltonian dynamics and symplectic geometry that are used in the subsequent sections.

Consider a particle of unit mass moving in~$\RR^n$, 
subject to a potential force $-\nabla V_t(q)$ that may depend on time.
Here, $n$ may be large, since by ``a particle'' we may mean ``$k$ particles in the plane~$\RR^2$
or in space $\RR^3$'', and then $n=2k$ or $n=3k$.
According to Newton's law, the evolution curve $q(t)$ of the particle 
satisfies the second order ordinary differential equation on~$\RR^n$
\[
\ddot q(t) = -\nabla V_t(q(t)) .
\]
Rewrite this equation as the first order differential equation for the curve~$(q(t),p(t))$
on~$\RR^{2n}$:
\[
\left\{\begin{array} {lcl}
\dot q(t)      &= &   p(t), \\ [0.2em]
\dot p(t)      &= &   -\nabla V_t(q(t)) .
\end{array}\right.
\]
Introducing the function $H_t(q,p) = \frac 12 |p|^2 + V_t(q)$ on $\RR \times \RR^{2n}$,
that represents the total energy, 
this system, in turn, becomes
\begin{equation} \label{e:Ham}
\left\{\begin{array} {lcl}
\dot q(t)      &= &   \phantom{-} \displaystyle \frac{\pp H_t}{\pp p}(q(t), p(t)), \\ [1 em]
\dot p(t)      &= &   -           \displaystyle \frac{\pp H_t}{\pp q}(q(t), p(t)) .
\end{array}\right.
\end{equation}
The beautiful skew-symmetric form of this system leads to a geometric reformulation:
Define the constant differential 2-form $\go_0$ on $\RR^{2n}$ by
\[
\go_0 \,=\, \sum_{i=1}^n dp_i \wedge dq_i .
\]
This 2-form is non-degenerate in the sense that $\go_0(u,v) = 0$ for all $v \in \RR^{2n}$ implies $u=0$.
Hence with $z = (q,p) \in \RR^{2n}$, the equation
\[
\go_0(X_{H_t}(z), \cdot) \,=\, -dH_t (z)
\]
defines a unique time-dependent vector field $X_{H_t}$ on $\RR^{2n}$,
and one easily sees that 
\[
X_{H_t} = \left( \frac{\pp H_t}{\pp p} , - \frac{\pp H_t}{\pp q}\right).
\]
Therefore, the flow of $X_{H_t}$ yields the solution curves of~\eqref{e:Ham}.

Note that 
\[
X_{H_t} = -J_0 \1 \nabla H_t,
\] 
where $J_0 \cong i$ is the standard complex structure on $\RR^{2n} \cong \CC^n$, $(q,p) \mapsto q+ip$ .
This suggests that complex structures may be relevant to Hamiltonian dynamics.

The above geometrisation of Hamilton's equations~\eqref{e:Ham} has the advantage that it readily generalizes to manifolds:
A symplectic manifold is a smooth manifold $M$ endowed with a non-degenerate closed differential 2-form~$\go$.
Again, non-degenerate means that $\omega_x(u,v)=0$ for all $v \in T_xM$ implies $u=0$, and closed means
that the exterior derivative vanishes, $d \go =0$.
Given a smooth function $H \colon \RR \times M \to \RR$ (the Hamiltonian function), 
the vector field $X_{H_t}$ is defined by $\go (X_{H_t}, \cdot ) = -dH_t( \cdot )$,
and its flow $\phi_H^t$ is called the Hamiltonian flow of~$H$.

\begin{examples*}
{\rm 
{\bf 1.\ Cotangent bundles.}\
At least historically, the most important examples of symplectic manifolds are cotangent bundles.
Let $Q$ be a smooth manifold and $T^*Q = \bigcup_q T_q^*Q$ its cotangent bundle. 
Given local coordinates $q = (q_1, \dots, q_n)$ on~$Q$,
an element $\ga \in T_q^*Q$ of this manifold can be written as $\ga = \sum_i p_i \,dq_i$ for suitable real numbers~$p_i$. 
This defines local coordinates $(q,p)$ on~$T^*Q$. 
The {\bf Liouville form} $\lambda_{\can} := \sum_i p_i \,dq_i$ on~$T^*Q$ does not depend on the choice of these coordinates.
Its derivate $\go_{\can} =: d \lambda_{\can}$ is the canonical symplectic form on~$T^*Q$. 
It agrees with $\go_0$ if $Q = \RR^n$.

Fix a Riemannian metric $g$ on~$Q$, a function $V \colon \RR \times Q \to \RR$
and a 1-form~$\ga$ on~$Q$. 
Then the Hamiltonian system on $(T^*Q,\go_{\can})$ of the function
\begin{equation}  \label{e:classical}
H_t(q,p) \,=\, \tfrac 12\, \| p-\ga (q) \|^2 + V_t(q)
\end{equation}
describes the dynamics of a particle subject to the potential force $-\nabla_g V_t$ (such as a gravitational or electrostatic force) and to the force of the vector potential~$A$ that is $g$-dual to~$\ga$ (such as a magnetic force).
Here, $\| \cdot\|$ is the norm on~$T^* Q$ which is dual to the norm on $TQ$ defined by the metric $g$.
If $V$ and~$\ga$ vanish, this system is the (co-)geodesic flow of the Riemannian manifold~$(Q,g)$.

\m \ni
{\bf 2.}\
Any orientable surface endowed with an area-form is a symplectic manifold.

\m \ni
{\bf 3.}\
The product $(M_1 \times M_2, \go_1 \oplus \go_2)$ of symplectic manifolds $(M_i,\go_i)$ is symplectic. 
The torus $\TT^2 \times \dots \times \TT^2 = \TT^{2n} = \RR^{2n}/\ZZ^{2n}$ with the symplectic structure induced by $\go_0$
is an important example.

\m \ni
{\bf 4.}\
All K\"ahler manifolds are symplectic manifolds, but there are many smooth manifolds
that admit a symplectic form but no K\"ahler form. 
\diam 
}
\end{examples*}

The natural transformations of symplectic manifolds are {\bf symplectomorphisms}, that is, smooth diffeomorphisms 
$\varphi \colon M_1 \rightarrow M_2$ such that
\[
\varphi^* \omega_2 = \omega_1
\]
where $\omega_1$ and $\omega_2$ are the symplectic structures on $M_1$ and~$M_2$, respectively. 
On the one hand, symplectomorphisms preserve the structure of Hamilton equations: 
Any symplectomorphism $\varphi \colon (M_1,\omega_1) \rightarrow (M_2,\omega_2)$ maps the orbits of the Hamiltonian system on~$M_1$ which is induced by the Hamiltonian $H \colon \RR \times M_1 \rightarrow \RR$ to the orbits of the Hamiltonian system on~$M_2$ which is induced by the Hamiltonian 
$H \circ (\mathrm{id} \times \varphi^{-1}) \colon \RR \times M_2 \rightarrow \RR$. 
On the other hand, the flow of any (possibly time-dependent) Hamiltonian vector field~$X_{H_t}$ 
consists of symplectomorphisms.

The symplectomorphisms $\varphi \colon (M,\omega) \rightarrow (M,\omega)$ which can be obtained as time-one maps of a time-dependent Hamiltonian vector field are called 
{\bf Hamiltonian diffeomorphisms}. 
Not all symplectomorphisms~$\varphi$ of $(M,\omega)$ are Hamiltonian: 
An obvious necessary condition is that $\varphi$ should be isotopic to the identity map
through a path $\varphi^t$, $0 \leqslant t \leqslant 1$, of symplectomorphisms.
In this case, $\varphi$ is Hamiltonian if and only if one can find such an isotopy whose 
{\bf flux} vanishes:
The integral of $\omega$ over the ``cylinder'' $\bigcup_{0 \leqslant t \leqslant 1} \varphi^t(C)$
vanishes for every closed curve~$C$ in~$M$.
This integral depends only on the homology class of~$C$.
On simply connected manifolds there is therefore no difference between Hamiltonian diffeomorphisms and symplectic diffeomorphisms that are isotopic to the identity 
through symplectomorphisms.

\begin{figure}[h]   
 \begin{center} 
  \psfrag{a}{$-$}
  \psfrag{+}{$+$}
  \psfrag{0}{$0$}
  \psfrag{1}{$1$}
  \psfrag{q}{$q$}
  \psfrag{p}{$p$}
  \psfrag{fS}{$\phi (S^1)$}
  \leavevmode\includegraphics{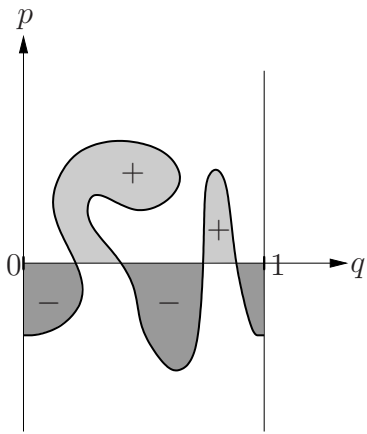}
 \end{center}
 \caption{$\phi$ is Hamiltonian if and only if $+ = -$}
 \label{fig.flux}
\end{figure}

As an example consider the cylinder $T^*S^1$. 
For $\tau \in \RR$ the translation $\psi^\tau (q,p) = (q,p+\tau)$ is symplectic, 
but Hamiltonian only if $\tau =0$.
For a Hamiltonian diffeomorphism~$\phi_H$ of $T^*S^1$ the area of the region 
between $\phi_H(S^1)$ and~$S^1$ that lies over~$S^1$ must be equal to the area of the region under~$S^1$, 
see Figure~\ref{fig.flux}.
Hence for a Hamiltonian diffeomorphism, $\phi_H (S^1)$ and $S^1$ intersect in at least
two points.

%

For historical and bibliographical notes to this section see the end of Section~\ref{s:basicsII}. 


\section{Hamiltonian Floer homology} \label{s:Ham}

\subsection{The Hamiltonian Arnol'd conjecture} \label{s:hamiltonian}
Consider a Hamiltonian function $H \colon \RR \times M \to \RR$ on a symplectic manifold~$(M,\go)$ with flow~$\phi_H^t$.
We would like to find ``points that are back at a prescribed time~$T$''.
Without loss of generality we assume that~$T=1$.
We thus look for $x \in M$ such that 
\[
\phi_H^1(x) = x .
\]
With the language introduced in the previous section: We are looking for fixed points of a Hamiltonian diffeomorphism of~$M$. We first observe that we may assume that the Hamiltonian~$H_t$ 
is 1-periodic in time: Choose a  smooth function~$\chi$ on~$[0,1]$ such that $\chi =0$ near~$0$ and $\chi =1$ near~$1$,
and set $H_t^\chi(x) = \chi'(t) \, H_{\chi (t)}(x)$. 
Then $\phi_{H^{\chi}}^t = \phi_{H}^{\chi(t)}$ for $t \in [0,1]$, 
and in particular $\phi_{H^{\chi}}^1 = \phi_{H}^1$. 
Since $H_t^\chi$ vanishes for $t$ near~$0$ and~$1$, 
we can extend $H_t^\chi$ to a function on $\TT \times M$, where $\TT = \RR/\ZZ$.

Therefore, looking for fixed points of a Hamiltonian diffeomorphism is equivalent to looking for 1-periodic orbits of a Hamiltonian system defined by a Hamiltonian $H \colon \TT \times M \to \RR$. 
The search for periodic orbits is a very traditional topic in dynamical systems. 
After stationary points, periodic orbits are the simplest invariant sets. 
In Hamiltonian systems, that grew out of celestial mechanics, periodic orbits are particularly 
important, 
since one would like to know ``whether the moon comes back''.

On a non-compact symplectic manifold, Hamiltonian diffeomorphisms may have no fixed points. 
For instance, the flow of $H(q,p) = p$ on~$\RR^2$ moves all points to the right. 
We thus assume that $M$ is closed. 
If $H$ does not depend on time, then the critical points of~$H$ are stationary points, 
since there $X_H$ vanishes. If $H$ is also Morse, Theorem~\ref{t:Morse} already guarantees that 
the number of fixed points of~$\phi_H$ is at least the sum of the Betti numbers of~$M$. 
In the sixties, Arnol'd conjectured that the same lower bound should hold for every Hamiltonian diffeomorphism of a closed symplectic manifold with only non-degenerate fixed points:

\begin{Arnold.Ham} 
\label{arnold}
The number of fixed points of a non-degenerate Hamiltonian diffeomorphism 
of a closed symplectic manifold~ $(M,\omega)$ is at least the sum of the Betti numbers of~$M$.
\end{Arnold.Ham}

Here, a fixed point~$x$ of~$\phi_H$ is said to be {\it non-degenerate}\/ 
if the differential of~$\phi_H$ at~$x$ does not have 1 as an eigenvalue. 
By the inverse mapping theorem, non-degenerate fixed points are isolated. 
And a Hamiltonian diffeomorphism is said to be {\it non-degenerate}\/ 
if all its fixed points are non-degenerate.
Equivalently, its graph in $M \times M$ is transverse to the diagonal, 
that is, to the graph of the identity mapping.

Actually, Arnol'd's original formulation of Conjecture~\ref{arnold} is stronger and gives as lower bound the minimal number of critical points of a Morse function on~$M$, but the homological version stated above is commonly accepted as a more treatable statement. 
Notice that the sum of the Betti numbers is usually much larger than the lower bound given by the Euler characteristic of~$M$, which would hold for every diffeomorphism of~$M$ isotopic to the identity 
thanks to the Lefschetz fixed point theorem. For instance, a translation on the $2n$-torus $\TT^{2n}$ has no fixed points, in accordance to the fact that the Euler characteristic of~$\TT^{2n}$ is zero, but Arnol'd's conjecture requires a non-degenerate Hamiltonian diffeomorphism on~$\TT^{2n}$ to have at least~$2^{2n}$ fixed points.

The first instance of Floer homology which we describe in this section was born in order to give a positive answer to this conjecture.

\m \ni
{\bf The space and the function.}
The space of candidates for $1$-periodic orbits of~$\phi_H^t$ is the space $C^\infty(\TT,M)$
of smooth loops in~$M$. 
We look for a function on this space whose critical points are the 1-periodic solutions 
of the Hamilton equation. For $(M,\omega) = (\RR^{2n},\omega_0)$, the action functional 
of classical mechanics 
\[
\ca_H (x) \,=\, \int_{\TT} \bigl(p(t)\cdot \dot q(t) - H_t(x(t)) \bigr) \,dt, \qquad x(t) = (q(t), p(t)),
\]
does the job.
Similarly, if $\go$ has a primitive $1$-form, $\omega = d \lambda$, we can take the function
\begin{equation}
\label{withprimitive}
\ca_H (x) = \int_{\TT} \bigl( \lambda(\dot x(t)) - H_t (x(t)) \bigr) \,dt.
\end{equation}
But on a closed symplectic manifold, $\go$ is certainly not exact,
because the non-degeneracy of~$\go$ implies that $\go^n$ is a volume form on~$M$, 
whence $[\go]^n \neq 0$ in the de Rham cohomology $\H^{2n}(M;\RR)$, and so~$[\go] \neq 0$ in~$\H^2(M;\RR)$.

Notice that the first term in the latter function is the integral of the 1-form~$\lambda$ over the closed curve~$x$. When this closed curve bounds a disc, this integral coincides with the integral of~$\omega$ on this disc. This suggests that if we restrict the attention to the space of {\it contractible}\/ 
loops $C^\infty_{\mathrm{contr}}(\TT,M)$, we still have a well-defined action functional, 
in which we replace the integral of~$\lambda$ over~$x$ by the integral of~$\omega$ over~$\tilde{x}(\DD)$, where $\tilde{x} \colon \DD \rightarrow M$ is a smooth map on the closed unit disc $\DD\subset \CC$ 
such that $\tilde{x}(e^{2\pi i t}) = x(t)$. 
In order to avoid integrating over the possibly bad set $\tilde x (\DD)$, a good way of writing the action functional is then
\begin{equation} \label{e:actionfu}
\ca_H (x) \,=\int_{\DD} \tilde{x}^* \omega - \int_{\TT} H_t(x(t)) \,dt, \qquad x \in C^\infty_{\mathrm{contr}}(\TT,M).
\end{equation}
For this functional to be well-defined, we need that the integral of $\tilde{x}^*\omega$ does not depend on the choice of the extension to the disc~$\tilde{x}$. Two different extensions can be patched together to form a map $u \colon S^2 \rightarrow M$, and what we need is that the integral of~$\omega$ on it vanishes. 
This is guaranteed by, and is actually equivalent to, the following assumption on 
the symplectic manifold~$(M,\omega)$:

\begin{assumption}
\label{aspherical}
The cohomology class $[\go]$ vanishes on the image
of the second homotopy group~$\pi_2(M)$ in $\H_2(M;\ZZ)$, 
or $[\go] |_{\pi_2(M)} =0$ for short.
\end{assumption}

This assumption holds for the torus $\TT^{2n}$, for all orientable surfaces of genus at least one, and hence for products of these symplectic manifolds. It does not hold for the 2-sphere and for the complex projective space. Throughout this section, we will make this assumption.

The tangent space to $C^\infty_{\mathrm{contr}}(\TT,M)$ at a loop~$x$ can be thought of as 
the space of tangent vector fields along $x$: An element of $T_x C^\infty_{\mathrm{contr}}(\TT,M)$ is a smooth map $\xi \colon \TT \rightarrow TM$ such that $\xi (t)\in T_{x(t)} M$ for all $t \in \TT$. 
The differential of~$\ca_H$ at~$x$ can be computed to be
\begin{equation}
\label{diff}
d\ca_H(x)[\xi] \,=\, - \int_{\TT} \omega( \dot{x} - X_{H_t}(x),\xi )\, dt,
\end{equation}
and the non-degeneracy of $\omega$ shows that the critical points of~$\ca_H$ are precisely the 1-periodic orbits of~$X_{H_t}$. It is not difficult to check that the elements in the kernel of the second differential of~$\ca_H$ at a critical point~$x$ are the 1-periodic vector fields~$\xi$ along~$x$ which satisfy the linearisation of the Hamilton equation, that is, the vector fields of the form
\[
\xi(t) := d\phi_H^t(x(0)) \,\xi_0,
\]
where $\xi_0 \in T_{x(0)}M$ is an eigenvector of $d\phi_H^1(x(0))$ with eigenvalue~1. 
Therefore, the fact that all fixed points of~$\phi_H^1$ are non-degenerate 
is equivalent to the fact that the Hamiltonian action functional is Morse.

In order to understand the analytical properties of the functional~$\ca_H$, 
it is convenient to rewrite it in the case $(M,\omega) = (\RR^{2n},\omega_0)$: 
In this case, $\omega_0$ is the differential of the 1-form 
$\lambda_0(x)[v] = \frac 12  \omega_0(x,v) = -\frac 12  J_0 x\cdot v$ and we obtain 
\[
\ca_H (x) \,=\, -\frac{1}{2} \int_{\TT} J_0 x\cdot \dot{x} \, dt - \int_{\TT} H_t(x(t)) \,dt.
\]
The identification $(\R^{2n},J_0) \cong (\CC^n,i)$ gives us the following Fourier series expansion 
for an arbitrary smooth loop $x \colon \TT \rightarrow \R^{2n}$,
\[
x(t) = \sum_{k\in \ZZ} e^{2\pi k t J_0} \hat{x}_k, \qquad \hat{x}_k\in \R^{2n},
\]
and $\ca_H(x)$ takes the form
\begin{equation}
\label{fourier}
\ca_H (x) \,=\, -\pi \sum_{k\in \ZZ} k |\hat{x}_k|^2 - \int_{\TT} H_t(x(t)) \,dt.
\end{equation}
Therefore, the leading term of this functional is a quadratic form which is 
positive definite on the space of loops whose Fourier expansion involves only 
coefficients~$\hat{x}_k$ with $k<0$ and is negative definite on the space of loops 
whose Fourier expansion involves only coefficients~$\hat{x}_k$ with~$k>0$. 
Both of these spaces are infinite dimensional. The second differential of~$\ca_H$ 
at~$x$ has the form
\begin{eqnarray*}
d^2\ca_H (x)[\xi,\xi] 
&=& 
-2\pi \sum_{k\in \ZZ} k |\hat{\xi}_k|^2 - \int_{\TT} d^2H_t(x)[\xi,\xi] \,dt \\
&=& 
- \int_{\TT} J_0 \xi \cdot \dot \xi \,dt - \int_{\TT} d^2H_t(x)[\xi,\xi] \,dt .
\end{eqnarray*}
The second integral involves no derivatives of~$\xi$, hence is a lower order term.
It follows that $d^2\ca_H (x)$ is also positive definite on an infinite dimensional space 
and negative definite on an infinite dimensional space. 
Therefore, all critical points of~$\ca_H$ have infinite Morse index and co-index, 
so the classical extension of Morse theory to infinite dimensional Hilbert manifolds 
does not give any information: passing a critical value would change the sublevels 
of~$\ca_H$ by the attachment of an infinite dimensional ball along its boundary, 
and since infinite dimensional Hilbert balls are retractable onto their boundary, 
the homotopy type of the sublevels never changes.

It is a general fact in nonlinear analysis that the ``right'' functional space 
for studying a functional by topological methods is the space with the weakest 
topology on which this functional is ``regular'', for instance continuously differentiable. 
Indeed, regularity is needed to have a negative gradient equation giving 
a well-posed~ODE on the infinite dimensional space;
the choice of a finer topology keeps the regularity 
(and sometimes even improves it), 
but usually causes the failure of an important compactness condition known as the Palais--Smale condition. The leading term in~\eqref{fourier} shows that
in the case of $(\RR^{2n}, \omega_0)$ or of the torus $(\TT^{2n},\omega_0)$, in which the action functional can again be expressed in terms of the Fourier coefficients of the loops and has the same form, 
the ``right'' functional space is the fractional Sobolev space~$H^{\frac 12}$. Indeed, the quantity
\[
|\hat{x}_0|^2 + \sum_{k\in \ZZ} |k| |\hat{x}_k|^2
\]
is the square of the $H^{\frac 12}$ norm of the map $x \colon \TT \rightarrow \R^{2n}$. 
It is by using this function space, together with suitable finite dimensional reductions 
in order to overcome the problem of the infiniteness of the Morse indices, 
that Conley and Zehnder proved the Arnol'd conjecture for the $2n$-torus in~1983, see~\cite{CoZe83}.

Besides for the infiniteness of the Morse indices, the study of the action functional on loops 
in a more general symplectic manifold has the extra difficulty that loops of Sobolev class $H^{\frac 12}$ 
taking values in a manifold do not form a Hilbert manifold. 
Indeed, in general the space of Sobolev maps between manifolds has a good differentiable structure 
only when the Sobolev space embeds into the space of continuous loops, 
and the Sobolev space $H^{\frac 12}(\TT)$ does not have this property 
(the exponent~$\frac 12$ is critical for the embedding into the space of continuous functions).

\subsection{Floer's idea}  \label{s:idea}
In the mid eighties, Floer came up with a new idea which solves both difficulties at one stroke, 
and which we now wish to describe. The starting point is to replace the $H^{\frac 12}$-gradient of~$\ca_H$ 
by the $L^2$-gradient. Since the functional is not even well-defined in~$L^2$, 
this seems to be a hazardous move. However, if the $L^2$ metric is chosen properly, the resulting negative gradient equation turns out to be very nice. 

In order to describe this metric, we need to introduce some further structure on~$(M,\omega)$. 
An {\em almost complex structure}~$J$ on~$M$ is a smooth collection of linear maps 
$J_x \colon T_xM \to T_xM$ such that $J_x^2 = -\id$. Such a $J$ is said to be compatible with~$\go$ if 
\begin{equation} \label{e:g}
g_J(\xi,\eta) \,=\, \omega (J\xi, \eta)
\end{equation}
defines a Riemannian metric on $M$. The standard complex structure~$J_0$ on~$\RR^{2n}$ is compatible 
with~$\omega_0$, and $g_{J_0}$ is the Euclidean inner product. 
A symplectic manifold $(M,\omega)$ always admits $\omega$-compatible almost complex structures. 
These form a contractible space.

\begin{remark}
{\rm 
We warn the reader that the usual convention for $\omega$-compatibility has the opposite sign: 
$g_J(\xi,\eta) = \omega (\xi, J\eta)$. In the literature about Floer homology and the theories 
originated from it, sign conventions differ from paper to paper. 
The reason is that in this field three different sets of conventions come together: 
The convention that in Morse theory, and more generally in the calculus of variations, 
one preferably works with {\em negative}\/ gradient flows; 
the sign conventions from classical mechanics, which require the symplectic form to be 
$dp\wedge dq$, and not the opposite, and the Hamiltonian to appear with a negative sign in the Hamiltonian action; 
the sign conventions from complex geometry, which require the Riemannian structure~$g$ and the complex structure~$J$ on a K\"ahler manifold to be related by the identity $g(\xi,\eta)=\omega(\xi,J\eta)$. 
If one adopts all these conventions, the Floer equation, which we derive below, will have as 
leading term the anti-Cauchy--Riemann operator. The fact that holomorphic functions are more 
natural than anti-holomorphic ones makes it preferable to have the Cauchy--Riemann operator 
instead. In order to have this, at least one of the above sign conventions has to be changed. 
Our choice here is to sacrifice the sign convention from complex geometry. 
Notice that with our convention the Hamiltonian vector field has the form $X_H = - J \nabla H$, 
where $\nabla H$ is the gradient of~$H$ with respect to the metric~$g_J$ defined in~\eqref{e:g}.
}
\end{remark}

Given a loop $x \in C^{\infty}_{\rm contr}(\TT,M)$ and two smooth vector fields $\xi, \eta$ along~$x$,
define their $L^2$-inner product by
\[
\langle \xi, \eta \rangle_{L^2} \,=\, \int_{\TT} g_J (\xi(t), \eta(t)) \,dt. 
\]
By \eqref{diff}, the differential of $\ca_H$ at~$x$ can be rewritten as
\[
d\ca_H(x)[\xi] \,=\, \int_{\TT} g_J(J(\dot{x} - X_{H_t}(x)),\xi)\, dt \,=\, \langle J(\dot{x} - X_{H_t}(x)),\xi \rangle_{L^2},
\]
so the $L^2$-gradient of $\ca_H$ at $x$ is the following vector field along~$x$:
\[
\nabla_{L^2} \ca_H(x) = J(\dot{x} - X_{H_t}(x)).
\]
Then the negative $L^2$-gradient equation for a curve $u \colon \R \to C^{\infty}_{\rm contr}(\TT,M)$,
\[
\frac{d}{ds} u(s) = - \nabla_{L^2} \ca_H(u(s))
\]
becomes the following PDE, 
\begin{equation}
\label{floer}
\frac{\partial u}{\partial s} + J(u) \left( \frac{\partial u}{\partial t} - X_{H_t}(u) \right) = 0,
\end{equation}
once we view $u$ as a map from the cylinder $\R \times \TT$ into $M$.

For $M=\RR^{2n}$, $J = J_0$ and $H=0$, this is the Cauchy--Riemann equation 
determining holomorphic mappings from the cylinder $\RR \times \TT$ to~$\CC^n$. 
On an arbitrary $(M,\omega,J)$, this is a zero-order perturbation of the equation for pseudo-holomorphic curves that Gromov introduced in~1985, see~\cite{Gro85}. 

This equation does not define a flow on $C^{\infty}_{\rm contr}(\TT,M)$, 
nor on any other reasonable function space. However, as Floer first realized,
this is not much of a problem if we want to use it as a negative gradient equation for defining a boundary operator using a formula analogous to~\eqref{bdry}. Indeed, \eqref{bdry} involves only gradient flow lines which are asymptotic to two critical points. Therefore, what we need is that the space of solutions~$u$ 
of~\eqref{floer} which for $s \rightarrow \pm \infty$ converge to two loops that are critical points 
of~$\ca_H$ is a finite dimensional manifold with good compactness properties.

We now review the steps which lead to the definition of the Floer homology groups for the action 
functional~$\ca_H$. The assumptions for now are that $(M,\omega)$ is a closed symplectic manifold 
with $[\omega]|_{\pi_2(M)}=0$ and $H \colon \TT \times M \rightarrow \RR$ is a smooth Hamiltonian all of whose 1-periodic orbits are non-degenerate 
(that is, all fixed points of $\phi^1_H$ are non-degenerate).

Fix two critical points $x$ and $y$ of~$\ca_H$. Denote by $\widehat{\cm}_{H,J}(x,y)$ the set of smooth solutions $u \colon \RR \times \TT \to M$ of the Floer equation~\eqref{floer} which satisfy the 
asymptotic conditions
\[
\lim_{s\rightarrow -\infty} u(s,\cdot) = x \qquad \mbox{and} \qquad \lim_{s\rightarrow +\infty} u(s,\cdot) = y \qquad \mbox{in } C^{\infty}(\TT,M).
\]
The set $\widehat{\cm}_{H,J}(x,y)$ replaces the intersection of unstable and stable manifolds of two critical points in Morse homology. 

\subsection{The Fredholm index of the linearized operator}  \label{s:fredholm}
Let $u$ be an element of $\widehat{\cm}_{H,J}(x,y)$. 
By choosing a unitary trivialisation of the tangent bundle of~$M$ along~$u$ and linearizing the Floer equation~\eqref{floer} along~$u$ 
one gets an operator of the form
\begin{equation}
\label{operator}
v \mapsto \frac{\partial v}{\partial s} + J_0 \frac{\partial v}{\partial t} - A(s,t)\, v
\end{equation}
for maps $v \colon \RR \times \TT \rightarrow \RR^{2n}$, where $A \colon [-\infty,+\infty] \times \TT \rightarrow L(\R^{2n},\R^{2n})$ is a continuous map such that the asymptotic systems
\begin{equation}
\label{aslin}
\dot{w}(t) = J_0 A(-\infty,t) \,w(t) \qquad \mbox{and} \qquad \dot{w}(t) = J_0 A(+\infty,t) \,w(t)
\end{equation}
are the linearisations of the Hamilton equation along $x$ and~$y$, respectively. 
It is convenient to see~\eqref{operator} as a bounded linear operator from the Sobolev Space 
$W^{1,p}(\RR\times \TT,\RR^{2n})$ to $L^p(\RR\times \TT,\RR^{2n})$ for $p \in (2,\infty)$, 
and the first step is to prove that this operator is Fredholm, that is, 
it has finite dimensional kernel and cokernel. This fact follows from the ellipticity of the Cauchy--Riemann operator $\frac{\partial}{\partial s} + J_0 \frac \partial{\partial t}$ 
and from the fact that the linear systems~\eqref{aslin} do not have non-trivial 1-periodic solutions, 
due to the non-degeneracy assumption. 
In order to discuss the Fredholm index of the operator~\eqref{operator}, that is, the difference of the dimension of its kernel and its cokernel, we need to introduce the Conley--Zehnder index.

The space of linear Hamiltonian systems with 1-periodic coefficients and no non-trivial 1-periodic solutions has countably many connected components. These components are labeled by an integer which is called 
the Conley--Zehnder index after~\cite{CoZe83,CoZe84}, 
but which appeared also in an older paper of Gelfand and Lidski~\cite{GL58}. 
This integer counts the ``windings'' of the fundamental solution of a linear Hamiltonian system, 
which is a path of linear symplectic mappings starting at the identity, 
within the symplectic linear group (which has $\ZZ$ as fundamental group). 
By looking at the spectral flow of the path of operators
\[
s \mapsto -J_0 \frac d{dt} + A(s,t) ,
\] 
one can prove that the Fredholm index of the operator (\ref{operator}) is the Conley--Zehnder index of the first linear system in~\eqref{aslin} minus the Conley--Zehnder index of the second one.

One now would like to define the Conley--Zehnder index of a non-degenerate Hamiltonian periodic orbit~$x$ 
as the integer label attached to the linearisation of the Hamilton equation along~$x$. 
This is well-defined provided that the label of the linearized equation does not depend on the choice of 
the trivialisation of the tangent bundle of~$M$ along~$x$ which produces the linearized equation. 
Since we are working with contractible periodic orbits, we can restrict the attention to trivialisations 
of~$TM$ along~$x$ which extend to a trivialisation of~$TM$ on a disc bounded by~$x$. 
In order for the Conley--Zehnder index associated to these trivialisations to be well-defined 
we need one more assumption on the symplectic manifold~$(M,\omega)$:

\begin{assumption}
The first Chern class $c_1(\omega) \in \H^2(M)$ of $(M,\omega)$ vanishes on spheres: $c_1(\omega)|_{\pi_2(M)}=0$.
\end{assumption} 

Here $c_1(\omega)$ is the first Chern class of the complex vector bundle~$(TM,J)$, 
where $J$ is an $\omega$-compatible almost complex structure on~$M$. 
The fact that the space of such almost complex structures is contractible implies that this characteristic class depends only on the symplectic form~$\omega$ and justifies the notation. 
This assumption, together with Assumption~\ref{aspherical}, is usually summarized by saying that 
$(M,\omega)$ is {\em symplectically aspherical}. Again, the $2n$-torus, surfaces of genus at least~1 and products thereof are symplectically aspherical.

When $(M,\omega)$ is symplectically aspherical, the Conley--Zehnder index of every non-dege\-nerate periodic Hamiltonian orbit~$x$ is a well-defined integer, which we denote by~$\mu_{\sCZ}(x)$. 
By the above discussion, the Fredholm index of the operator~\eqref{operator} 
which is obtained by linearizing the Floer equation~\eqref{floer} along some solution~$u$ 
in $\widehat{\cm}_{H,J}(x,y)$, is the integer
\[
\mu_{\sCZ}(x) - \mu_{\sCZ}(y).
\]

\paragraph{\bf Transversality.} 
When we can guarantee that for every $u \in \widehat{\cm}_{H,J}(x,y)$ the operator~\eqref{operator}, 
which is obtained by linearizing the Floer equation~\eqref{floer} along~$u$, is surjective, 
the implicit function theorem implies that $\widehat{\cm}_{H,J}(x,y)$ is a smooth manifold 
whose dimension equals the dimension of the kernel of this operator, which by the above discussion 
and the fact that now the cokernel is zero is $\mu_{\sCZ}(x) - \mu_{\sCZ}(y)$:
\[
\dim \widehat{\cm}_{H,J}(x,y) \,=\, \mu_{\sCZ}(x) - \mu_{\sCZ}(y).
\]
This shows that the Conley--Zehnder index in Floer theory plays the role of the Morse index in Morse homology. The surjectivity of the operator~\eqref{operator} is not to be expected in general, 
but it can be achieved by a generic choice of the $\omega$-compatible almost complex structure~$J$, 
in analogy with the Morse homology case, in which the Morse--Smale condition was achieved 
by choosing a generic metric~$g$. Here actually one has to work with almost complex structures~$J$ 
which depend smoothly and periodically on time. Replacing the $t$-independent~$J$ by a $t$-dependent one 
does not make the Floer equation~\eqref{floer} more complicated, and one can prove that for a generic choice of such a~$J$ the linearized operator at any $u \in \widehat{\cm}_{H,J}(x,y)$ and for any pair of periodic orbits~$x$ and~$y$ is surjective.

The proof of this fact uses standard techniques in transversality theory, but needs also to face 
a difficulty which is not present in the proof of the genericity of the Morse--Smale condition for 
gradient flows. Indeed, here the space of admissible perturbations is quite small: 
If we wish to modify a given solution~$u$ of the Floer equation near a point $(s,t) \in \RR \times \TT$, 
what we can do is modify $J$ near $(t,u(s,t))$. But this change will affect the solution~$u$ also at every point $(s',t) \in \RR \times \TT$ such that $u(s',t)=u(s,t)$. 
This difficulty can be solved thanks to two useful properties of solutions of the Floer equation: 
On the one hand, the complement of the set of points $(s,t) \in \RR \times \TT$ such that $u(s',t)=u(s,t)$ for some $s'\neq s$ has non-empty interior. On the other hand, the solutions of~\eqref{floer}, 
although in general not analytic, satisfy a unique continuation property which is analogous to that of analytic functions. 

\subsection{Compactness} 
\begin{figure}[h]   
 \begin{center}
  \psfrag{x}{$x$}
  \psfrag{x'}{$x'$}
  \psfrag{y}{$y$}
  \psfrag{y'}{$y'$}
  \leavevmode\includegraphics{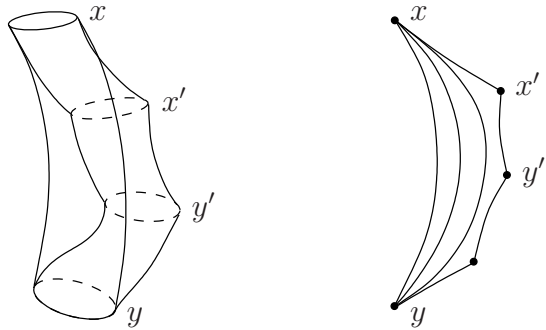}
 \end{center}
 \caption{Breaking, seen in $M$ and in $C^\infty (\TT, M)$}  \label{fig.breaking}
\end{figure}

The starting point for proving compactness results about the space $\widehat{\cm}_{H,J}(x,y)$ is the following identity, which follows from the fact that the Floer equation is a negative gradient equation 
for the action functional:
\[
\int_{\RR} \left\| \frac{\partial u}{\partial s} \right\|_{L^2}^2\, ds \,=\, 
\int_{\RR} \int_{\TT} g_J\left( \frac{\partial u}{\partial s}, \frac{\partial u}{\partial s} \right) \, dt \,ds \,=\, \ca_H(x) - \ca_H(y).
\]
The left-hand side is called the {\em energy}\/ of the solution~$u$, and the above identity tells us that solutions of the Floer equation which are asymptotic to periodic orbits have uniformly bounded energy 
(recall that $\ca_H$ has finitely many critical points).
Starting from this, one can prove that the space $\widehat{\mathcal{\cm}}_{H,J}(x,y)$ is pre-compact 
in~$C^{\infty}_{\mathrm{loc}}$: Any sequence $(u_h) \subset \widehat{\cm}_{H,J}(x,y)$ 
has a subsequence~$(u_h')$ such that for every~$S>0$ the restriction of~$u_h'$ to the compact set 
$[-S,S] \times \TT$ converges to some solution~$u$ of the Floer equation, 
uniformly together with all its derivatives. 
Exactly as in Morse homology, the limiting solution~$u$ needs not be in $\widehat{\cm}_{H,J}(x,y)$, 
but will be in $\widehat{\cm}_{H,J}(x',y')$ for some periodic orbits~$x'$ and~$y'$ with
\[
\ca_H(y) \leqslant \ca_H(y') \leqslant \ca_H(x') \leqslant \ca_H(x),
\]
see Figure~\ref{fig.breaking}. 

The first step in the proof of this compactness result is to show that the differential of~$u$ is 
uniformly bounded when $u$ varies in $\widehat{\mathcal{M}}_{H,J}(x,y)$. This follows from the assumption that $[\omega]$ vanishes on~$\pi_2(M)$ and would not be true in general. 
It is proved by contradiction, by considering a sequence of solutions 
$(u_h) \subset \widehat{\mathcal{M}}_{H,J}(x,y)$ and a sequence of points $(z_h)\in \RR \times \TT$ 
such that $\|du_h(z_h)\|$ diverges. Then a blowing-up argument allows to show that the restriction 
of~$u_h$ to a disc of radius~$r_h$ centered at~$z_h$, with $r_h$ a suitable infinitesimal sequence, 
produces in the limit a non-constant $J$-holomorphic sphere, that is, a non-constant map
\[
v \colon (S^2,j) \rightarrow (M,J)
\]
from the sphere $S^2$ with its standard complex structure~$j$ to the almost complex manifold~$(M,J)$ 
such that
$$
dv \circ j = J \circ dv .
$$
This phenomenon is known as {\em bubbling off of $J$-holomorphic spheres}.
By the compatibility of~$J$ with~$\omega$, the integral of $\|dv\|^2$ over~$S^2$ coincides with twice the integral of~$v^* \omega$, but the latter integral vanishes because of 
the assumption $[\omega]|_{\pi_2(M)}=0$. This shows that $v$ must be constant, 
producing the desired contradiction. If one drops the symplectic asphericity assumption, 
$(M,\omega)$ might have non-constant $J$-holomorphic spheres and this argument fails.

\begin{figure}[h]   
 \begin{center}
  \psfrag{u1}{$u_1$}
  \psfrag{u10}{$u_{10}$}
  \psfrag{u100}{$u_{100}$}
  \psfrag{v}{$v$}
  \leavevmode\includegraphics{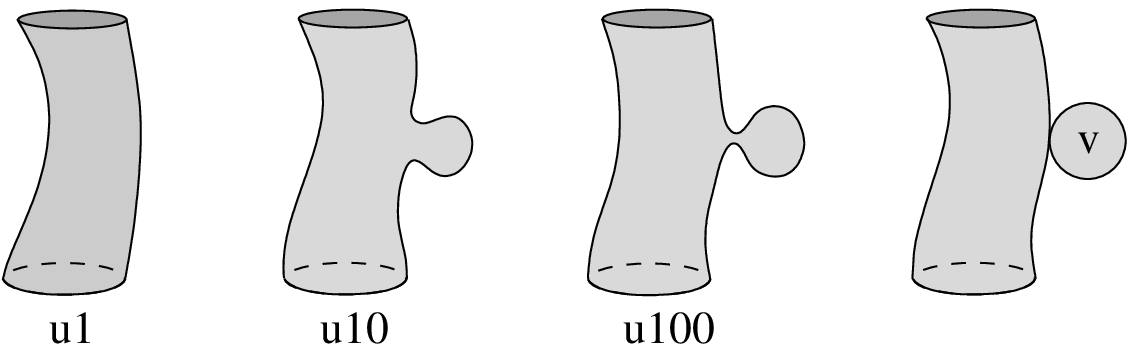}
 \end{center}
 \caption{Bubbling}  \label{fig.bubbling}
\end{figure}

By the theorem of Ascoli--Arzel\`a, the uniform bound on the differential of~$u$ implies uniform convergence on compact subsets of the cylinder, and allows one to localize the analysis by using charts on~$M$ and to consider sequences of $\R^{2n}$-valued maps instead of $M$-valued ones. 
Then a bootstrap argument using the Calderon--Zygmund inequality
\[
\|du\|_{L^p(\RR\times \TT)} \leqslant C_p \left\| \frac{\partial u}{\partial s}  + J_0 \frac{\partial u}{\partial t} \right\|_{L^p(\RR\times \TT)} \qquad 
\forall \, u \in C^{\infty}_c(\RR\times \TT,\RR^{2n}), \; \forall \, p\in (1,\infty),
\]
gives the desired compactness in $C^{\infty}_{\mathrm{loc}}$.
  
\subsection{The Floer complex} 
Now we proceed like in Morse homology and define the analogue of the Morse complex. 
We postpone the discussion about orientations and deal with $\ZZ_2$ coefficients for now. 
For every $k \in \ZZ$, we define $\CF_k(H)$ to be the $\ZZ_2$-vector space generated by the (finitely many) 
1-periodic orbits of Conley--Zehnder index~$k$. When 
\[
\mu_{\sCZ}(x) - \mu_{\sCZ}(y) = 1,
\]
the space $\widehat{\mathcal{M}}_{H,J}(x,y)$ is a one-dimensional manifold which is invariant under the free action of~$\RR$ given by translations of the $s$ variable: 
$(\sigma \cdot u)(s,t) = u(s+\sigma,t)$. 
The quotient $\mathcal{M}_{H,J}(x,y)$ of $\widehat{\mathcal{M}}_{H,J}(x,y)$ by this action is 
a zero-dimensional manifold. Compactness and transversality imply that $\mathcal{M}_{H,J}(x,y)$ is 
a finite set, by an argument which is completely analogous to what we have seen in Morse homology. 
We then define $\nu(x,y) \in \ZZ_2$ to be the parity of this set, and define the boundary homomorphism exactly as in~\eqref{bdry}:
\[
\partial_k \colon \CF_k(H) \to \CF_{k-1}(H), \quad 
\partial_k x = \sum_{\substack{y \in \Crit \ca_H \\ \mu_{\sCZ}(y)=k-1}} \nu(x,y)\, y, \quad 
\forall x \in \Crit \ca_H \mbox{ with } \mu_{\sCZ}(x)=k.
\]
The proof of the identity $\partial_k \circ \partial_{k+1}=0$ is based on the same cobordism argument 
that we have seen in Morse homology. The main point is to show that when $x$, $y$, and $z$ are periodic orbits with
\[
\mu_{\sCZ}(x)=k+1, \qquad \mu_{\sCZ}(y) = k, \qquad \mu_{\sCZ}(z)=k-1,
\]
and $[u] \in \mathcal{M}_{H,J}(x,y)$ and $[v] \in \mathcal{M}_{H,J}(y,z)$, then one can find a unique connected component of $\mathcal{M}_{H,J}(x,z)$ having the pair $([u],[v])$ as a boundary point 
(see again Figure~\ref{fig.gluing}). In the case of a Morse--Smale negative gradient flow on a finite dimensional manifold, this fact can be deduced by the hyperbolic behaviour of the flow near the 
stationary point~$y$. In the case of the Hamiltonian action functional, we do not have a gradient flow 
at our disposal, and the proof is more complicated. It is based on a Newton iteration scheme for finding elements of $\widehat{\mathcal{M}}_{H,J}(x,z)$ starting from the approximated solution which is given by 
``gluing together'' the solutions~$u$ and~$v$. Once this is done, we have a chain complex of finite dimensional $\ZZ_2$-vector spaces $\{\CF_k(H),\partial_k\}$, 
which is called the {\em Floer complex}\/ of~$(H,J)$.

\subsection{Orientations} 
It is sometimes desirable to have a homology theory with integer coefficients, rather than coefficients 
in~$\ZZ_2$. As we have seen, in finite dimensional Morse homo\-lo\-gy this can be done by using orientations of the unstable manifolds. In the infinite dimensional generalisation, we do not have unstable manifolds. 
And even if we could define them somehow, they would be infinite dimensional objects. The strategy in this case is to adopt a common escamotage in global analysis: In infinite dimensional situations, 
orient maps rather than spaces. Indeed, our aim is to orient the finite dimensional manifolds $\widehat{\mathcal{M}}_{H,J}(x,y)$, which means orienting their tangent spaces. 
Each of these tangent spaces is seen as the kernel of the surjective Fredholm operator~\eqref{operator}. 
The space of Fredholm operators between two real Banach spaces is the base space of a real line bundle, called the {\em determinant bundle}, which was introduced by Quillen in~\cite{qui85}, 
whose fiber at the Fredholm operator~$T$ is the space of alternating forms of top degree on the finite dimensional space $\ker T \times (\mathrm{coker}\, T)^*$. In particular, the determinant line at a surjective Fredholm operator is the space of alternating forms of top degree on its kernel, so any non-zero element of this line defines an orientation of the kernel of~$T$. The determinant bundle over the space of all Fredholm operators does not have non-vanishing sections, but its restriction to the space of operators 
of the form~\eqref{operator} with fixed asymptotics~\eqref{aslin} does. 
Starting from these observations one can orient all spaces $\widehat{\mathcal{M}}_{H,J}(x,y)$ 
in a coherent way, and this allows one to define a Floer complex over the integers. 

\subsection{Invariance}  
\label{s:invariance}
As we have seen, a non-degenerate 1-periodic Hamiltonian function on a closed symplectically aspherical 
manifold~$(M,\omega)$ has a Floer complex $\{\CF_k(H),\partial_k\}$, where the boundary operator depends 
on the choice of a generic 1-periodic $\omega$-compatible almost complex structure~$J$. 
Natural questions are now: How does this chain complex vary when we change the almost complex structure? 
And what happens if we change the Hamiltonian? 

If $J_1$ and $J_2$ are two $\omega$-compatible almost complex structures which give us transversality, 
we can join them by a path of $\omega$-compatible almost complex structures and define 
a family~$J(s,t,\cdot)$ of $\omega$-compatible almost complex structures such that 
$J(s,t,\cdot) = J_1(t,\cdot)$ for $s \leqslant 0$ and $J(s,t,\cdot) = J_2(t,\cdot)$ for $s\geqslant 1$. 
Then we can look at the $s$-dependent Floer equation
\[
\frac{\partial u}{\partial s} + J(s,t,u) \left( \frac{\partial u}{\partial t} - X_{H_t}(u) \right) = 0.
\]
When the homotopy $J$ is chosen generically, an algebraic count of the solutions of this equation which are asymptotic to pairs of periodic orbits of~$X_H$ with the same Conley--Zehnder index produces a chain map 
from the Floer complex of~$(H,J_1)$ to that of~$(H,J_2)$. The proof of the chain map property is analogous 
to the proof of~$\partial^2=0$. This map can be shown to be a chain complex isomorphism, 
so the Floer complex changes by an isomorphism when varying the almost complex structure.

Changing the Hamiltonian involves a more drastic change, because in this case also the groups~$\CF_k(H)$ change, and not only the boundary operators. But the same strategy works: 
Given non-degenerate 1-periodic Hamiltonians~$H_1$ and~$H_2$, one interpolates them by a family of Hamiltonians $H(s,t,\cdot)$ with $H(s,t,\cdot) = H_1(t,\cdot)$ for $s \leqslant 0$ and 
$H(s,t,\cdot) = H_2(t,\cdot)$ for $s \geqslant 1$ and uses solutions of the equation
\begin{equation} \label{floer.s}
\frac{\partial u}{\partial s} + J_t(u) \left( \frac{\partial u}{\partial t} - X_{H_{s,t}}(u) \right) = 0
\end{equation}
to build a chain map from the Floer complex of $(H_1,J)$ to that of~$(H_2,J)$:
\[
\Phi_{H_1,H_2} \colon \CF_k(H_1) \rightarrow \CF_k(H_2) .
\]
By considering ``homotopies of homotopies'', one can prove that the chain map $\Phi_{H_2,H_1}$ 
is a homotopy inverse to $\Phi_{H_1,H_2}$, meaning that their compositions in both orders are 
chain homotopic to the identity. This implies that the homomorphisms induced by~$\Phi_{H_1,H_2}$ 
and~$\Phi_{H_2,H_1}$ at the level of homology,
that are called {\it continuation maps},
are inverse to each other. We conclude that the homology 
of the Floer complex of~$(H,J)$ depends on the choice of the data~$(H,J)$ only up to an isomorphism. 
Because of this, one can refer to it as to the {\em Floer homology of $(M,\omega)$}\/ 
and denote it by~$\HF_*(M,\omega)$.

\subsection{Computation of the Floer homology} \label{s:comp}
The last step is to determine the Floer homology of $(M,\omega)$. Notice that for now we have discussed 
the differential structure and the compactness properties of the spaces 
$\widehat{\mathcal{M}}_{H,J}(x,y)$ without 
ever bothering to show that these spaces are not empty. But now that we know that the Floer homology 
does not depend on the choice of the Hamiltonian, we can choose particular ones, for which we can hope 
to say more about the Floer complex. A good idea is to choose an autonomous and $C^2$-small 
Hamiltonian $H \colon M \rightarrow \R$ which is a Morse function. In general, autonomous Hamiltonians are not suited for the version of Floer homology we are considering, because a non-constant 1-periodic orbit cannot be non-degenerate, since its velocity is a periodic solution of the linearized equation. 
But if the autonomous Hamiltonian~$H$ is Morse and $C^2$-small, one can show that the only 1-periodic orbits are the constant ones, which are given by the critical points of~$H$, and these are non-degenerate as periodic orbits because they are non-degenerate as critical points of~$H$. 
Moreover, it is not difficult to compute the Conley--Zehnder index of these constant orbits and to relate 
it to their Morse index:
\[
\mu_{\sCZ}(x) = n-\ind_H(x) \qquad \forall \, x \in \Crit H.
\]
Therefore, we have
\[
\CF_k(H) \,\cong\, \CM_{n-k}(H) \,=\, \CM_{k+n}(-H).
\]
Furthermore, one can show that the $C^2$-smallness of $H$ implies, if we choose~$J$ to be 
$t$-independent, that all the finite energy solutions~$u$ of~\eqref{floer} do not depend on~$t$. 
Therefore, they are of the form $u=u(s)$ and the Floer equation~\eqref{floer} for them reduces to
\[
u'(s) - \nabla H(u(s)) = 0.
\]
This means that the Floer trajectories $u$ are in this case the negative gradient flow lines of~$-H$, 
and hence the Floer complex of~$(H,J)$ is nothing else but the Morse complex of~$(-H,g_J)$, 
after a degree shift of~$n$. From the fact that the Morse homology of~$-H$ is isomorphic to the singular homology of~$M$, we conclude that the same is true for the Floer homology of~$(M,\omega)$:
\begin{equation} \label{iso:HF}
\HF_k(M,\omega) \cong \H_{k+n}(M).
\end{equation}

\subsection{The dictionary}
We summarize the translation of Morse homology to Hamiltonian Floer homology of a closed 
symplectically aspherical manifold $(M,\omega)$ in the following dictionary.

\begin{center}
    \begin{tabular}{ | p{2.4cm} || p{4cm} | p{6.5cm} |}
    \hline
      & Morse homology $\HM_*$ & Floer homology $\HF_*$  \\ \hhline{|=#=|=|}
    manifold & $M$ closed manifold &  $C^\infty_{\mathrm{contr}} (\TT,M)$ \newline space of contractible loops in~$M$
                                            \\ \hline
    function & $f \colon M \to \RR$ \newline Morse function & $\ca_H \colon C^\infty_{\mathrm{contr}} (\TT,M) \to \RR$ \newline action functional of non-degenerate Hamiltonian $H \colon \TT \times M \to \RR$
                                            \\ \hline
    grading & Morse index $\ind_f$ & Conley--Zehnder index $\mu_{\sCZ}$
                                            \\ \hline						
    Riemannian metric & $g$ such that $(f,g)$ is \newline Morse--Smale & generic $\go$-compatible \newline almost complex                                          structure $J_t$    \\ \hline
    gradient flow lines & $\dot u = - \nabla_g f(u)$ \newline negative gradient flow  & $u_s + J_t(u) \left( u_t-X_{H_t}(u)\right) =0$ \newline Floer equation
                                            \\ \hline

    \end{tabular}
\end{center}

\medskip

\subsection{Proof of the Arnol'd conjecture} 
In the case of a closed symplectically aspherical manifold $(M,\omega)$, 
the Arnol'd conjecture~\ref{arnold} follows at once from the fact that the homology of the Floer complex 
of~$(H,J)$, where $H$ is any non-degenerate 1-periodic Hamiltonian, is isomorphic to the singular homology 
of~$M$. Indeed, the total rank of a chain complex is not smaller than the total rank of its homology. Actually, a more precise relationship between the number~$p_k(H)$ of 1-periodic orbits of Conley--Zehnder index $k$ and the Betti numbers~$b_k(M)$ of~$M$ is expressed by the following relation between Laurent polynomials 
\[
\sum_{k\in \ZZ} p_k(H) \, z^k \,=\, \sum_{k=-n}^{n} b_{k+n}(M)z^k + (1+z) Q(z),
\]
where $Q$ is a Laurent polynomial with integer non-negative coefficients.

\subsection{The Conley conjecture}

Let $(M,\go)$ be a symplectically aspherical manifold.
By the proof of the Arnol'd conjecture, 
for a non-degenerate Hamiltonian diffeomorphism $\phi \colon M \rightarrow M$ 
there is a non-trivial, though finite, lower bound for the number of contractible fixed points. 
Here, the adjective ``contractible'' refers to the fact that the fixed points correspond to contractible 1-periodic orbits of the Hamiltonian system whose time-1 flow is~$\phi$. 
Usually weaker but still non-trivial lower bounds hold if one drops the non-degeneracy assumption, see the historical notes below.

We now look for periodic points of $\phi$, that is, fixed points of some iterate~$\phi^k$.
The following result was conjectured by C.\ Conley in~1984
for the special case of the standard torus $\RR^{2n} / \ZZ^{2n}$.

\begin{theorem} \label{t:Conley}
Assume that $(M,\go)$ is closed and symplectically aspherical. 
Then every Hamiltonian diffeomorphism $\phi \colon M \rightarrow M$ has infinitely many 
contractible periodic points.
\end{theorem}

Notice that here no non-degeneracy assumption is required. An equivalent reformulation 
of the above theorem is that any 1-periodic Hamiltonian system on~$M$ has infinitely many contractible closed orbits of integer period.

\medskip
\ni
{\it Idea of the proof of Theorem~\ref{t:Conley}.}
The proof is by contradiction: We assume that there is a Hamiltonian diffeomorphism~$\phi \colon M \rightarrow M$ with only finitely many periodic points. 
Up to replacing $\phi$ with a suitable iterate, we may assume that all the periodic points 
of~$\phi$ are fixed points of~$\phi$. We will find a contradiction by showing that there is 
a high iterate of~$\phi$ with a fixed point which is not a fixed point of~$\phi$.

The argument uses results on the behaviour of the Conley--Zehnder index under iteration. 
If $x$ is a fixed point of~$\phi$, then it has a Conley--Zehnder index $\mu_{\sCZ}(x,\phi)$. 
But being also a fixed point of~$\phi^k$, it also has a Conley--Zehnder index $\mu_{\sCZ}(x,\phi^k)$. 
These Conley--Zehnder indices grow linearly in~$k$, meaning that the limit
\[
\hat{\mu}(x) \,:=\, \lim_{k \rightarrow \infty} \frac{\mu_{\sCZ}(x,\phi^k)}{k}
\]
exists and is finite. The real number $\hat{\mu}(x)$ is called mean Conley--Zehnder index of the fixed point~$x$. It satisfies the inequalities 
\begin{equation}
\label{indexbd}
-n \,\leqslant\,  \mu_{\sCZ}(x,\phi^k) - k \, \hat{\mu}(x) \,\leqslant\, 
\mu_{\sCZ}(x,\phi^k) + \nu(x,\phi^k) - k \, \hat{\mu}(x)  \,\leqslant\, n \qquad \forall\, k\in \NN, 
\end{equation}
where $\dim M=2n$ and $\nu(x,\phi^k)$ denotes the geometric multiplicity of~1 as an eigenvalue 
of~$d\phi^k(x)$. Moreover, the two extremal inequalities are strict when the spectrum 
of~$d\phi^k(x)$ contains eigenvalues different from~1.

We now make the extra assumption that all the fixed points $x$ of~$\phi$ are 
``weakly non-degenerate'', meaning that the spectrum of~$d\phi(x)$ contains eigenvalues 
different from~1. In this case, we can find a positive integer~$k$ such that:
\begin{enumerate}[(i)]
\item 
the spectrum of $d\phi^k(x)$ contains eigenvalues different from~1 for all 
fixed points~$x$ of~$\phi$;
\item 
if $x$ is a fixed point of $\phi$ with $\hat{\mu}(x) \neq 0$, then $|\mu(x,\phi^k)| \geqslant 3n$.
\end{enumerate}
Indeed, (i) holds for any $k \in \NN$ which is not a multiple of the order of all
the (finitely many) non-trivial roots of~1 in the spectrum of the differential of~$\phi$ 
at any fixed point of~$\phi$, 
so it holds for any $k \in \NN$ in the complement of finitely many proper ideals. 
Condition (ii) holds if $k$ is large enough, by~\eqref{indexbd} and the fact that $\phi$ 
has finitely many fixed points.

Now we would like to consider the Floer complex of a Hamiltonian~$H$ whose time-one map is~$\phi^k$. This is not well-defined in general, since $\phi^k$ might have degenerate fixed points. 
However, all fixed points of~$\phi^k$ are isolated 
(they are just the finitely many fixed points of~$\phi$) and this implies that we can perturb~$\phi^k$ and obtain a non-degenerate Hamiltonian diffeomorphism~$\psi$ such that each fixed point~$x$ 
of~$\phi^k$ is replaced by finitely many fixed points $x_1,\dots,x_h$ of~$\psi$ which are near~$x$ 
and whose Conley--Zehnder indices satisfy
\[
\mu_{\sCZ}(x,\phi^k) \,\leqslant\, \mu_{\sCZ}(x_j,\psi) \,\leqslant\, 
\mu_{\sCZ}(x,\phi^k) + \nu(x,\phi^k) \qquad \forall \, j=1,\dots,h.
\]
This fact is completely analogous to the classical fact that a smooth function with isolated critical points can be perturbed and made Morse by replacing each critical point~$x$ by finitely many critical points with Morse index between the Morse index of~$x$ and the Morse index plus the nullity of~$x$.

Choose a  Hamiltonian $K$ whose time-one map is~$\psi$. The Floer complex of~$K$ is well-defined 
and as we know has the homology of~$M$ shifted by~$-n$. In particular, the fact that $\H_{2n}(M)$ does not vanish implies that the $n$-th Floer homology group does not vanish, 
and hence $\psi$ has a fixed point~$y$ of Conley--Zehnder index~$n$. 
Since all fixed points of~$\psi$ are bifurcating from fixed points of~$\phi^k$, 
$y$ originates from a fixed point~$x$ of~$\phi^k$ for which we have
\[
\mu_{\sCZ}(x,\phi^k) \,\leqslant\, \mu_{\sCZ}(y,\psi) \,\leqslant\, 
\mu_{\sCZ}(x,\phi^k) + \nu(x,\phi^k).
\]
If the mean Conley--Zehnder index $\hat\mu (x)$ vanishes, then our assumption that $x$ is a weakly non-degenerate fixed point of~$\phi^k$, by condition (i) above, implies that the 
inequalities~\eqref{indexbd} are strict and in particular that
\[
\mu_{\sCZ}(x,\phi^k) + \nu(x,\phi^k) < n.
\]
But then $\mu_{\sCZ}(y,\psi) < n$, contradicting the fact that $\mu_{\sCZ}(y,\psi) =n$. 
If $\hat \mu (x) \neq 0$, 
then condition~(ii) implies that
\[
0 = |\mu_{\sCZ}(y,\psi)| - n \,\geqslant\, |\mu_{\sCZ}(x,\phi^k)| - \nu(x,\phi^k) - n 
\,>\, 3n - 2n - n = 0,
\]
which is a contradiction. In both cases we get a contradiction, which shows that $\phi^k$ must have 
a fixed point which is not fixed by~$\phi$. This proves the theorem under the extra assumption that all fixed points of~$\phi$ are weakly non-degenerate.

A closer inspection of the above argument shows that the case which remains to be considered is the following: there is a fixed point~$x$ of~$\phi$ with $\hat \mu (x)=0$ such that for every natural number~$k$ the point $x$ is isolated in the fixed point set of~$\phi^k$ and produces, 
after perturbing $\phi^k$ to a non-degenerate Hamiltonian diffeomorphism~$\psi$, 
a fixed point~$y$ with $\mu_{\sCZ}(y,\psi)=n$. A fixed point with these properties is called 
a ``symplectically degenerate maximum''. A careful investigation of fixed points of this kind, together with a delicate argument relating the local and global behaviours of Floer homology groups under iteration, shows that the presence of a symplectically degenerate maximum implies 
the existence of infinitely many periodic points with diverging periods. 
This concludes the proof of the Conley conjecture. 

\subsection*{Historical notes and bibliography} 
The one presented here is the first version of Floer homology, which Floer created in the second half 
of the eighties, see~\cite{Flo88:index, Flo88:unreg, Flo89:hol, Flo89:Witten}. 
Hamiltonian Floer theory for symplectically aspherical manifolds was essentially complete by the time 
of Floer's tragic death in~1991, but certain aspects of it have been clarified in the period immediately afterwards. The interpretation of the grading in terms of Conley--Zehnder indices is due to Salamon and Zehnder~\cite{SaZe92}, and a more comprehensive treatment is due to Robbin and Salamon~\cite{RobSal95}. 
The first full account of the transversality theory for the Floer equation is contained in a paper of Floer, Hofer and Salamon~\cite{FHS:95}. The orientation theory which leads to Floer homology with integer coefficients is due to Floer and Hofer~\cite{FH:93}. The book~\cite{AuDa14} gives a complete and detailed account of the whole theory over $\ZZ_2$ coefficients. Previous lecture notes through which generations of students have become familiar with Hamiltonian Floer homology are~\cite{Sal99}.

If one drops the symplectic asphericity assumption, two main difficulties occur: The Hamiltonian action functional is not well-defined anymore and bubbling off of $J$-holomorphic spheres may happen. 
The first difficulty can be addressed by extending to the infinite dimensional setting Novikov's ideas 
for developing a Morse theory for zeroes of closed one-forms. The second difficulty can be overcome relatively easily when one can show that bubbling off of $J$-holomorphic spheres does not occur for moduli spaces $\mathcal{M}_{H,J}(x,y)$ of dimension~0 and~1, which are the only ones which are used in the definition of~$\partial$ and in the proof of $\partial^2=0$. 
These ideas have lead to the definition of Floer homology and to the proof of the Arnol'd conjecture 
first for monotone symplectic manifolds (by Floer~\cite{Flo89:hol}) and then for weakly monotone 
symplectic manifolds (by Hofer and Salamon~\cite{hs95} and Ono~\cite{ono95}), 
two classes of symplectic manifolds generalizing symplectically aspherical ones and including interesting examples such as the complex projective space. The case of a general closed symplectic manifold requires incorporating $J$-holomorphic spheres in the compactifications of the moduli spaces. 
This in turn makes the transversality theory much more delicate. Various solutions have been proposed 
(see \cite{lt98, fo99}) and are currently being developed 
(see \cite[Part~II]{FOOO}, \cite{Ho04, HWZ17} and the references therein, 
\cite{McWe17, McWe18}, and~\cite{Par16})
in order to solve these difficulties, which also appear in other versions of Floer homology and in 
the theory of Gromov--Witten invariants, 
but there is yet no consensus among the experts on what should be the ultimate neat solution. 

Arnol'd generalized and reformulated the first versions of his Conjecture~\ref{arnold} from~1965 
in~\cite{Arn65} several times, 
see for instance~\cite[Problem~20 on p.\ 66]{Arn76} and~\cite[Appendix~9]{Arn78}, 
and~\cite{Aud14} for the history of this and other Arnol'd conjectures.
There is also a degenerate version of the Hamiltonian Arnol'd conjecture, 
in which one does not assume the fixed points of the Hamiltonian diffeomorphism to be non-degenerate. 
In this case, Arnol'd conjectured as lower bound the minimal number of critical points of an arbitrary smooth function on~$M$, and this conjecture is usually weakened by asking as lower bound the cuplength of~$M$,
namely 1 plus the maximal length of a non-vanishing product of elements of positive degree in the cohomology
ring~$\H^*(M)$.  
When $[\omega]$ vanishes on~$\pi_2(M)$, the weak form of this conjecture has been proved by Floer 
in~\cite{Flo89:cup}, see also~\cite[\S 6.4]{HoZe94} for a very accessible presentation of this proof. 
The degenerate Arnol'd conjecture has been proved also for the complex projective space 
by Fortune~\cite{for85}, using methods more similar to the $H^{\frac 12}$ approach introduced by Conley and Zehnder for the $2n$-torus, but it is open for general closed symplectic manifolds.

The Conley conjecture for symplectically aspherical closed manifolds was proved by Salamon and Zehnder in~\cite{SaZe92} under the extra assumption that all the fixed points of the Hamiltonian diffeomorphism are weakly non-degenerate, and in general by Ginzburg~\cite{Gi10}. 
For the standard torus $\RR^{2n}/ \ZZ^{2n}$, the conjecture had been proved earlier by Hingston~\cite{Hi09}.
The asphericity assumption in Theorem~\ref{t:Conley} cannot be omitted:
Take the sphere~$S^2 \subset \RR^3$ of radius~$1$ with the Euclidean area form. 
Then the flow~$\phi_H^t$ of the Hamiltonian $H = 2 \pi \ga x_3$ is 
given by rotation by angle $2\pi \ga$ about the vertical axis.
If $\ga$ is irrational, the only periodic orbits of integral period are thus 
the two fixed points at the poles.

\begin{figure}[h]   
 \begin{center} 
  \psfrag{H}{$H = 2 \pi \ga x_3$}
  \leavevmode\includegraphics{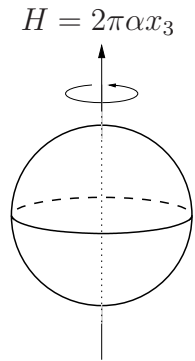}
 \end{center}
 \caption{An irrational rotation of the $2$-sphere}
 \label{fig.rot}
\end{figure}
%

\ni
This example is quite exceptional. 
In fact, the asphericity assumption in Theorem~\ref{t:Conley} has been much relaxed: 
It is enough to assume that the closed symplectic manifold~$(M,\go)$ 
contains no homology class~$A$ represented by a 2-sphere such that $[\go] (A)$ and $c_1 (A)$
are positive~\cite{GiGu16}. 
For a survey on the Conley conjecture, that also describes several variations on the theme,
see~\cite{GiGu15}.

\begin{figure}[h]   
 \begin{center}
  \psfrag{Z1}{$Z_1$}
  \psfrag{Z2}{$Z_2$}
  \psfrag{Z3}{$Z_3$}
  \psfrag{S}{$\Sigma_{\mbox{\tiny int}}$}
  \leavevmode\includegraphics{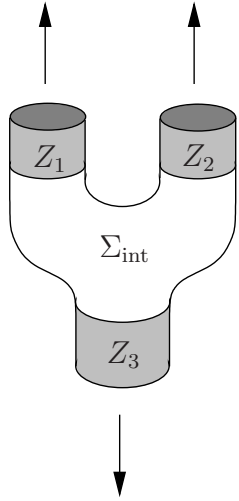}
 \end{center}
 \caption{A pair of pants}  \label{fig.pants}
\end{figure}

Algebraic structures on Floer homology that go beyond its vector space structure
are defined by extending the Floer equation to Riemann surfaces more general than 
the cylinder. 
Take, for instance, a thrice punctured sphere~$\Sigma$ with a complex structure~$j$. 
After fixing holomorphic polar coordinates $z = s+it$ with $s \in (-\infty,0)$, $t \in \TT$ 
near the two punctures~$Z_1$ and~$Z_2$ and with $s \in (0,+\infty)$, $t \in \TT$ near the third 
puncture~$Z_3$, we can think of~$\Sigma$ as the pair of pants drawn in Figure~\ref{fig.pants}.

The Floer equation on cylinders can be extended to the domain~$\Sigma$ in the following way: 
Let $\beta$ be a closed 1-form on~$\Sigma$ whose local expression is~$dt$ near~$Z_1$ and~$Z_2$ and $2 \, dt$ 
near~$Z_3$. The existence of such a 1-form is guaranteed by the fact that the sum~$1+1$ of the weights in front of~$dt$ at the negative punctures~$Z_1$ and~$Z_2$ agrees with the weight~$2$ at the positive puncture~$Z_3$. 
The Floer equation for a map $u \colon \Sigma \rightarrow M$, 
where $(M,\go)$ is a symplectic manifold equipped with an almost complex structure~$J$, 
can now be written as
\[
(du - X_H \otimes \beta)^{0,1} = 0,
\]
where $\eta^{0,1}$ denotes the anti-holomorphic part of a $u^*(TM)$-valued differential form~$\eta$ on~$\Sigma$, 
that is
\[
\eta^{0,1} = \eta - J \circ \eta \circ j.
\]
Counting the solutions of the above equation which are asymptotic to periodic Hamiltonian orbits at the three punctures defines a chain map 
\[
\CF_k (H) \otimes \CF_\ell (H) \to \CF_{k+\ell-n} (2H),
\]
where $\dim M =2n$, and hence a product 
\begin{equation} \label{e:pp}
\HF_k (M,\go) \otimes \HF_\ell (M,\go) \,\to\, \HF_{k+\ell-n} (M,\go) 
\end{equation}
that endows the Floer homology $\HF (M,\go)$ with a ring structure. This product is known as {\it pair-of-pants product}\/ and was first studied by Schwarz in~\cite{Schw95}, who used a slightly different definition involving cut-off of the Hamiltonian function away from the punctures. The formulation via the 1-form~$\beta$ outlined above is due to Seidel~\cite{sei08}. When the closed symplectic manifold~$(M,\omega)$ is symplectically aspherical, this product corresponds to the cup product in singular cohomology~$\H^*(M)$, after identifying the Floer homology of~$(M,\omega)$ with the singular cohomology of~$M$ by the isomorphism~\eqref{iso:HF} and Poincar\'e duality.
If one drops the symplectic asphericity condition, 
then $M$ can contain $J$-holomorphic spheres and the pair of pants product
is a finer symplectic invariant, which agrees with the quantum cup product from Gromov--Witten theory, 
see~\cite{PiSaSch94},~\cite{LT99} and~\cite[\S 12.2]{McSal12}. 
An extensive study of finer algebraic structures on Floer homology is contained in~\cite{Rit13}.


\section{Basics on Hamiltonian dynamics and symplectic geometry, II}
\label{s:basicsII}

In order to discuss more advanced applications of Floer's ideas, 
we need to introduce further concepts from Hamiltonian dynamics and symplectic geometry.
We start with contact type hypersurfaces and their characteristic flows,  
which naturally lead to the more general notion of Reeb flows on contact manifolds. 
We then introduce Lagrangian and Legendrian submanifolds, and finally have a first look at symplectic embedding
problems.

\m \ni
{\bf From hypersurfaces to contact manifolds.}
Let $(M,\omega)$ be a $2n$-dimensional symplectic manifold. 
If a Hamiltonian function $H \colon M \rightarrow \R$ does not depend on time, 
then its level sets are preserved by the Hamiltonian flow:
\[
H(\phi^t_H(x)) = H(x).
\]
In particular, if $c$ is a regular value of $H$ then $\Sigma := H^{-1}(c)$ is a smooth hypersurface in~$M$ 
which is invariant under the Hamiltonian flow, and the Hamiltonian vector field~$X_H$ restricts to a nowhere vanishing tangent vector field on~$\Sigma$. 
One can thus restrict $\phi_H^t$ to $\Sigma$
and study the properties of this restricted dynamical system.
Actually, the direction of~$X_H$ depends only on~$\Sigma$, 
and not on the particular Hamiltonian~$H$ of which $\Sigma$ is a regular level. 
Indeed, $\RR X_H$ is precisely the kernel of the restriction of~$\omega$ to~$\Sigma$, 
which is called the {\bf characteristic line field} of~$\Sigma$.  
Therefore, the choice of a different Hamiltonian~$H$ still having~$\Sigma$ as regular level determines 
a flow on~$\Sigma$ which is a time reparametrisation of the original flow. 
When one ignores the time parametrisation, Hamiltonian orbits on~$\Sigma$ are just integral flow lines 
of the characteristic line field and they are called {\bf characteristics}.
Many properties of the flow $\phi_H^t$ depend only on the characteristic foliation on~$\Sigma$,
but not on the specific choice of~$H$. Examples are the set of closed orbits (which are in bijection with the closed characteristics) and (at least if $\Sigma$ is compact) also the asymptotics 
for $T \rightarrow +\infty$ of the number of closed orbits with period not exceeding~$T$, 
and the positivity of the topological entropy 
(cf.\ \S \ref{s:entropy}).

\s
The hypersurface $\Sigma$ is said to be of {\bf contact type} when the restriction of~$\omega$ to a neighbourhood of~$\Sigma$ admits a primitive~$\lambda$ such that its 
{\bf Liouville vector field}~$X_\lambda$ defined by
\[
\omega(X_{\lambda},\cdot) = \lambda
\]
is transverse to~$\Sigma$. 
Equivalently, $\lambda \wedge d\lambda^{n-1}$ restricts to a volume form on~$\Sigma$. 
This second formulation suggests that we study a certain class of flows
on an odd-dimensional manifold without considering an ambient symplectic manifold explicitly:
A smooth 1-form on a $(2n-1)$-dimensional manifold~$Y$ is said to be a {\bf contact form} 
if $\alpha \wedge d \alpha^{n-1}$ is a volume form on~$Y$. The kernel of~$\alpha$ is called the {\bf contact structure} on~$Y$, 
and defines~$\alpha$ only up to multiplication by a non-vanishing real function. 
The non-vanishing vector field $R_{\alpha}$ defined by
\[
d\alpha(R_{\alpha},v) = 0 \quad \forall v \in TY, \qquad \alpha(R_{\alpha}) = 1,
\]
is called the {\bf Reeb vector field} of $\alpha$. 

If $Y \subset M$ is a contact type hypersurface of $(M,\omega)$ with respect to a primitive $\lambda$ 
of $\omega$ near~$Y$, then $\alpha := \lambda|_Y$ is a contact form on~$Y$ and $R_{\alpha}$ 
generates the characteristic line field. Therefore, the flow of $R_{\alpha}$ agrees, up to a time reparametrisation,  with the flow of any autonomous Hamiltonian vector field which is defined by a Hamiltonian having $Y$ as regular level.

Conversely, if $\alpha$ is a contact form on $Y$ then there is a standard way of seeing $Y$ as a contact type hypersurface inside a symplectic manifold: We consider the product $M := (0,+\infty) \times Y$ with the one-form
\[
\lambda = \rho \alpha,
\]
where $\rho$ is the coordinate on the first component in the product $(0,+\infty) \times Y$, 
and endow~$M$ with the symplectic form $\omega = d \lambda$. 
The pair $(M,\omega)$ is called the {\bf symplectisation} of $(Y,\alpha)$. 
Notice that the Liouville vector field induced by~$\lambda$ has the form
\[
X_{\lambda} = \rho \frac{\partial}{\partial \rho}.
\]

A {\bf contact manifold} is a manifold on which a contact structure has been fixed. 
One then has the pair~$(Y,\xi)$, where $\xi = \ker \alpha$ is the kernel of one
and hence any other contact form $f \ga$ on~$(Y,\xi)$.
The Reeb flow $\phi_{f\ga}^t$ of $R_{f \ga}$ strongly depends on the function~$f$,
see Examples~\ref{ex:star} and~\ref{ex:spherisation} below.
One can nevertheless ask for properties of the Reeb flows $\phi_{f\ga}^t$ 
that do not depend on the specific choice of the contact form:

\m \ni
{\bf 1.\ Weinstein conjecture.}
Every Reeb flow on a closed contact manifold has a closed orbit.

\s \ni
{\bf 2.\ Positive topological entropy.} 
Are there contact manifolds $(Y,\xi)$ for which every Reeb flow 
has positive topological entropy?

\m
We shall see that Floer homologies have much to say about both questions
(Sections~\ref{s:SH}, \ref{s:disc} and~\ref{s:ech} for the the Weinstein conjecture and 
Section~\ref{s:entropy} for positive entropy).
The Weinstein conjecture in particular holds for the boundary of many Liouville domains.

\m \ni
{\bf Liouville domains.}\  
A Liouville domain is a compact symplectic manifold~$(M,\omega)$ with boundary~$\pp M$ 
such that $\omega$ is exact and admits a primitive $\lambda$ whose induced Liouville vector field $X_{\lambda}$ 
is transverse to~$\partial M$. 
Therefore, a Liouville domain is a symplectic manifold with contact type boundary, 
such that the primitive of $\omega$ giving the contact condition near the boundary 
extends to a primitive on the whole~$M$.

\begin{example}[\bf Starshaped domains in $\RR^{2n}$] \label{ex:star}
{\rm
Let $U \subset \R^{2n}$ be a bounded open set containing the origin and such that $\partial U$ is a smooth hypersurface transverse to the radial direction. Then $(\overline{U},\omega_0|_{\overline U})$ is a Liouville domain with respect to the following primitive of~$\omega_0$:
\[
\lambda_0 := \frac{1}{2} \sum_{j=1}^n \bigl( p_j \, dq_j - q_j \, dp_j \bigr).
\]
Indeed, the Liouville vector field induced by $\lambda_0$ is the radial vector field
\[
X_0 := \frac{1}{2} \sum_{j=1}^n \left( p_j \frac{\partial}{\partial p_j} + q_j \frac{\partial}{\partial q_j} \right) \,=\, \frac 12 \,r \frac{\partial}{\partial r} .
\]
The boundary $\pp U$ carries the contact structure $\xi_{\pp U} = \ker \lambda_0 |_{\pp U}$.
Given another starshaped domain~$V$, the radial projection $\pp U \to \pp V$ maps $\xi_{\pp U}$ to $\xi_{\pp V}$.
All these contact manifolds are thus diffeomorphic, 
and define the canonical contact structure on the sphere~$S^{2n-1}$.
On the other hand, the Reeb flows on~$\pp U$ and $\pp V$ with respect to~$\lambda_0$ can be
very different. Consider, for instance, the ellipsoids
\begin{equation}
\label{ellipsoid}
\E (a_1, \dots, a_n) \,=\, 
\left\{ (z_1, \dots, z_n) \in \CC^n \mid \sum_{j=1}^n \frac{\pi |z_j|^2}{a_j} < 1 \right\} 
\end{equation}
whose projection to the $j$-th complex coordinate plane is the disc of area $a_j$.
Here, we describe an ellipsoid by its areas $a_j= \pi r_j^2$, not its radii~$r_j$, 
since lengths and distances have no intrinsic meaning in symplectic geometry.
Then the Reeb flow on the boundary of $\E (a_1, \dots, a_n)$ has only $n$ simple closed orbits
if the $a_j$ are rationally independent 
(namely the orbits $\gamma_j = \{ z_k = 0 \mbox{ if } k \neq j \})$, 
but has infinitely many simple closed orbits otherwise.
}
\end{example}

\begin{example}[\bf Starshaped domains in cotangent bundles]  \label{ex:spherisation}
{\rm 
Let $Q$ be a compact manifold and let $U \subset T^* Q$ be a bounded open set containing the zero section 
and such that for every $q \in Q$ the set $\partial U \cap T_q^* Q$ is a smooth hypersurface in~$T^*_q Q$ transverse to the radial direction. Then $\overline{U}$ is a Liouville domain with respect to the Liouville form $\lambda_{\can} = \sum_j p_j\, dq_j$, which induces the canonical symplectic 
form $\omega_{\can}$ on~$T^* Q$. 
Indeed, the Liouville vector field induced by $\lambda_{\can}$ is the fiberwise radial vector field
\[
X_{\can} = \sum_{j=1}^n p_j \frac{\partial}{\partial p_j}.
\]

Interesting examples in this class are given by the cotangent disc bundles which are determined by a Riemannian structure on~$Q$ and a Hamiltonian modelling a magnetic and potential force:
Specializing~\eqref{e:classical} we consider the autonomous Hamiltonian
\begin{equation} \label{e:AV}
H(q,p)  \,=\, \tfrac 12\, \| p-\ga (q) \|^2 + V(q) .
\end{equation}
When $c> \max V$, the intersection of $H^{-1}(c)$ with $T_q^*Q$ is a sphere centered 
at~$\ga (q)$, for every $q \in Q$. If $c$ is even larger, and precisely if
\begin{equation} \label{threshold}
c \,>\, \max_{q \in Q} \left( \tfrac 12 \| \ga (q) \|^2 + V(q)\right) \,=:\, \hat \kappa,
\end{equation}
the origin of $T_q^* Q$ belongs to $\{H<c\}$ and hence
the energy level $H^{-1}(c)$ is a fiberwise starshaped hypersurface,
see the middle drawing in Figure~\ref{fig.star}.
Actually, even for energies below the threshold~$\hat \kappa$ the energy level~$H^{-1}(c)$ can be mapped by a symplectomorphism to a 
fiberwise starshaped hypersurface: 
This is certainly the case if $c$ is larger than the number
\[
\kappa \,:=\, \inf_{\theta} \max_{q\in Q} H(q,\theta(q)) \,=\, \inf_{\theta} \max_{q \in Q} \left( \tfrac 12\, \| \theta(q) - \ga (q) \|^2 + V(q) \right) ,
\]
where the infimum is taken over all closed 1-forms $\theta$ on~$Q$. Indeed, if $c>\kappa$ then there is a closed 1-form $\theta$ such that the hypersurface $H^{-1}(c)\cap T_q^* Q$ is starshaped with respect to the point~$\theta(q)$, for every $q \in Q$. 
The diffeomorphism
\[
T^*Q \rightarrow T^* Q, \quad (q,p) \mapsto \bigl(q,p-\theta(q)\bigr),
\]
which is symplectic because the 1-form $\theta$ is closed, maps $H^{-1}(c)$ to a fiberwise 
starshaped hypersurface~$\Sigma$.
More precisely $\Sigma$ is fiberwise convex: Its intersection~$\Sigma_q$ with~$T^*_q Q$ is a sphere centred
at $\alpha(q)-\theta(q)$ and with the origin in its interior,
see the right drawing in Figure~\ref{fig.star}.
Therefore, $\Sigma$ is the unit sphere cotangent bundle associated to a smooth 
(possibly non-reversible) Finsler metric on~$Q$.
Since symplectomorphisms preserve the characteristic foliations, 
we conclude that for $c>\kappa$ the Hamiltonian flow on~$H^{-1}(c)$ is conjugated, 
after a time reparametrisation, to a Finsler geodesic flow.

The number $\kappa$ is called Ma\~n\'e critical value.
It is between the maximum of~$V$ and the threshold~$\hat \kappa$ in~\eqref{threshold},
but often smaller than~$\hat \kappa$, 
and it has many other interesting characterisations
in terms of the Lagrangian which is Legendre dual to~$H$ 
(see~\eqref{e:Leg} for the definition)
and in terms of Lagrangian graphs in $T^*Q$
(defined below).
The Ma\~n\'e critical value $\kappa$ is therefore a stronger and more natural energy threshold 
than $\hat \kappa$.
Many results on this value can be found in~\cite{CieFraPat10, CoItPa98}, 
where the former paper uses the Rabinowitz--Floer homology outlined in~Section~\S \ref{s:RF}.

\begin{figure}[h]   
 \begin{center} 
  \psfrag{X}{$X$}
  \psfrag{aq}{$\alpha(q)$}
  \psfrag{tq}{$\theta(q)$}
  \psfrag{at}{$\alpha(q)-\theta(q)$}
  \psfrag{S}{$\Sigma_q$}
  \psfrag{p1}{$p_1$}
  \psfrag{p2}{$p_2$}
  \leavevmode\includegraphics{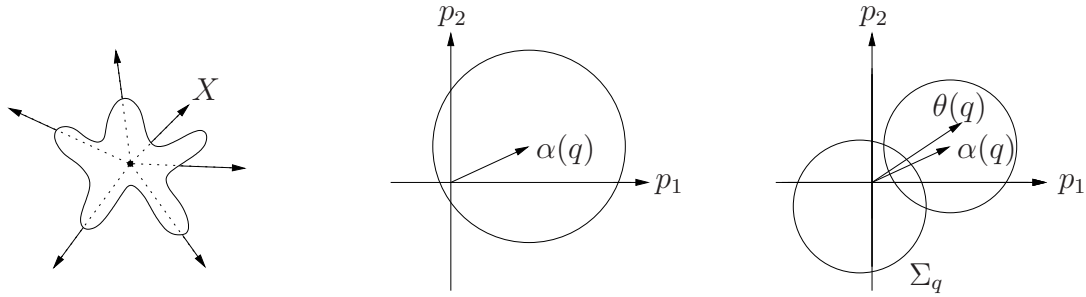}
 \end{center}
 \caption{A starshaped domain in $\RR^2$, 
a fiber sphere $H^{-1}(c) \cap T_q^*Q$ for the Hamiltonian system~\eqref{e:AV}, 
and a fiber sphere $\Sigma_q$ obtained after translation by $-\theta(q)$}
 \label{fig.star}
\end{figure}
%

Given two starshaped domains $U$ and $V$ in~$T^*Q$, the fiberwise radial projection 
identifies the contact structures $\ker \lambda_{\can}$ on $\pp U$ and~$\pp V$.
The resulting contact manifold is called the {\bf spherisation} $(S^*Q, \xi_{\can})$ of~$Q$. 
The Reeb flows on $\pp U$ and~$\pp V$ with respect to~$\lambda_{\can}$, 
and hence different Reeb flows on $(S^*Q, \xi_{\can})$,
can be very different. 
For instance, the geodesic flow on the round sphere is periodic, while 
for most other Riemannian metrics there are only finitely many closed geodesics below a given length.
}
\end{example}

\begin{example}[\bf Stein domains]
{\rm
Let $(W,J)$ be a Stein manifold, that is a complex ma\-ni\-fold endowed with a coercive plurisubharmonic function $h \colon W \rightarrow \R$. Then each regular sublevel $\{h \leqslant c\}$ is a Liouville domain with respect to the one-form $\lambda = - d^c h = dh\circ J$. Indeed, in this case the Liouville vector field is the gradient of~$h$ with respect to the metric determined by the K\"ahler form $d\lambda = -dd^c h$.
}
\end{example}

Algebraic geometry is another rich source of Liouville domains, see \cite{sei08}. 

\m \ni
{\bf Completion.}\
Let $(M,\lambda)$ be a Liouville domain. The Liouville vector field $X_{\lambda}$ is now defined 
on all of~$M$, and since its flow~$\psi^t$ expands~$\omega$,
$$
(\psi^t)^* \omega \,=\, e^t \omega,
$$
it also expands the volume form and hence $X_\lambda$ is outward-pointing along the boundary. 
Therefore, its flow $\phi_{X_{\lambda}}^s$ is defined for all $s \leqslant 0$ and induces an embedding
\[
j \colon  (0,1] \times \partial M \rightarrow M, \quad (\rho,x) \mapsto \phi^{\log \rho}_{X_{\lambda}}(x),
\] 
such that $j^* \lambda = \rho \,\alpha$ and $j^* X_{\lambda} = \rho \frac{\partial}{\partial \rho}$,
where again $\alpha = \lambda |_{\pp M}$ and 
$\rho$ denotes the coordinate on the first factor of $(0,1] \times \partial M$. 
The {\bf completion} of~$M$ is the manifold
\[
\widehat{M} \,=\, M \cup_{\partial M} ( [1+\infty) \times \partial M).
\]
The Liouville form~$\lambda$ and hence the symplectic form~$\omega$ extend naturally to the completion by setting 
\[
\lambda|_{[1+\infty) \times \partial M} = \rho \,\alpha .
\]
In other words, $\widehat{M}$ is obtained by gluing to $M$ the positive part of the symplectisation 
of~$\partial M$.

\begin{figure}[h]   
 \begin{center}
  \psfrag{M}{$M$}
  \psfrag{dM}{$\pp M$}
  \psfrag{1}{$1$}
  \psfrag{0}{$0$}
  \psfrag{r}{$\rho$}
  \leavevmode\includegraphics{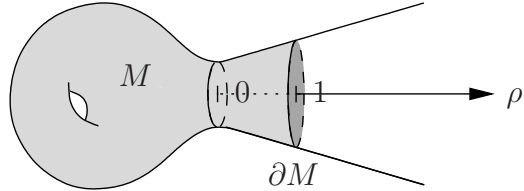}
 \end{center}
 \caption{The completion $\widehat M$ of $M$}  \label{fig.comp}
\end{figure}

One easily checks that the completion of a starshaped domain in $\R^{2n}$ is $(\R^{2n},\lambda_0)$, 
while the completion of a starshaped domain in~$T^* Q$ is $(T^* Q,\lambda_{\can})$. 
In general, the completion of a regular sublevel $\{h \leqslant c\}$ in the Stein manifold $(W,J)$ with plurisubharmonic function~$h$ has a simpler topology than~$W$, 
unless $c$ is larger than the largest critical value of~$h$.

\m \ni
{\bf Lagrangian submanifolds.}\
A Lagrangian submanifold of a $2n$-dimensional symplectic manifold $(M,\go)$ is an $n$-dimensional submanifold~$L$ on which the restriction of~$\go$ vanishes:
$\go_x (v,w) = 0$ for all $v,w \in T_x L$ and $x \in L$.

\begin{examples}
{\rm
(1)
Let $T^*Q$ be the cotangent bundle over a smooth manifold~$Q$,
endowed with the canonical symplectic form $\omega_{\can} = d \lambda_{\can}$.
Then the zero section $Q = \{ (q,p) \mid p = 0\}$ is Lagrangian, 
and every cotangent fiber $T_q^*Q$ is Lagrangian.

\s 
(2)
Every embedded circle on a surface with an area form is Lagrangian. In particular, 
any circle $S^1$ in $\CC \cong \RR^2$ is Lagrangian, 
and so the torus $T^n = S^1 \times \dots \times S^1 \subset \CC^n$ is also Lagrangian.

\s 
(3)
Let $V = \{ z \in \CP^n \mid p(z)=0 \}$ be a smooth submanifold in complex projective space
given as the zero set of a homogenous polynomial (a smooth projective variety). 
Then the real part $V_{\RR} = \{ x \in \RP^n \mid p(x)=0 \}$
is a Lagrangian submanifold of $V$.
An example is $\RP^n \subset \CP^n$ itself.  

\s 
(4)
If $L \subset (M,\go)$ is Lagrangian and $\gf$ is a symplectomorphism of $(M,\go)$,
then $\gf (L)$ is Lagrangian.
In particular, the graph $\Gamma_{\gf} = \{ (x, \gf (x)) \mid x \in M \}$  
of any symplectomorphism~$\gf$ of~$M$ is a Lagrangian submanifold of $(M \times M, \go \oplus (-\go))$.
\diam
}
\end{examples}

The Lagrangian submanifolds are (together with energy surfaces) the most interesting submanifolds of symplectic manifolds for several reasons. 
One reason is that these are the submanifolds that exhibit 
``symplectic rigidity''.
Another reason is that many problems in symplectic geometry and dynamics can be reformulated as problems 
on Lagrangian submanifolds,
\begin{equation} \label{dictum:Wein}
\mbox{``Everything is a Lagrangian submanifold''~\cite{Wei81}}
\end{equation}

We shall see examples for both reasons in \S~\ref{s:Lag}.

\m
Recall from the end of Section~\ref{s:basicsI}
that for every Hamiltonian diffeomorphism~$\phi_H$ of $T^*S^1$,
the image $\phi_H(S^1)$ of the zero section intersects~$S^1$ in at least two points
(and there need be no intersection points for an arbitrary symplectic diffeomorphism).
Now let $Q$ be any closed smooth manifold, 
and let $\phi_H$ be a Hamiltonian diffeomorphism of~$T^*Q$.
If $\phi_H$ is $C^1$-small, then $\phi_H (Q)$ is the graph
$\{ (q,\ga(q))\}$ of a 1-form~$\ga$ on~$Q$.
Since $\phi_H$ is symplectic, the $1$-form $\ga$ is closed, 
and since $\phi_H$ is Hamiltonian, $\ga$ is even exact: 
$\ga = df$ for a smooth function $f \colon Q \to \RR$.
Now assume that $\phi_H(Q)$ and $Q$ intersect transversally.
This is equivalent to $f$ being a Morse function. 
Hence the Morse inequalities from Theorem~\ref{t:Morse} imply that 
\begin{equation} \label{e:MorseLag}
\# (\phi_H(Q) \cap Q) \,\geqslant\, \sum_{k=1}^n b_k(Q).
\end{equation}
Arnol'd conjectured that this estimate 
holds even without the smallness assumption on~$\phi_H$.

\m \ni
{\bf Arnol'd conjecture for Lagrangian intersections in $T^*Q$.}\
{\it
Let $Q$ be a closed $n$-dimensional manifold, and let $\phi_H$ be any Hamiltonian diffeomorphism of~$T^*Q$ such that $\phi_H(Q)$ and $Q$ intersect transversally. 
Then~\eqref{e:MorseLag} holds true.}

\m 
Let $L$ be a Lagrangian submanifold of an exact symplectic manifold $(M,\omega=d\lambda)$. 
Then $\omega|_L=d\lambda|_L=0$, and hence the restriction of~$\lambda$ to~$L$ is closed. 
If $\lambda |_L$ happens to be exact, then $L$ is said to be an {\bf exact Lagrangian submanifold}. 
Every Lagrangian submanifold~$L$ with vanishing first cohomology group is necessarily exact.
Examples in~$(T^*Q,\omega_{\can}=d\lambda_{\can})$ are graphs of exact 1-forms and, more generally, 
the image $\phi_H (Q)$ of~$Q$ under a Hamiltonian diffeomorphism.
Arnol'd conjectured in~\cite{Arn86} that there are no other examples:

\m \ni
{\bf Arnol'd nearby Lagrangian conjecture.}\
{\it Every closed exact Lagrangian submanifold~$L$ in the cotangent bundle 
$(T^*Q,\omega_{\can}=d\lambda_{\can})$ of a closed manifold~$Q$ is the image of the zero section~$Q$ under a Hamiltonian diffeomorphism.}

\m
It is a nice and simple exercise to check the conjecture for $Q=S^1$
(cf.\ Figure~\ref{fig.flux}).

\m \ni
{\bf Legendrian submanifolds.}\
These are the analogues of Lagrangian submanifolds 
for contact manifolds. 
A Legendrian submanifold of a contact manifold $(Y,\xi)$ of dimension~$2n+1$ 
is an $n$-dimensional submanifold that is everywhere tangent to~$\xi$.
Examples are the fibres~$S_q^*Q$ of a spherisation $(S^*Q, \xi_{\can})$.
A compact connected Legendrian submanifold in a three-dimensional contact manifold
is called a {\bf Legendrian knot}.
Given a contact form $\lambda$ for~$\xi$, a {\bf Reeb chord}\/ of a Legendrian knot~$\kk$
is a trajectory of the Reeb vector field starting and ending on~$\kk$,
i.e., a path $\gg \colon [0,T] \to Y$ for some $T>0$ with $\dot \gg (t) = R_{\lambda}(\gg(t))$
and $\gg(0), \gg(T) \in \kk$.
Here are two examples.

\s
{\bf 1.}\
For any point $q$ on a closed surface~$\Sigma$ with Riemannian metric~~$g$ 
there exists a (non-constant) geodesic starting and ending at~$q$.
In contact dynamical terms, this statement generalizes to:
Given a point $q$ in the base of the contact manifold $(S^* \Sigma, \xi_{\can})$
and given any contact form~$\lambda$ for~$\xi_{\can}$ there exists a Reeb chord for
the Legendrian knot $S_q^* \Sigma$.

\s
{\bf 2.}\
Take a Legendrian knot $\kk$ in the standard contact manifold $(\RR^3, dz + xdy)$,        
with a regular parametrisation on~$[0,1]$.
The Legendrian condition $(dz + xdy) (\dot \kk (t)) =0$ is $\dot z (t) = -x(t) \dot y(t)$ for all $t$.
The projection ${\kk}_\pi (t) = (x(t), y(t))$ to the $xy$-plane 
is thus an immersed curve, whose signed area $\int_{\kk} x dy = - \int_{\kk} dz =0$ vanishes.
Let $q := {\kk}_\pi(t_1) = {\kk}_\pi(t_2)$ be a self-intersection point,
with $t_1 \neq t_2$ such that $z(t_1) < z(t_2)$.
Since the Reeb vector field is just $\frac{\pp}{\pp z}$,
the straight line from $(q, z(t_1))$ to $(q, z(t_2))$ is then a Reeb chord for~$\kk$.

\m
Based on these examples,
Arnol'd conjectured in~\cite[\S 8]{Arn86} the following relative version of the Weinstein conjecture: 

\m \ni
{\bf Arnol'd chord conjecture.}\
{\it Every Legendrian knot in $\RR^3$ with the standard contact structure $\ker (dz+xdy)$ has a Reeb chord, for every contact form for this contact structure.}

\m \ni
{\bf Symplectic embedding problems.}\
By ``symplectic rigidity'' one means results that show that (certain) symplectic mappings are 
more rigid than volume preserving or smooth mappings. 
One example is the Lagrangian intersection phenomenon addressed in the previous section. 
Another one is the rigidity of mappings between symplectic manifolds.
For simplicity we restrict to simply connected domains in~$(\RR^{2n},\go_0)$.
When are two such sets $U$ and $V$ symplectomorphic? 
An obvious necessary condition is that $U$ and $V$ have the same total volume 
$\Vol (U) = \frac 1{n!} \int_U \go_0^n$,
since symplectomorphisms are volume-preserving, 
$$
\gf^* (\go_0^n) \,=\,(\gf^* \go_0)^n \,=\, \go_0^n .
$$
Are there other conditions?
Similarly, a symplectic embedding $\gf \colon U \to V$ can only exist if
$\Vol (U) \leqslant \Vol (V)$.
It is not hard to see that for volume preserving embeddings, this is the only obstruction~\cite{Sch03}.
Symplectic mappings, however, are much more rigid, 
as was first shown by Gromov in his famous non-squeezing theorem from~\cite{Gro85}.

\begin{theorem}
[\bf Non-squeezing] \label{t:Gromov}
Endow $\RR^{2n} \cong \CC^n$ with the standard symplectic structure $\omega_0$. 
If $0<a<A$ then one cannot symplectically embed the ball 
\[
\B^{2n} (A) := \{(z_1, \dots, z_n) \in \CC^n \mid \sum_{j=1}^n \pi |z_j|^2 < A \}
\]
into the cylinder 
\[
\Z^{2n} (a) := \{(z_1, \dots, z_n)\in \CC^n \mid \pi |z_1|^2 < a \}.
\]
\end{theorem}

Hence there is no symplectic embedding of a ball into a cylinder that ``does a better job'' than the inclusion.
Gromov's original proof of this result is based on $J$-holomorphic curves, the same objects that 
a few years later Floer modified in order to attack the Arnol'd conjecture, 
as we have seen in Section~\ref{s:Ham}. 
Other proofs are based on existence results for periodic Hamiltonian orbits, 
and in Section~\ref{s:emb} we shall give a proof based on Floer homology, 
that intertwines $J$-holomorphic curves with periodic orbits.

In order to get a feeling on how periodic Hamiltonian orbits can be used to prove rigidity results 
as the above one, let us prove here the following weaker version of Gromov's non-squeezing theorem: 
There is no smooth path $\{\phi^t\}_{t\in [0,1]}$ of Hamiltonian diffeomorphisms of $(\R^{2n},\omega_0)$ 
such that $\phi^0 = \mathrm{id}$ and
\[
\phi^t(\B^{2n} (\pi)) \subset \Z^{2n} (a(t)) \qquad \forall \,t\in [0,1]
\]
for some function $a \colon [0,1] \rightarrow (0,+\infty)$ with $a(0)=\pi$ and $a'(0)<0$. 
In other words, we cannot instantaneously squeeze a ball of radius~1 into a thinner cylinder 
by means of a Hamiltonian deformation.

Let $H_t$ be the time-dependent Hamiltonian on $\RR^{2n}$ whose flow is the path~$\phi^t$ 
and consider the function
\[
K \colon \CC^{n} \rightarrow \R, \qquad K(z_1,\dots,z_n) = \pi |z_1|^2.
\] 
If $z \in \overline{\B^{2n}(\pi)}$, then the condition on~$\phi^t$ implies that 
$K(\phi^t(z)) \leqslant a(t)$ for all $t \in [0,1]$. If $z$ belongs to the circle 
\[
\Gamma = \{ (z_1,0,\dots,0) \mid |z_1|=1 \}
\]
then $K(\phi^0(z)) = K(z) = \pi$. For such a $z$, the real function $K(\phi^t(z))$ does not exceed 
$a(t)$ and coincides with it for $t=0$. Therefore
\[
\frac{d}{dt} K(\phi^t(z)) \Big|_{t=0} \leqslant \frac{d}{dt} a(t) \Big|_{t=0} = a'(0) < 0.
\]
This quantity can be written as 
\[
\begin{split}
\frac{d}{dt} K(\phi^t(z)) \Big|_{t=0} &= dK(z)[X_{H_0}(z)] = 
-\omega_0(X_K(z), X_{H_0}(z)) = \omega_0(X_{H_0}(z), X_K(z)) \\ 
&= -dH_0(z)[X_K(z)] .
\end{split}
\]
Now notice that $\Gamma$ is the image of a periodic orbit of the Hamiltonian flow given by 
the Hamiltonian~$K$. 
Then the fact that the above quantity is negative tells us that the function~$H_0$ increases strictly 
along this periodic orbit, which is of course impossible. This contradiction proves the claim.

Here, the presence of a periodic orbit for the Hamiltonian~$K$ obstructs the possibility of instantaneously squeezing~$\B^{2n}(\pi)$ into a cylinder~$\Z^{2n}(a)$ with $a<\pi$. 
So it will be not so surprising that more powerful tools to deal with Hamiltonian periodic orbits, such as Floer homology, 
are capable to prove much deeper rigidity results.

\subsection*{Historical notes and bibliography}
Hamiltonian flows and even the concept of a symplectic form were known to Lagrange, 
see~\cite{Alb13} for a discussion.
After the introduction of important new concepts by Poincar\'e, 
the geometrisation of classical mechanics 
(by the introduction of notions such as Lagrangian and Legendrian submanifolds) 
is largely due to Arnol'd \cite{Arn78,AG85}.
Gromov's introduction of $J$-holomorphic curves to symplectic geometry in~\cite{Gro85}
and Floer's creation of his homology 
can be seen as higher levels of geometrisation 
of classical mechanics.
Contact geometry has its origins in Huygens' work on geometric optics.
The study of contact structures was initiated by Eliashberg in~\cite{El89},
and Hofer's compactness theorem in symplectisations from~\cite{Ho93} 
opened the door for studying contact dynamics by $J$-holomorphic curves and Floer homology.

Excellent books covering all basic aspects of Hamiltonian dynamics and symplectic geometry
are~\cite{HoZe94} and~\cite{McSa15},
where the focus of the first text is more on dynamics while the one of the second text
is more on topology.
The book~\cite{Gei08} is the classic reference for contact geometry.

The Weinstein conjecture goes back to the old question in Hamiltonian dynamics 
whether an energy level $H^{-1}(c)$ admits periodic orbits. 
This conjecture was (together with the Arnol'd conjecture) an important driving force 
in the development of symplectic and contact dynamics.
%
%
There are examples of compact hypersurfaces in $(\RR^{2n},\omega_0)$ without periodic orbits, 
see~\cite{Gi97, He99}, but
Weinstein~\cite{Wei78} and Rabinowitz~\cite{Rab78} proved that on every  
convex resp.\ starshaped hypersurfaces in $\RR^{2n}$ there is a closed characteristic 
(see also \cite{cla81} for a more elegant proof of the Weinstein conjecture for convex hypersurfaces).
`Starshaped' is not a symplectically invariant notion. 
Weinstein thus looked for a symplectically invariant notion and formulated his conjecture 
in~\cite{Wei79} for compact hypersurfaces of contact type in arbitrary symplectic manifolds.
Weinstein also made the assumption $\H^1(S;\RR)=0$ on the hypersurface, 
which however played no role in any of the subsequent proofs.
This form of the Weinstein conjecture was proved by Viterbo~\cite{Vi87} in~$\RR^{2n}$, 
see Section~\ref{s:SH} for a Floer homology proof,
and by Hofer and Viterbo~\cite{HoVi89} in spherisations $(S^*Q, \xi_{\can})$,
cf.\ Section~\ref{s:disc},
thereby generalizing the Lusternik--Fet theorem on the existence of a closed geodesic on any closed Riemannian manifold.  

Now, by `Weinstein conjecture' one usually means the existence of a periodic orbit for any Reeb vector field on a closed contact manifold.
This extended Weinstein conjecture has been proved in several cases, 
but in full generality only in dimension three:
Hofer~\cite{Ho93} proved the conjecture for $S^3$ and for all closed contact 3-manifolds whose second homotopy group does not vanish, and the latter condition was removed by Taubes~\cite{Ta07},
cf.\ the notes to Section~\ref{s:contact}.
A survey on the Weinstein conjecture is given in~\cite{Gi05}.

\section{Symplectic invariants from Floer homology: symplectic homology}
\label{s:symphom}

Shortly after the birth of Floer homology, Floer and Hofer realized that this theory can be used not only to prove existence statements about periodic Hamiltonian orbits, 
but also to construct symplectic invariants and hence to address rigidity questions in symplectic geometry. 
In this section we review the simplest version of {\em symplectic homology}, a theory which uses Floer homology to attach non-trivial invariants to certain symplectic manifolds with boundary. 
We will present the version of symplectic homology due to Viterbo which is defined for Liouville domains. 
We then apply this Floer homology to prove the Weinstein conjecture in $\RR^{2n}$, to show that there are no closed exact Lagrangian
submanifolds in~$\RR^{2n}$, and to prove the non-squeezing theorem.

\subsection{Outline of the construction} 
\label{s:SH}
Let $(M,\lambda)$ be a Liouville domain. We shall also assume that the first Chern class of~$TM$ vanishes identically, so that every Hamiltonian periodic orbit on~$M$ has a well-defined Conley--Zehnder index. We would like to define the symplectic homology of~$(M,\lambda)$ as the homology
of the Hamiltonian Floer complex given by an autonomous Hamiltonian~$H$ on the completion~$\widehat M$
of~$(M,\lambda)$ which vanishes identically on~$M$ and is a function of 
the radial variable~$\rho$ on~$\widehat{M}\setminus M$,
\[
H = h(\rho) \qquad \forall \rho \geqslant 1,
\]
where $h \colon [1,+\infty) \rightarrow \R$ is a function such that 
$h''(\rho) > 0$ for $\rho >1$,
\begin{equation}
\label{condh}
h(1)=h'(1) = 0 \qquad \mbox{and} \qquad \lim_{\rho\rightarrow +\infty} h'(\rho)=+\infty.
\end{equation}
Notice that by the exactness of $\omega$, we can use the expression~\eqref{withprimitive} 
for the action functional and deal with the whole loop space of~$\widehat{M}$, 
and not just with the component consisting of contractible loops. 
Let us see what the 1-periodic orbits of~$X_H$ and their action are. All the points in~$M$ are critical points of~$H$ with $H=0$, and hence stationary solutions of~$X_H$ with vanishing action. On the set $\widehat{M} \setminus M = (1,+\infty) \times \partial M$, the Hamiltonian vector field of~$H$ has the form
\[
X_H(\rho,x) = h'(\rho)\, R_{\alpha}(x),
\]
where $R_{\alpha}$ is the Reeb vector field of the contact form $\alpha = \lambda|_{\partial M}$. Therefore, a periodic orbit~$x$ of~$R_{\alpha}$ of (not necessarily minimal) period $T>0$ 
gives a 1-periodic orbit of~$X_H$ of the form
\begin{equation}
\label{perorb}
t \mapsto \bigl(\rho,x(Tt) \bigr),
\end{equation}
where $\rho>1$ is the unique solution of the equation $h'(\rho)=T$. Indeed, the existence of such a solution~$\rho$ follows from~\eqref{condh} and its uniqueness from the strict monotonicity 
of~$h'$. An easy computation shows that this 1-periodic orbit has action
\[
\ca_H = \rho \, h'(\rho) - h(\rho).
\]
The derivative of this function with respect to $\rho$ is $\rho \, h''(\rho)$, 
and we deduce that $\ca_H$ is positive at these 1-periodic solutions and increases monotonically with~$\rho$. We conclude that there is a one-to-one correspondence between the periodic orbits 
of~$R_{\alpha}$ on~$\partial M$ of arbitrary period and the 1-periodic orbits of~$X_H$ 
in the set $\widehat{M} \setminus M$, and the function Period $\mapsto$ Action associated to this correspondence is a strictly increasing positive function.

The Hamiltonian $H$ does not quite fit in the class of Hamiltonians for which we defined the Floer complex in Section~\ref{s:Ham}. The first missing property is the non-degeneracy assumption: The 1-periodic orbits belonging to the continuum of constant solutions in~$M$ are certainly degenerate, and also all the non-constant ones, since $H$ is autonomous. 
This issue can be fixed by perturbing~$H$ by a small time-periodic function. Here we shall simply ignore this issue and pretend that we can work with the Hamiltonian~$H$ defined above. We just mention that if one assumes the Reeb vector field~$R_{\alpha}$ to be non-degenerate, 
meaning that $1$ is not an eigenvalue of the restriction of the linearisation of the autonomous Reeb flow to a transverse section, then one can arrange this perturbation in such a way that every $T$-periodic orbit of~$R_{\alpha}$ produces exactly two 1-periodic orbits of the Hamiltonian system on $\widehat{M}$ which are close to~\eqref{perorb}.

Furthermore, the underlying symplectic manifold $\widehat{M}$ in this case is not compact. 
However, the above computation of the space of 1-periodic orbits of~$X_H$ shows that the critical set 
of~$\ca_H$ intersected with an arbitrary sublevel $\{\ca_H \leqslant a\}$ is compact. 
Indeed, an upper bound for the action of a periodic orbit of the form~\eqref{perorb} 
gives an upper bound for the period of the corresponding closed Reeb orbit~$x$. 
So, after perturbing~$H$ in order to make it non-degenerate, we will have finitely many 
1-periodic orbits in any action sublevel.

What remains to be checked is that for every pair of critical points~$x$ and~$y$ 
of~$\ca_H$ on~$C^{\infty}(\TT,\widehat{M})$, the images of the maps in the moduli spaces 
$\widehat{\mathcal{M}}_{H,J}(x,y)$ are contained in some compact subset of~$\widehat{M}$. 
Indeed, once this is guaranteed, the compactness theory works as in the case of 
a closed manifold~$M$, because $\omega$ certainly vanishes on $\pi_2(\widehat{M})$, 
since it is exact. 

This compactness result can be achieved by a proper choice of behaviour of the almost complex structure~$J$ 
on~$\widehat{M} \setminus M$: 
The condition is that the family of time-dependent $\omega$-compatible 
almost complex structures~$J_t$ satisfies
\begin{equation} \label{e:rJ}
d\rho \circ J_t = \lambda,
\end{equation}
for $\rho \geqslant \rho_0$, where $\rho$ is again the coordinate of the first component in the product $(1,+\infty) \times \partial M \cong \widehat{M} \setminus M$ and $\rho_0$ is some positive number. Indeed, this assumption implies that for every solution~$u$ 
of the Floer equation~\eqref{floer} the real valued function $\rho = \rho \circ u$ satisfies
\[
\Delta \rho = \left\| \frac{\partial u}{\partial s} \right\|_{J_t}^2 + \rho \, h''(\rho) \, \frac{\partial \rho}{\partial s}
\]
in the region of $\RR \times \TT$ in which $\rho \geqslant \rho_0$. Then, the maximum principle implies that this function cannot have interior maxima in this region, and hence $\rho \circ u$ is everywhere not larger than the maximum of~$\rho_0$ and the values of~$\rho$ at the two
asymptotic periodic orbits (in case that they lie in~$\widehat M \setminus M$). 
This is the point where the crucial assumption that the boundary of~$M$ is of contact type 
is used. 

Now we can define the Floer complex $\{\CF_k(H),\partial_k\}$ exactly as in the case of a 
closed symplectic manifold, with the only difference that $\CF_k(H)$ may be infinitely generated for some~$k$: Even after perturbing~$H$ in order to achieve non-degeneracy, 
it may happen that~$X_H$ has a sequence of 1-periodic orbits with the same Conley--Zehnder index and action going to~$+\infty$. The boundary $\partial x$ of every 1-periodic orbit~$x$ is nevertheless well-defined, because $\mathcal{M}_{H,J}(x,y)$ can be non-empty only for the finitely many 1-periodic orbits~$y$ with action less than~$\ca_H(x)$.
 
The homology of the Floer complex of~$H$ is independent of the choices made. It is called the 
{\em symplectic homology}\/ of the Liouville domain~$M$ and is denoted by~$\SH_*(M)$. 
This homology is rather stable: For instance, it does not change if we smoothly vary the symplectic form 
in such a way that~$\partial M$ remains of contact type along this deformation.

Moreover, symplectic homology comes with an important homomorphism
\[
c_* \colon \H_{*+n}(M,\partial M) \rightarrow \SH_*(M),
\]
whose definition we now briefly discuss. Recall that the critical set of the unperturbed action functional~$\ca_H$ consists of the constant loops in~$M$, which have action zero, 
and non-constant 1-periodic orbits, which are produced by the closed orbits of~$R_{\alpha}$ 
and have positive action bounded away from zero, say larger than~$2\eps$. When we perturb the Hamiltonian~$H$ in order to achieve non-degeneracy, we may assume the new Hamiltonian~$\widetilde{H}$ to be autonomous, $C^2$-small on~$M$, 
where we are actually perturbing the zero function, negative on the interior of~$M$ and vanishing on~$\partial M$. If the perturbation is small enough, all the new critical points of $\mathcal{A}_{\widetilde{H}}$ with action less than~$\eps$ 
are just critical points of~$\widetilde{H}$ in~$M$. By the same argument which leads to the computation of 
Floer homology for $C^2$-small autonomous Hamiltonians on closed manifolds, the subcomplex 
of $\CF_*(\widetilde{H})$ generated by critical points of $\mathcal{A}_{\widetilde{H}}$ with action 
less than~$\eps$ is just the Morse complex of $-\widetilde{H}$, shifted by~$n$. 
By the properties of~$\widetilde{H}$, the homology of this Morse complex is isomorphic to the relative singular homology $\H_*(M,\partial M)$, since the negative gradient of~$-\widetilde{H}$ is pointing in the outward direction on~$\partial M$. 
The natural homomorphism induced at the homology level by an inclusion of chain complexes
gives us the map~$c_*$.

Notice that the above homomorphism is an isomorphism if and only if the above subcomplex coincides with the whole complex $\CF_*(\widetilde{H})$, and this happens if and only if 
$R_{\alpha}$ has no closed orbits on~$\partial M$. Therefore, whenever $c_*$ is not an isomorphism, the Weinstein conjecture holds true:

\begin{theorem}
\label{vit}
Let $(M,\lambda)$ be a Liouville domain such that $c_*$ is not an isomorphism. 
Then the Reeb flow of $\lambda|_{\partial M}$ on $\partial M$ has closed orbits.
\end{theorem}

It is not difficult to show that the symplectic homology of a ball in~$\RR^{2n}$ vanishes, 
and from the stability properties of symplectic homology one deduces that the same is true 
for any starshaped domain in~$\RR^{2n}$. In particular, $c_*$ is the zero map when $M$ is such a domain, 
and since $\H_*(M,\partial M) \neq 0$, the Weinstein conjecture holds true for starshaped domains 
in~$\RR^{2n}$, a result which as already mentioned was first proved by Rabinowitz~\cite{Rab78} by variational methods. 

Actually, the above argument is quite flexible and permits to prove the Weinstein conjecture in greater generality, see \cite[Theorem 4.1]{vit99}. For instance, symplectic homology vanishes 
for all {\em subcritical}\/ Stein domains, that is, $2n$-dimensional Stein domains whose defining plurisubharmonic function has only critical points of Morse index strictly less than~$n$ 
(while an arbitrary plurisubharmonic function might have critical points of Morse index 
up to~$n$). From this fact one can deduce that the Weinstein conjecture holds true for any contact type hypersurface bounding a compact set in a subcritical Stein manifold. 

\subsection{Viterbo functoriality} 
Another important feature of symplectic homology for Liouville domains is that it is functorial with respect to Liouville embeddings. A Liouville embedding of the Liouville domain 
$(M',\lambda')$ into the Liouville domain $(M,\lambda)$ is 
an embedding $\imath \colon M' \hookrightarrow M$ of codimension zero such that
\[
\imath^*(\lambda) = c \lambda' + df
\]
for a positive constant $c$ and a smooth function $f \colon M' \rightarrow \RR$. 
A Liouville embedding $\imath \colon M' \hookrightarrow M$ induces a homomorphism
\[
\SH_*(\imath) \colon \SH_*(M) \rightarrow \SH_*(M')
\]
which is functorial with respect to composition of embeddings and 
invariant under isotopies of Liouville embeddings.
This homomorphism should be considered as a {\em transfer}\/ or {\em shriek}\/ map, 
analogous to the homomorphism which arises in the Gysin sequence of a fiber bundle. 
It is called {\em Viterbo transfer map}. The Viterbo transfer map behaves well with respect to the map~$c_*$, meaning that the diagram
\[
\begin{CD}
\H_{*+n}(M,\partial M) @>{\imath !}>> \H_{*+n}(M',\partial M') \\ @V{c_*}VV @VV{c_*}V \\ 
\SH_*(M) @>{\SH_*(\imath)\,}>> \SH_*(M') 
\end{CD}
\]
commutes. Here, the shriek map $\imath !$ is the composition
\[
\H_{*+n}(M,\partial M) \cong \H^{n-*}(M) \stackrel{\imath^*}{\longrightarrow} \H^{n-*}(M') \cong \H_{*+n}(M',\partial M'),
\]
where the two isomorphisms are given by Poincar\'e duality. The construction of the map 
$\SH_*(\imath)$ uses suitable 1-parameter families of Hamiltonians and an algebraic counting 
of the solutions of the corresponding $s$-dependent Floer equation.
 
\subsection{The symplectic homology of cotangent disc bundles}  \label{s:disc}
So far we have seen only examples of Liouville domains with vanishing symplectic homology. Now we wish to compute the symplectic homology of a fiberwise starshaped domain in the cotangent bundle of a closed manifold~$Q$ and show that it is far from being zero. 
By the invariance properties of symplectic homology, all these domains have the same symplectic homology and it is enough to consider the unit cotangent disc bundle~$D^*Q$ which is induced by a Riemannian metric on~$Q$. The relevant Hamiltonian now is a convex radial function which is identically zero on~$D^*Q$ and grows superlinearly outside of~$D^*Q$. After a time-dependent perturbation, which is in any case necessary to have non-degeneracy of 1-periodic orbits, we may assume that $H \in C^{\infty}(\TT\times T^* Q,\RR)$ is fiberwise strictly convex and superlinear. This is an interesting class of Hamiltonians on cotangent bundles, 
which are called {\em Tonelli Hamiltonians}\/ and behave well with respect to the Legendre transform. It is actually useful to strengthen the fiberwise convexity by asking that
 \[
d_{pp}H(t,q,p)[\xi,\xi] \geqslant \eps |\xi|^2
\]
for some $\eps>0$ and that $H$ does not grow more than quadratically in~$p$. Under these assumptions, the Legendre dual Lagrangian $L \colon \TT \times TQ \rightarrow \RR$, which is defined by
\begin{equation}  \label{e:Leg}
L(t,q,v) := \max_{p\in T_q^* Q} \bigl( p(v) - H(t,q,p) \bigr),
\end{equation}
has the same qualitative properties. The Legendre transform
\[
(t,q,p) \mapsto \bigl(t,q,d_p H(t,q,p)\bigr)
\]
induces a one-to-one correspondence between the 1-periodic orbits of~$X_H$ and those of 
the Euler--Lagrange equations of the Lagrangian action functional
\[
\mathcal{S}_L(q) = \int_{\TT} L(t,q,\dot q)\, dt, \qquad q\in C^{\infty}(\TT,Q).
\]
In other words, there is a one-to-one correspondence between the critical points of~$\ca_H$ 
and~$\mathcal{S}_L$. The critical points of~$\mathcal{S}_L$ have finite Morse index, 
which agrees with the Conley--Zehnder index of the corresponding Hamiltonian orbits. 
Moreover, $\mathcal{S}_L$ is bounded from below and satisfies the Palais--Smale condition on the Hilbert manifold $H^1(\TT,Q)$ of loops of Sobolev class~$H^1$ on~$Q$. 
Therefore, $\mathcal{S}_L$ has a well-defined Morse complex $\{\CM_*(\mathcal{S}_L),\partial\}$, whose homo\-lo\-gy is the singular homo\-lo\-gy of~$H^1(\TT,Q)$, or equivalently the singular homo\-lo\-gy of the homotopically equivalent free loop space $C^{\infty}(\TT,Q)$. The natural identifications
\[
\CM_k(\mathcal{S}_L) \cong \CF_k(H)
\]
make it plausible that the Floer homology of $(T^*Q,H)$ is also isomorphic to the singular homology of $C^{\infty}(\TT,Q)$. But this is by no means obvious, since the boundary operators in the Floer homology of~$(T^*Q,H)$ and in the Morse homology of~$\mathcal{S}_L$ on~$H^1(\TT,Q)$ have quite different definitions.
Nevertheless, it is possible to prove that the Floer complex of~$(T^*Q,H)$ and the Morse complex of~$\mathcal{S}_L$ on~$H^1(\TT,Q)$ are isomorphic.

The proof of this fact is based on a better understanding of the relationship between the Hamiltonian and the Lagrangian action functionals. In order to compare these two functionals, it is convenient to read the former on the loop space of~$TQ$ by means of the global diffeomorphism
\[
\mathcal{L} \colon C^{\infty}(\TT,TQ) \rightarrow C^{\infty}(\TT,T^*Q), \qquad (q,v) \mapsto \bigl(q,d_vL(\cdot,q,\dot q + v)\bigr),
\]
which is induced by the Legendre transform. Here $q$ denotes a loop in~$Q$ and $v$ is a loop of tangent vectors with $v(t) \in T_{q(t)} Q$ for all $t \in \TT$. Legendre duality and Taylor's formula imply the identity
\[
\ca_H(\mathcal{L}(q,v)) = \mathcal{S}_L(q) - \mathcal{U}(q,v),
\]
where $\mathcal{U}$ is the following functional on $C^{\infty}(\TT,TQ)$:
\[
\mathcal{U}(q,v) = \int_{\TT} \int_0^1 s\, d_{vv}L(t,q,\dot{q}+sv)[v,v]\, ds\, dt.
\]
The functional $\mathcal{U}$ vanishes on the loops of the form $(q,0)$, that is, 
loops taking values in the zero section of~$TQ$. Everywhere else it is positive and grows quadratically in~$v$. Therefore, $\mathcal{U}$ behaves like the functional 
$(q,v) \mapsto \|v\|^2_{L^2(\TT)}$.

We may think of $C^{\infty}(\TT,TQ)$ as a vector bundle over $C^{\infty}(\TT,Q)$: 
\[
\Pi \colon C^{\infty}(\TT,TQ) \rightarrow C^{\infty}(\TT,Q), \qquad (q,v)\mapsto q.
\]
In this picture, $\ca_H \circ \mathcal{L}$ agrees with $\mathcal{S}_L$ on the zero section of this bundle, and on each fiber it restricts to a concave functional achieving its maximum at zero. 
It is not difficult to define a negative gradient flow for $\ca_H \circ \mathcal{L}$ which leaves the zero section invariant, agrees on the zero section with the negative gradient flow 
of~$\mathcal{S}_L$, and commutes with the bundle projection. Actually, for this we need to replace the space $C^{\infty}(\TT,Q)$ by its Sobolev completion~$H^1(\TT,Q)$, but let us ignore this issue and pretend that we can work with spaces of smooth loops. 
All the critical points of $\ca_H \circ \mathcal{L}$ are contained in the zero section, their stable manifolds are also contained there and agree with the (infinite dimensional) stable manifolds of $-\nabla \mathcal{S}_L$, while their unstable manifolds are the inverse images 
by~$\Pi$ of the (finite dimensional) unstable manifolds of $-\nabla \mathcal{S}_L$. 
This implies that the Morse complex of this negative gradient flow is 
well-defined and agrees with the Morse complex of~$\mathcal{S}_L$, since intersections between stable and unstable manifolds occur only in the zero section.

Now we can push this negative gradient flow forward to $C^{\infty}(\TT,T^*Q)$, and we end up in the following situation: The functional~$\ca_H$ on~$C^{\infty}(\TT,T^*Q)$ has two negative gradient equations. The first is the one just described: it gives a genuine flow and the corresponding Morse complex is just the Morse complex of $\mathcal{S}_L$. We may call this flow 
\[
\psi^s \colon C^{\infty}(\TT,T^*Q) \rightarrow C^{\infty}(\TT,T^*Q).
\]
The second one is the Floer negative gradient equation, which induces no flow but defines nevertheless the Floer complex of~$H$.

When a functional has two different gradient equations, it is to be expected that their Morse complexes (or Floer complexes, if these equations have no well-defined flow) are isomorphic. 
A standard way of seeing this is to build a chain map between them which, given two critical points~$x$ and~$y$, counts the continuous paths~$u(s)$, $s\in \RR$, from $x$ to~$y$ 
which on $(-\infty,0)$ satisfy the first negative gradient equation and on $(0,+\infty)$ the second one. Indeed, this construction defines a chain map~$\Phi$ between the two chain complexes which is of the form
\[
\Phi x = x + \sum_y \nu(x,y) \, y,
\]
where the critical points $y$ appearing in the sum have function level lower than the function level of~$x$. Homomorphisms of this form are readily seen to be isomorphisms.

\begin{figure}[h]   
 \begin{center}
  \psfrag{x}{\footnotesize $x$}
  \psfrag{q}{\footnotesize $q$}
  \psfrag{y}{\footnotesize $y$}
  \psfrag{u}{\footnotesize $u$}
  \psfrag{Q}{\footnotesize $Q$}
  \psfrag{p}{\footnotesize $\Pi$}
  \psfrag{pu}{\footnotesize $\Pi (u(0))$}
  \psfrag{TQ}{\footnotesize $T_q^*Q$}
  \psfrag{u0}{\footnotesize $u(0)$}
  \psfrag{-A}{\footnotesize $-\nabla \ca_H$}
  \psfrag{-S}{\footnotesize $-\nabla \cs_L$}
  \psfrag{c}{\footnotesize $C^\infty (\TT, Q)$}
  \psfrag{C}{\footnotesize $C^\infty (\TT, T_q^*Q)$}
  \leavevmode\includegraphics{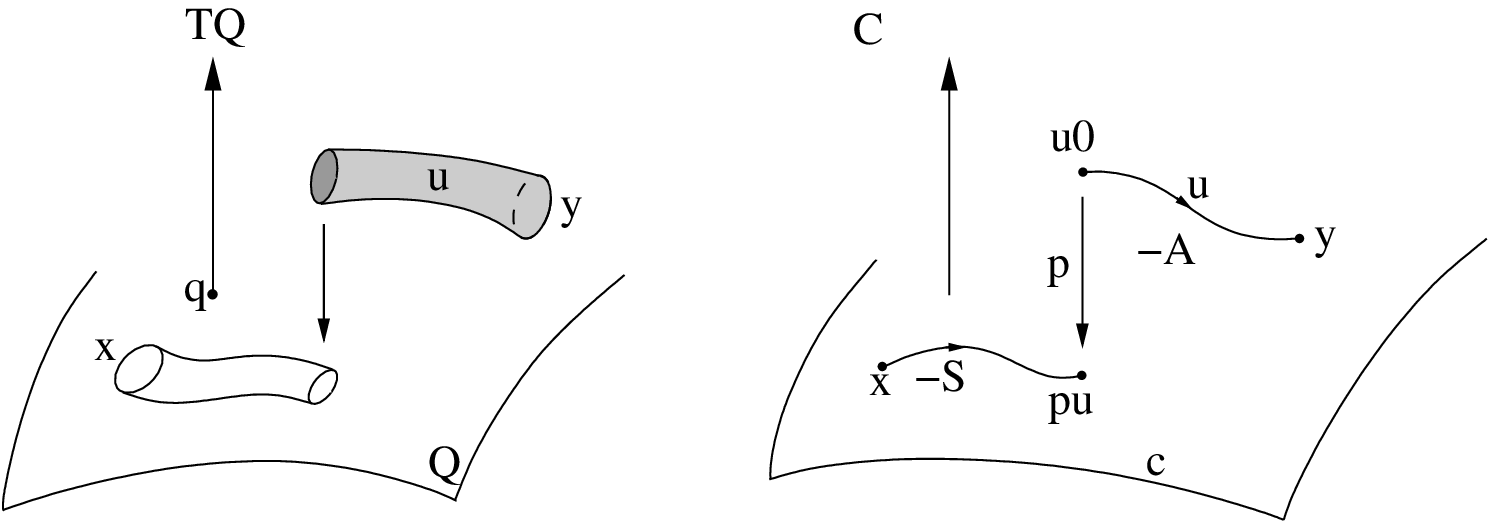}
 \end{center}
 \caption{An element of $\cm_\Phi (x,y)$, seen in $T^*Q$ and in $C^\infty (\TT, T^*Q)$}  \label{fig.AS}
\end{figure}

Let us see what this procedure gives us when the first negative gradient equation is the one inducing the flow~$\psi^s$ and the second is the Floer equation. Let $x$ and~$y$ be critical points of~$\ca_H$. The above recipe tells us that we should consider the space of paths 
$u \colon [0,+\infty) \rightarrow C^{\infty}(\TT,T^*Q)$ which solve the Floer equation, 
converge to~$y$ for $s \rightarrow +\infty$ and for $s=0$ belong to the unstable manifold of~$x$ with respect to the flow~$\psi^s$. By construction, this is the image under~$\mathcal{L}$ 
of the sets
\[
\Pi^{-1}\bigl(W^{\u}(x;-\nabla \mathcal{S}_L)\bigr).
\]
By the form of $\mathcal{L}$, this is just the inverse image of $W^{\u}(x;-\nabla \mathcal{S}_L)$  under the bundle projection $C^{\infty}(\TT,T^*Q) \mapsto C^{\infty}(\TT,Q)$. 
In other words, the relevant moduli spaces for the construction of the isomorphism
\begin{equation} \label{cm:Phi}
\Phi \colon \CM_*(\mathcal{S}_L) \rightarrow \CF_*(H)
\end{equation}
are the sets
\[
\begin{split}
\mathcal{M}_{\Phi}(x,y):= \{u \colon [0,+\infty) \times \TT \rightarrow T^*Q \,\mid\; & 
\mbox{$u$ solves the Floer equation for $H$}, \\ & 
\lim_{s \rightarrow +\infty} u(s,\cdot)=y,  \; \pi \circ u(0,\cdot) \in 
                                                  W^{\u}(x;-\nabla \mathcal{S}_L)\},
\end{split}
\]
where $\pi \colon T^* Q \rightarrow Q$ is the bundle projection. Actually, at this point one can forget about the heuristic arguments which brought us to the above definition of 
the moduli space~$\mathcal{M}_{\Phi}(x,y)$ and, after a careful analysis of their differentiable and compactness properties, use them to define an isomorphism~$\Phi$ from the Morse complex 
of~$\mathcal{S}_L$ to the Floer complex of~$H$. This implies that the symplectic homology of the cotangent disc bundle~$D^* Q$ is isomorphic to the singular homology of the free loop space of~$Q$:
\begin{equation}
\label{iso}
\SH_*(D^*Q) \cong \H_*(C^{\infty}(\TT,Q)).
\end{equation}
By looking at low values of the action, one can prove that the map 
\[
c_* \colon \H_{*+n}(D^*Q,\partial D^*Q) \longrightarrow \SH_*(D^*Q) \cong \H_*(C^{\infty}(\TT,Q))
\]
is the injective homomorphism given by the composition
\[
\H_{*+n}(D^*Q,\partial D^*Q) \cong \H_{*+n}(T^*Q,\partial T^*Q\setminus Q) \longrightarrow \H_*(Q) \longrightarrow \H_*(C^{\infty}(\TT,Q)),
\]
where the first arrow is the Thom isomorphism
(the cap product with the Thom class of the vector bundle~$T^*Q$) 
and the second one is induced by the inclusion of~$Q$ into its loop space as the set of constant loops. The fact that this inclusion has a left inverse implies that the homomorphism which it induces in homology is injective.

Actually, there is an issue with orientations in Floer theory which makes the existence of the 
isomorphism~\eqref{iso} true only if one works with $\ZZ_2$ coefficients, or if one assumes that the second Stiefel--Whitney class of~$Q$ vanishes on tori. In general, the symplectic homology of a starshaped domain in~$T^* Q$ is isomorphic to the singular homology of the free loop space 
of~$Q$ with a suitable local system of coefficients. 

\subsection{An obstruction to the existence of Lagrangian submanifolds}

In the same paper \cite{Gro85} where he proved the non-squeezing theorem, Gromov used $J$-holomorphic disks to prove 
that $(\RR^{2n},\omega_0=d\lambda_0)$ admits no compact exact Lagrangian submanifold. 
Using what we have seen of symplectic homology we can easily prove the following generalisation of Gromov's result:

\begin{theorem}
A Liouville domain with vanishing symplectic homology admits no compact exact Lagrangian submanifold.
\end{theorem}

This result applies, for instance, to subcritical Stein domains, and hence to subcritical Stein manifolds. The proof goes as follows: If $L$ is a compact Lagrangian submanifold of the Liouville domain $(M,\lambda)$, 
then a theorem of Weinstein tells us that there is 
a symplectic embedding~$\imath$ of a (small) cotangent disc bundle~$D^*L$ onto a neighbourhood of~$L$, 
mapping the zero section onto~$L$. By the exactness of~$L$, the 1-forms
$\imath^*(\lambda)$ and $\lambda_{\can}$ on~$D^*L$ differ by a differential, so $\imath$ is a Liouville embedding. The induced Viterbo transfer map determines the commutative diagram
\begin{eqnarray} \label{e:diaVit}
\begin{CD}
\H_{*+n}(M,\partial M) @>\imath !>> \H_{*+n}(D^*L,\partial D^*L) \\ @V{c_*}VV @VV{c_*}V \\ \SH_*(M) @>{\SH_*(\imath)\,}>> \SH_*(D^*L). 
\end{CD}
\end{eqnarray}

Here and in what follows we use $\ZZ_2$ coefficients. 
Since we assume that the symplectic homology of~$M$ vanishes, 
the composition of the left vertical map with the bottom horizontal one is zero. 
However, $\imath !$ is an isomorphism in degree~$*=n$, 
and the right vertical map does not vanish in that degree, since, as we have seen, 
it corresponds to the injective homomorphism
\[
\H_{2n}(D^*L,\partial D^*L) \cong \H_n(L) \longrightarrow \H_n(C^{\infty}(\TT,L)) \cong \SH_n(D^*L),
\] 
which is induced by the inclusion $L \hookrightarrow C^{\infty}(\TT,L)$.
This contradiction proves the theorem. Again, the fact that symplectic homology 
can be defined in a more general setting implies more general versions of 
the above theorem, see e.g.\ \cite[Theorem 4.3]{vit99}.

\subsection{Obstructions to symplectic embeddings} \label{s:emb}
In order to extract from symplectic homology invariants that are monotone 
with respect to Liouville embeddings, and thus have the potential to detect 
obstructions to symplectic embeddings, 
we filter symplectic homology by action:
Given a real number $a>0$ and a non-degenerate Hamiltonian $H \colon \TT \times M \to \RR$
on the Liouville domain $M$ let $\Crit^a \ca_H = \{ x \in \Crit \ca_H \mid \ca_H(x) \leqslant a \}$
and let $\CF^a (H)$ be the vector space generated by~$\Crit^a \ca_H$.
Since the action~$\ca_H$ does not increase along solutions~$u$ of its 
negative gradient equation~\eqref{floer}, 
the boundary operator~$\pp$ maps $\CF^a (H)$ to itself.
We thus obtain for every~$a$ the filtered Floer homology groups $\HF^a (H)$,
and for $H$ close to zero on~$M$ and a sufficiently steep function of~$\rho$ on the
symplectisation part we obtain the symplectic homology groups~$\SH_*^a(M)$.
They can be seen as the Floer homology built from the critical points of~$H$ in~$M$ and pairs of
Reeb orbits on~$\pp M$ of action $\leqslant a$.
The maps $c_*$ and the Viterbo transfer maps $\SH_*(\iota)$ are compatible with
the action filtration. In particular, the commutative diagram~\eqref{e:diaVit}
for Liouville embeddings $M' \hookrightarrow M$
can be refined to the diagrams
\[
\begin{CD}
\H_{*+n}(M,\partial M) @>{\imath !}>> \H_{*+n}(M',\partial M') \\ @V{c_*}VV @VV{c_*}V \\ 
\SH_*^a(M) @>{\SH_*(\imath)\,}>> \SH_*^a(M') 
\end{CD}
\]
Now define $\gg (M) \in [0,+\infty]$ by 
\begin{equation} \label{def:gg}
\gg (M) \,=\, \inf 
\left\{ a>0 \mid c_n \colon \H_{2n} (M,\pp M) \to \SH_n^a (M) \mbox{ vanishes} \right\} ,
\end{equation}
where the infimum of the empty set is $+\infty$. Notice that $\gamma(M)$ is finite when the natural homomorphism $c_*$ vanishes, and in particular when $\SH_*(M)=0$.
Since $\iota!$ is an isomorphism for $*=n$ and since $\gg (M)$ is a symplectic invariant,
the above diagram with $*=n$ shows that $\gg$ is monotone with respect to Liouville embeddings.

Recall from the construction of the maps $c_*$ that they are isomorphisms
if $a$ is smaller than the action of the shortest Reeb orbit on~$\pp M$.
Hence $\gg (M)$ is not smaller than this number.
As we shall see shortly,  
for ellipsoids the converse is also true.

To better see the meaning of the number $\gg (M)$, 
take as $M$ a starshaped domain $U \subset \RR^{2n}$,
and in the construction of $\SH_*^a (U)$ take Hamiltonians that have exactly one
critical point~$x_0$ (the minimum) at~$0$.
This critical point generates $\H_{2n} (U,\pp U)$ and his action is close to zero, 
while all other periodic orbits of~$H$ have larger action. 
Hence $\gg (U)$ is the smallest value~$a$
such that $x_0$ is the boundary of a linear combination of orbits of~$H$ 
of Conley--Zehnder index~$n+1$ all of whose actions are $\leqslant a$.

Now take as $U$ an ellipsoid $\E (a_1, \dots, a_n)$ with $a_1 < \dots < a_n$ rationally independent, 
see~\eqref{ellipsoid}.
One can then choose the relevant Hamiltonians~$H$ such that there exists exactly 
one periodic orbit~$x$ of index~$n+1$, whose action is~$\approx a_1$, and such that
$\pp x = x_0$, see~\cite{fhw94}. 
(That $\pp x = x_0$ also follows from the general fact that $\SH_*(U)$ 
vanishes for starshaped domains in~$\RR^{2n}$.)
Hence 
\begin{equation} \label{e:E}
\gg (\E (a_1, \dots, a_n)) \,=\, a_1 .
\end{equation}

This yields a proof of Gromov's non-squeezing theorem~\ref{t:Gromov}:
Assume that there exists a symplectic embedding $\varphi \colon \B^{2n}(a) \to \Z^{2n}(b)$
and fix $\gve >0$. We then find irrational ellipsoids $\E := \E(a_1, \dots, a_n)$ 
and $\E' := \E (b_1, \dots, b_n)$ such that
$a-\gve < a_1 < \dots < a_n < a$ and $b=b_1 < b_2 < \dots < b_n$ and
$\varphi (\E) \subset \E'$.
By~\eqref{e:E} and by the monotonicity of~$\gamma$, 
$a_1 = \gg (\E) \leqslant  \gg (\E') = b$, and so $a \leqslant b+\gve$.
Since this holds for all $\gve >0$, the claim follows.

\subsection{Rabinowitz--Floer homology}
\label{s:RF}

We consider a compact hypersurface $\Sigma$ bounding a compact domain~$U$ in an exact symplectic manifold 
$(M,\lambda)$.
Recall from Section~\ref{s:SH} that one can find closed characteristics
on~$\Sigma$ by taking Hamiltonian functions that are very flat on~$U$ and very steep on the symplectisation
part $\widehat U \setminus U$.
Rabinowitz--Floer homology is a variant of Floer homology that finds closed characteristics
on~$\Sigma$ in a different way, namely by using a Lagrangian multiplier.

Critical points of a function $f \colon \RR^n \to \RR$ restricted to a constraint hypersurface $\{ x \in \RR^n \mid g(x) =0\}$
can be found by looking at the function
$$
F(x,\eta) \,=\, f(x) + \eta \,g(x)
$$
on $\RR^n \times \RR$, whose critical points are those $(x,\eta)$ with $g(x) =0$
and $df(x) + \eta \,dg(x) =0$.
Given an autonomous Hamiltonian function 
$H \colon M \to \RR$ with $H^{-1}(0) =\Sigma$
we thus look at the action functional $\ca_H \colon C^\infty (\TT,M) \times \RR$ defined by
\begin{equation}  \label{e:unperturbed}
\ca_H (x,\eta) \,=\, \int_{\TT} x^*\lambda - \eta \int_{\TT} H (x(t)) \,dt .
\end{equation}
Its critical points $(x,\eta)$ are the solutions of the problem
$$
\dot x(t) = \eta \,X_H(x(t)), \qquad 0 =  \int_{\TT} H(x(t)) \, dt ,
$$
and therefore correspond to
%
closed orbits of~$X_H$ on the fixed energy surface 
$\Sigma = H^{-1}(0)$ of arbitrary period $|\eta| \geqslant 0$.

The Floer homology of this functional is called Rabinowitz--Floer homology $\RFH (\Sigma, M)$.
The new difficulty in its construction is to establish a uniform bound on the Lagrangian multiplier~$\eta$ 
along the solutions of the Floer equation between two critical points.
A particularly interesting application of~$\RFH$ are lower bounds for {\it leafwise intersection points}:
Let $\gg$ be a characteristic on $\Sigma = H^{-1}(0)$, and 
assume that an integral curve $x(t)$ of the Hamiltonian vector field~$X_H$ 
calmly moves along~$\gg$, but that at time $t=0$
an ``earthquake'' happens, during time 1, given by a time-dependent Hamiltonian function~$F$.
Under this new force, our curve $x(t)$ will in general leave its usual trajectory~$\gg$ 
and even its energy surface~$\Sigma$.
Is it possible that nonetheless $x(1)$ is back on the curve~$\gg$, 
as if the earthquake had not happened?
Such a point $x(1)$ is called a leafwise intersection.
Leafwise intersection points are the critical points of the action functional obtained by adding to~\eqref{e:unperturbed}
the term $-\int_{\TT} F_t(x(t)) \,dt$, and so the Rabinowitz--Floer homology of this functional leads to 
lower bounds of the number of leafwise intersection points, 
at least if the earthquake is not too strong.
Denote by $\cuplength (\Sigma,M)$ the maximal length of a cup product of elements of positive degree 
in~$\H^* (M;\ZZ_2)$ whose restriction to~$\Sigma$ does not vanish.

\begin{theorem} \label{t:RF}
Let $\Sigma$ be a compact and bounding hypersurface of contact type
in the exact symplectic manifold~$(M,\lambda)$.
Denote by $\frak p >0$ the smallest period of a closed Reeb orbit on $(\Sigma, \lambda |_{\Sigma})$
that is contractible in~$M$.
Let $F \colon M \times [0,1] \to \RR$ be a compactly supported function such that
$$ 
\int_0^1 \left( \max_{x \in M} F(x,t) - \min_{x \in M} F(x,t) \right) dt
\,<\, \frak{p} .
$$
Then the number of leafwise intersection points on~$\Sigma$ is 
\begin{itemize}
\item[(i)]
at least $1 + \cuplength (\Sigma,M)$;

\s
\item[(ii)]
at least $\dim \H (\Sigma;\ZZ_2)$ for generic~$F$.
\end{itemize}
\end{theorem}

The example of a circle in the plane, that can be displaced by Hamiltonian diffeomorphisms,
shows that the smallness assumption on~$F$ cannot be omitted.

\subsection*{Historical notes and bibliography} 
Symplectic homology was first introduced by Floer and Hofer in~\cite{fh94} and further developed, also jointly with Cieliebak and Wysocki, in \cite{fhw94,cfh95,cfhw96}. 
This first version of symplectic homology applies to open subsets of $(\R^{2n},\omega_0)$, 
or of more general symplectic manifolds with contact type boundary. 
Among its applications are the symplectic classification of ellipsoids and polydiscs:
$\E (a_1, \dots, a_n)$ and $\E (b_1, \dots, b_n)$ are symplectomorphic only if 
$\{ a_1, \dots, a_n \} = \{ b_1, \dots, b_n \}$, and similarly for polydiscs,
see~\cite{fhw94}.

The more qualitative version of symplectic homology which we discussed here was introduced by Viterbo in~\cite{vit99}. 
An alternative approach to time-dependent perturbations has been developed by Bourgeois and Oancea in~\cite{bo09a}. 
The survey~\cite{oan04} describes both versions of symplectic homology, 
while the survey~\cite{sei08} focuses on Viterbo's one. 

The analysis of the behaviour of symplectic homology under handle attachment has been carried out by Cieliebak in~\cite{cie02}, 
who in particular showed that the symplectic homology of subcritical Stein domains vanishes, see also~\cite{bc01}. 
The theory of Lefschetz fibrations can be extended to the symplectic setting, and a symplectic Lefschetz fibration on the completion of a Liouville domain~$M$ 
often allows one to compute the symplectic homology of~$M$, see \cite{eli97,sei08,mcl12}.

Bourgeois and Oancea~\cite{bo09b} have related the symplectic homology of a compact symplectic manifold with contact type boundary to the linearized contact homology 
of its boundary via an exact sequence, cf.\ Section~\ref{s:ch}.
Symplectic homology has an $S^1$-equivariant version, due to the fact that it is the Floer homology of an autonomous Hamiltonian. 
See \cite[\S 5]{vit99}, \cite{sei08}, \cite{bo13}. 
In the latter paper the positive part of the $S^1$-equivariant symplectic homology of~$M$, 
which is obtained from the Floer complex by modding out the subcomplex generated by critical points of action 
less than~$\eps$ (see the discussion around the definition of~$c_*$ above), 
is shown to be isomorphic to the linearized contact homology of~$\partial M$.

Seidel and Smith have used symplectic homology to prove the existence of exotic symplectic structures 
convex at infinity on~$\R^{4n}$ for all $n \geqslant 2$, 
see~\cite{ss05}. See also \cite{mcl09} and~\cite{sei11} for other results in this direction. 

The isomorphism between the symplectic homology of the cotangent disc bundle of a 
closed manifold~$Q$ and the singular homology of the free loop space of~$Q$ 
was first proved by Viterbo in~\cite{vit03} by using generating functions. 
A different proof involving the heat flow for loops on a Riemannian manifold was found by Salamon and Weber in~\cite{sw06}. 
The proof we sketched here is due to the first author and Schwarz, 
see \cite{as06} and~\cite{as15}. 
The orientation issue which requires the use of local coefficients when the second Stiefel--Whitney class does not vanish on tori had been overlooked in \cite{vit03,sw06,as06} 
and was discovered by Kragh in~\cite{kra07}, see also~\cite{kra13}, 
and corrected in~\cite{as14,as15}. 
See also~\cite{abo15} for another approach to this isomorphism and~\cite{coh10} for an extension of these ideas towards homotopy. 
This isomorphism is also the starting point of a fruitful interaction between symplectic geometry and string topology, see~\cite{as10,cl09}.
For instance, the pair of pants product described at the end of Section \ref{s:Ham}
can also be defined on the symplectic homology of~$D^*Q$,
and under the isomorphism~\eqref{iso} it corresponds to the Chas--Sullivan product on 
the homology of the free loop space of~$Q$, see~\cite{as10}. 

The spectral invariant $\gg (M)$ in \eqref{def:gg} goes back to~\cite{fhw94}
and in this form was defined in~\cite{vit99}.
This is just one of many similar spectral invariants that can be
constructed by Floer homologies.
We shall describe a whole sequence of such invariants 
for four-dimensional starshaped domains in Section~\ref{s:ech}.
We here briefly outline the construction of another important set of spectral invariants.

\begin{figure}[h]   
 \begin{center}
  \psfrag{x}{$x$}
  \psfrag{u}{$u$}
  \psfrag{g}{$\gamma$}
  \psfrag{ug}{$u_{\infty} = \gamma (0)$}
  \leavevmode\includegraphics{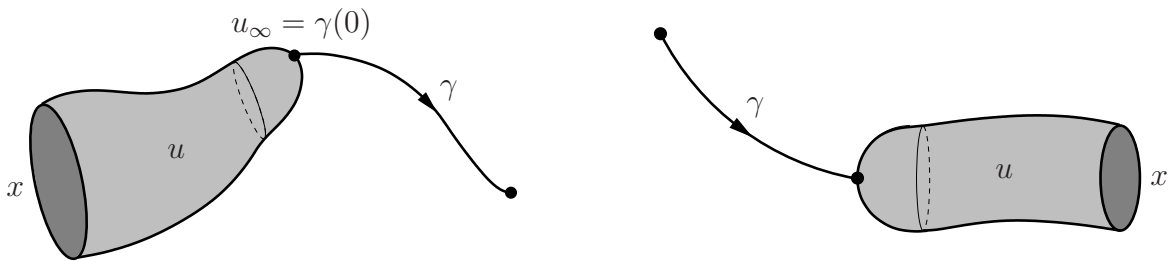}
 \end{center}
 \caption{Hybrid flow lines used to define $\Psi_{\PSS}$ and $\Phi_{\PSS}$}  \label{fig.pss}
\end{figure}

Take again a closed aspherical symplectic manifold $(M,\go)$.
As we have seen in Section~\ref{s:comp}
an isomorphism $\HF_k(M,\go) \to \H_{k+n}(M)$
is obtained by identifying the Floer chain complex of a $C^2$-small autonomous Hamiltonian
with its Morse chain complex.
A different chain homotopy equivalence between the Floer and the Morse chain complex
was constructed by Piunikhin, Salamon and Schwarz in~\cite{PiSaSch94}.
It is constructed by using hybrid flow lines, that are somewhat similar to the hybrid flow lines
used to construct the chain map~\eqref{cm:Phi}.
Choose a non-degenerate Hamiltonian function $H \colon \TT \times M \to \RR$ 
and a Morse function $f \colon M \to \RR$,
and an almost complex structure $J$ and a Riemannian metric $g$ such that the 
Floer complex of $(H,J)$ and the Morse complex of $(f,g)$ are defined.
Take a smooth cut-off function $\beta \colon \RR \to [0,1]$ with $\beta(s)=1$ for $s \leqslant 0$
and $\beta(s)=0$ for $s \geqslant 1$.
Let $\cm$ be the space of solutions of the $s$-dependent Floer equation~\eqref{floer.s} 
with finite energy $\int_{\RR \times \TT} |\frac{\partial u}{\partial s}|^2 \,ds dt$.
Then $u(s,\cdot)$ tends to a 1-periodic orbit~$x$ of~$H$ as $s \to -\infty$
and to a point $u_\infty \in M$ as $s \to +\infty$.
Now define a chain map $\CF_k(H) \to \CM_{k+n}(f)$ by counting pairs~$(u,\gg)$
where $u \in \cm$ and $\gg$ is a negative gradient flow line of~$f$ with $\gg(0) = u_\infty$,
see the left drawing in Figure~\ref{fig.pss}.

\ni
A chain map $\CM_{k+n}(f) \to \CF_{k}(H)$ is defined in a similar way, see the right drawing 
in Figure~\ref{fig.pss}. 
The two maps induced in homology 
$$
\Psi_{\PSS} \colon \HF_k(H) \to \HM_{k+n}(f),  \qquad
\Phi_{\PSS} \colon \HM_{k+n}(f) \to \HF_k(H)
$$
turn out to be inverses. In particular $\Psi_{\PSS}$ is an isomorphism.

The isomorphism $\Phi_{\PSS}$ can be used to define so-called action selectors, 
that associate to each homology class $\ga \in \H_* (M)$ the action of a 1-periodic orbit of~$H$.
The inclusion $\CF_*^a(H) \to \CF_*(H)$ descends to a map $\iota^a_* \colon \HF_*^a(H) \to \HF_*(H)$. 
For $\ga \in \H_* (M) \cong \HM_*(M)$ define the real number
\begin{equation} \label{def.cga}
c_\ga(H) \,=\, \inf \left\{ a \mid \Phi_{\PSS} (\ga) \in \im \iota^a_* \right\}.
\end{equation}
In words, $c_\ga(H)$ is the smallest $a$ such that the image of $\ga$ under the PSS isomorphism
can be represented by a sum of 1-periodic orbits of $H$ all of whose actions do not exceed~$a$.
For the generator $[\pt]$ of $\H_0(M)$, this definition is similar to the definition of~$\gg (M)$
in~\eqref{def:gg}.
The action selectors $c_\ga (H)$ have many natural properties, that make them quite computable. 
For instance, the pair of paints product~\eqref{e:pp} implies that
$$
c_{\ga * \gb}(H \# K) \,\leqslant\, c_{\ga} (H) + c_{\gb} (K)
$$
where $\ga * \gb$ is the product defined by Poincar\'e duality and the cup product, 
and where 
$$
(H \# K) (x,t) \,=\, H(x,t)+K( (\phi_{H}^t)^{-1}x,t)$$ 
is the Hamiltonian function 
that generates the path $\phi^t_H \circ \phi^t_K$.
Taking as $\ga$ the generator $[\pt]$ of~$\H_0(M)$ one defines for every open subset~$U$ of~$M$
the capacity
$$
c(U) \,=\, \sup \left\{ 
c_{\scriptscriptstyle [\pt]} (H) \mid  
H \colon \TT \times U \to \RR \mbox{ has compact support} 
\right\} .
$$ 

Action selectors as above have been first constructed for compactly supported Hamiltonian
functions on~$\RR^{2n}$ by Viterbo~\cite{Vi92} with the help of generating functions
and by Hofer--Zehnder~\cite{HoZe94} by a topological linking argument in the loop space $H^{\frac 12}$
that we encountered in Section~\ref{s:hamiltonian}.
The above Floer homological selector was constructed by Schwarz~\cite{Sch00}
and extended to arbitrary closed symplectic manifolds by Oh~\cite{Oh06} and Usher~\cite{Us08}.
A similar construction for a class of symplectically aspherical manifolds comprising Liouville domains
and their completions was given in~\cite{FrSch07}.
Action selectors and their capacities have many applications to Hamiltonian dynamics and symplectic topology,
see for instance the references above and~\cite{EP09, FrGiSch05, HRS16}.
Polterovich~\cite{po02} used the PSS action selector to show that certain groups have no interesting representations
in the symplectomorphism groups of closed symplectically aspherical manifolds.
For instance, every homomorphism of a finitely generated subgroup of $\SL (n;\ZZ)$ with $n \geqslant 3$
to the group of area preserving diffeomorphisms of a closed surface of genus~$\geqslant 2$
has finite image.
Polterovich et al.\ used the PSS action selector to detect new manifestations of symplectic rigidity,
which take place in certain function spaces associated to symplectic manifolds, see~\cite{PoRo14}.

The action functional~\eqref{e:unperturbed} was first considered by Rabinowitz~\cite{Rab78} in his proof of the
Weinstein conjecture for starshaped hypersurfaces in~$\RR^{2n}$.
The Floer homology~$\RFH$ of this functional was constructed by Cieliebak and Frauenfelder in~\cite{CieFra09}.
In Theorem~\ref{t:RF} one leafwise intersection point and assertion~(ii) was established
in~\cite{AlbFra10}, while the cuplength estimate was proven in~\cite{almo10}.
RFH was used in~\cite{CieFraPat10} to study the dynamics of a charged particle 
in a magnetic field at different energy levels.
It turns out that RFH is closely related to symplectic homology~\cite{CieFraOan10}, 
and that it can be interpreted as a Floer homology without Lagrangian multiplier~\cite{AbMe15, CiOa15}. 
For a survey of RFH and its applications see~\cite{AlbFra12}.

\section{Floer homology for Lagrangian intersections} \label{s:Lag}

In this section we first outline the construction of Lagrangian Floer homology in the simplest situation. 
Among its applications are the proof of the Arnol'd conjecture for Lagrangian intersections in~$T^*Q$,
and a variant of the construction gives another proof of the Hamiltonian Arnol'd conjecture~\ref{arnold}.
We then discuss the Lagrangian Floer homology of two fibers in a cotangent bundle, and show how it implies
lower complexity bounds (positive topological entropy) for Reeb flows on many spherisations.
Throughout this section we use $\ZZ_2$ coefficients for simplicity.

\subsection{Outline of the construction}
\begin{figure}[h]   
 \begin{center} 
  \psfrag{x}{$x$}
  \psfrag{y}{$y$}
  \psfrag{L}{$L$}
  \psfrag{L'}{$L'$}
  \psfrag{f}{$\phi_H(L)$}
  \psfrag{u}{$u$}
  \psfrag{0}{$0$}
  \psfrag{1}{$1$}
  \psfrag{s}{$s$}
  \psfrag{t}{$t$}
  \psfrag{na}{\footnotesize $-\nabla \ca_H$}
  \psfrag{F'}{\footnotesize $\overline \pp u =0$}
  \psfrag{F}{\footnotesize $\overline \pp u = J X_H$}
  \leavevmode\includegraphics{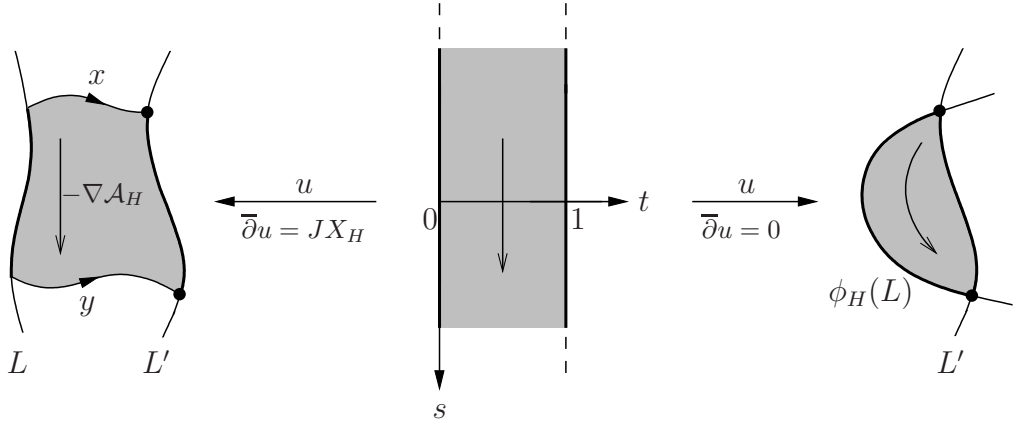}
 \end{center}
 \caption{A solution $u$ in $\widehat \cm_{H,J}(x,y)$ (left), corresponding to  an ``eye'' after a Hamiltonian coordinate change (right)}
 \label{fig.Lag}
\end{figure}

Let $L, L'$ be two Lagrangian submanifolds in a symplectic manifold~$(M,\go)$.
We consider the simplest case, in which $(M,\lambda)$ is the completion of a Liouville domain
and $L,L'$ are compact and exact.
The reader may wish to have a concrete example in mind and take
$L=L'=Q$ to be the zero section in $(T^*Q, \omega_{\can} = d \lambda_{\can})$.
The dictionary to Morse homology this time is as follows:
The space of candidates is the space of paths from $L$ to~$L'$,
$$
\cp (L,L') \,=\, \left\{ x \in C^\infty ([0,1], M) \mid x(0) \in L,\, x(1) \in L' \right\} .
$$
The ``Morse function'' is again the action functional from classical mechanics, 
$$
\ca_H (x) \,=\, \int_0^1 \bigl( \lambda(\dot x(t)) - H_t (x(t)) \bigr) \,dt,
$$
where $H \colon [0,1] \times M \to \RR$.
The critical points of~$\ca_H$ are the Hamiltonian chords of time~1 from $L$ to~$L'$, that is, 
the flow lines of $\phi_H^t$ starting on~$L$ at $t=0$ and arriving on~$L'$ at time~$t=1$.
Note that these chords can be identified with the points in $\phi_H^1(L) \cap L'$.
The Morse condition for~$\ca_H$ turns out to be that all intersection points are transverse:
$\phi_H^1(L) \pitchfork L'$.
As in \S \ref{s:idea} we define an $L^2$-inner product on the tangent spaces of~$\cp (L,L')$ 
with the help of an 
$\go$-compactible almost complex structures~$J$, and interprete negative gradient flow lines 
of~$\ca_H$ as solutions of the Floer equation:
For any pair $x,y$ of critical points we look at the space $\widehat \cm_{H,J}(x,y)$ 
of solutions $u \colon \RR \times [0,1] \to M$ of 
\begin{equation}
\label{floerLag}
\frac{\partial u}{\partial s} + J(u) \left( \frac{\partial u}{\partial t} - X_{H_t}(u) \right) = 0,
\end{equation}
with Lagrangian boundary conditions $u(s,0) \in L$ and $u(s,1) \in L'$ for all $s \in \RR$, and
asymptotic conditions
\[
\lim_{s\rightarrow -\infty} u(s,\cdot) = x \qquad \mbox{and} \qquad 
\lim_{s\rightarrow +\infty} u(s,\cdot) = y \qquad \mbox{in } \cp (L,L'),
\]
see the left part of Figure~\ref{fig.Lag}, in which we abbreviated 
$\overline \pp u = \frac{\pp u}{\pp s} + J(u) \frac{\pp u}{\pp t}$.
The Lagrangian boundary conditions crucially enter the proof that at every solution~$u$ 
the linearized operator~\eqref{operator} is Fredholm and in the compactness analysis of the 
spaces~$\widehat{\cm}_{H,J}(x,y)$.

%

In favorable situations (such as $c_1(M,\omega) =0$)
one can associate to the finitely many critical points of $\ca_H$ an index~$\mu$, 
with the property that the dimension of $\widehat \cm_{H,J} (x,y)$ is equal to $\mu(x)-\mu(y)-1$,
cf.\ \S \ref{s:twofibers}.
For defining an (ungraded) Floer homology, however, having the index is not necessary:
Just define $\CF (L,L';H)$ to be the $\ZZ_2$ vector space freely generated by the finitely many 
critical points of~$\ca_H$, and set
$$
\partial \1  x \,=\, \sum_{\substack{y \in \Crit \ca_H}} \nu (x,y)\, y
$$
where $\nu (x,y)$ is the finite number of 1-dimensional components of $\widehat \cm_{H,J}(x,y)$.
Then one proves $\pp^2 =0$ as before and can thus define the Lagrangian Floer homology $\HF (L,L';H,J)$
as the homology of the chain complex $\{ \CF (L,L';H), \pp \}$.
As in \S \ref{s:invariance} one shows that this group does not depend on the generic choice of~$J$
nor on~$H$ as long as $\phi_H^1(L) \pitchfork L'$; it is thus denoted by $\HF (L,L')$.

The invariance has two consequences:
First, $\HF (L,L')$ vanishes if $L$ and~$L'$ are disjoint, since then for a $C^1$-small Hamiltonian
there are no Hamiltonian chords from $L$ to~$L'$.
Second, taking the Hamiltonian $C^2$-small and also independent of time, 
one finds that in the case $L=L'$ the Floer homology $\HF (L) :=\HF (L,L)$ of~$L$
is isomorpic to its (ungraded) homology,
\begin{equation} \label{iso:HL}
\HF (L) \,\cong\, \H (L) .
\end{equation}

In the Morse theoretic picture above the generators of $\HF (L,L';H,J)$ are the Hamiltonian 
time-1-chords from $L$ to~$L'$, and the Floer gradient lines~$u$ are the solutions of the perturbed 
Cauchy--Riemann equation~\eqref{floerLag}.
We have already noticed that the set of generators can be identified with the intersection points 
$\phi_H^1(L) \cap L'$. Under this identification the Floer equation~\eqref{floerLag} becomes 
the unperturbed Cauchy--Riemann equation 
\begin{equation*} 
\frac{\partial u}{\partial s} + J(u) \frac{\partial u}{\partial t} = 0,
\end{equation*}
with Lagrangian boundary conditions $u(s,0) \in \phi_H^1(L)$ and $u(s,1) \in L'$ for all $s \in \RR$, and
asymptotic conditions 
\[
\lim_{s\rightarrow -\infty} u(s,\cdot) = x \qquad \mbox{and} \qquad 
\lim_{s\rightarrow +\infty} u(s,\cdot) = y
\]
where now $x,y$ are points in $\phi_H^1(L) \cap L'$.
One can thus view the Lagrangian Floer homology of $L$ and~$L'$ also as the homology 
of the chain complex freely generated by the Lagrangian intersections $\phi_H^1(L) \pitchfork L'$ 
and with boundary operator defined by counting ``$J$-holomorphic eyes'' 
as on the right of Figure~\ref{fig.Lag}.

\subsection{Applications to Lagrangian submanifolds}

Applying the above construction to $L=L'$ the zero section of $T^*Q$, we obtain at once 

\begin{theorem} 
[\bf Arnol'd conjecture for Lagrangian intersections in $T^*Q$] \label{t:lag}
Let $Q$ be a closed $n$-dimensional manifold, and let $\phi_H$ be any Hamiltonian diffeomorphism of~$T^*Q$ such that $\phi_H(Q)$ and $Q$ intersect transversally. 
Then $\# (\phi_H(Q) \cap Q) \,\geqslant\, \sum_{i=1}^n b_i(Q;\ZZ_2)$.
\end{theorem}

In a similar way Floer~\cite{Flo88:Lag} proved the following variant of Theorem~\ref{t:lag}.

\begin{theorem} \label{t:lagclosed}
Let $L$ be a closed Lagrangian submanifold in a compact symplectic manifold~$(M,\go)$
of dimension~$2n$ 
such that the relative homotopy group $\pi_2 (M,L)$ vanishes,
and let $\phi_H$ be any Hamiltonian diffeomorphism of~$M$ such that $\phi_H(L)$ and~$L$ 
intersect transversally. 
Then~$\# (\phi_H(L) \cap L) \geqslant \sum_{j=0}^n b_j(L;\ZZ_2)$.
\end{theorem}

Note that the topological assumption $\pi_2 (M,L) = 0$ cannot be omitted, 
as the example of a small contractible embedded circle on the torus $\RR^2 / \ZZ^2$ shows.
Theorem~\ref{t:lagclosed} on Lagrangian intersections can be used to prove 
the Hamiltonian Arnol'd conjecture~\ref{arnold}
(at least if $\pi_2(M)=0$ and over $\ZZ_2$ coefficients), 
giving one confirmation of Weinstein's dictum~\eqref{dictum:Wein}.
Indeed, 
the diagonal $\Gamma_0 = \{ (x,x) \mid x \in M \}$ is a Lagrangian submanifold of the product
$(M \times M, \omega \oplus (-\go))$, and the fixed points of~$\phi_H$ on~$M$
correspond to the intersections of $\Gamma_0$ with the graph $\Gamma_1 = \{ (x,\phi_H(x)) \mid x \in M \}$,
which are transverse if and only if all the fixed points of~$\phi_H$ are non-degenerate.
The isotopy of Lagrangian submanifold $\Gamma_t = \{ (x,\phi_{H}^tx) \mid x \in M \}$
from $\Gamma_0$ to $\Gamma_1$ is generated by a Hamiltonian flow on~$M \times M$.
Furthermore, the assumption $\pi_2 (M) =0$ implies that $\pi_2(M \times M, \Gamma_0) =0$,
since $\pi_2(M \times M) = \pi_2(M) \times \pi_2(M) =0$ and since the last arrow of the part
$$
0 = \pi_2(M \times M) \to \pi_2(M \times M,\Gamma_0) \to \pi_1(M) \to \pi_1(M \times M) 
$$
of the exact sequence of relative homotopy groups is injective.

\m
As we have seen in \S \ref{s:basicsII},
Arnol'd's nearby Lagrangian conjecture
holds for $Q = S^1$, and it is known also for $Q = S^2$ by a result of R.\ Hind~\cite{Hi12}.
For all other manifolds~$Q$, the strongest result known so far is 
the following recent theorem due to Abouzaid and Kragh~\cite{AbKr16}.

\begin{theorem} \label{t:nearby}
Let $Q$ be a closed manifold and $L \subset T^*Q$ a closed exact Lagrangian submanifold. 
Then the restriction of the canonical projection $\pi \colon T^*Q \rightarrow Q$ to~$L$ is 
a simple homotopy equivalence.
\end{theorem}

Here ``simple'' means that the Whitehead torsion of the map $\pi|_L$ vanishes.
(The hierarchy is: homotopy equivalent $\Leftarrow$ simple homotopy equivalent $\Leftarrow$ homeomorphic 
$\Leftarrow$ diffeomorphic.)
The main tools in the proof are a Floer theoretic definition of the Whitehead torsion,
and ingredients from the Fukaya category of~$T^*Q$,
which is a category whose objects are the exact Lagrangian submanifolds in~$T^*Q$ 
and whose morphisms are the Floer chain complexes of pairs of such Lagrangians. 

\s
Theorem~\ref{t:nearby} has a beautiful application to a question on the relation between the
differential topology of a manifold~$Q$ and the symplectic topology of its cotangent bundle:
Diffeomorphic manifolds have symplectomorphic cotangent bundles,
since a diffeomorphism $\gf \colon Q_1 \to Q_2$ induces
the symplectomorphism $\gf^* \colon (T^*Q_2, \omega_{\can}) \to (T^*Q_1, \omega_{\can})$,
where $\gf^* (q,p) = \left( \gf^{-1} (q), d\gf(q)^T (p) \right)$.
It is an open problem going back to Arnol'd and Eliashberg whether for homeomorphic manifolds 
the converse is also true.
Theorem~\ref{t:nearby} shows that this is so for Lens spaces.
Recall that for an integer $p \geqslant 2$ and integers $q_1, \dots, q_n$ relatively prime to~$p$,
the lens space $L(p;q_1,\dots,q_n)$ is the quotient $S^{2n-1}/G_p $ of the sphere
$S^{2n-1} = \left\{ (z_1, \dots, z_n) \in \CC^n \mid \sum_{i=1}^n |z_i|^2 =1 \right\}$
by the action
$$
\zeta (z_1, \dots, z_n) \,=\,  (\zeta^{q_1} z_1, \dots, \zeta^{q_n} z_n)
$$
of the cyclic group $G_p = \{ \zeta \in \CC \mid \zeta^p = 1\}$.
While the classification of Lens spaces up to simple homotopy type is finer than 
the one up to homotopy equivalence
(for instance, the 3-dimensional lens spaces $L(7;1,1)$ and~$L(7;1,2)$ are homotopy equivalent, but not simple homotopy equivalent),
it is classically known following Franz and Reidemeister that the simple homotopy type 
of lens spaces already determines their diffeomorphism type. 
(For modern accounts see \cite[\S 9]{Mi66} or \cite[\S 11]{Tu01}.)
Theorem~\ref{t:nearby} thus implies:

\begin{corollary} \label{c:lens}
Two lens spaces are diffeomorphic if and only if their cotangent bundles are symplectomorphic.
\end{corollary}

\subsection{The Lagrangian Floer homology of two fibers in $T^*Q$} \label{s:twofibers}

Take two fibers $L_q = T_q^*Q$ and $L_{q'} = T_{q'}^*Q$ in the cotangent bundle over a closed 
manifold~$Q$, and fix a Riemannian metric~$g$ on~$Q$.
The Lagrangian Floer homology $\HF (L_q, L_{q'}; G)$
for the geodesic Hamiltonian $G(q,p) = \frac 12 \|p\|^2$
is the Lagrangian analogue of the symplectic homology of~$T^*Q$:
The geodesic time-1 chords from $L_q$ to~$L_{q'}$ correspond to the intersection points~$\phi_G^1 (L_q) \cap L_{q'}$,
and the action $\ca_G(x)$ of such a chord $x(t) = (q(t),p(t))$ 
is the constant value $G(x(t)) = \frac 12 \|p(t)\|^2$.
For a generic $q'$, these intersections are transverse, 
and hence there are finitely many chords below a given action.
Since the action~$\ca_G$ does not increase along solutions~$u$ of its negative gradient 
equation~\eqref{floerLag}, the boundary operator~$\pp$ on the chain group $\CF (L_q,L_{q'}; G)$ 
is then well-defined.
The resulting Floer homology $\HF (L_q, L_{q'}; G)$ does not depend on the specific choice of~$G$.

\begin{figure}[h]   
 \begin{center} 
  \psfrag{q}{\footnotesize $q$}
  \psfrag{q'}{\footnotesize $q'$}
  \psfrag{1}{\footnotesize $1$}
  \psfrag{Lq}{$L_q$}
  \psfrag{Lq'}{$L_{q'}$}
  \psfrag{Q}{$Q$}
  \leavevmode\includegraphics{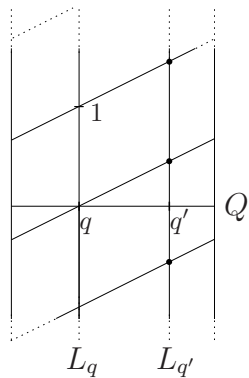}
 \end{center}
 \caption{$\phi_G (L_q) \cap L_{q'}$ in $T^*Q$ for $Q = S^1 = [0,1]/\sim$}
 \label{fig:cyl}
\end{figure}
%

For the application in the next section we need the refinement of this homology by its action filtration.
By taking only geodesic time-1 chords of action $\leqslant a$ we obtain 
the filtered Floer homology groups $\HF^a (L_q,L_q';G)$. 
In general, these groups certainly do depend on~$G$, 
but if $H$ is another Hamiltonian with $H \leqslant G$ that agrees with~$G$ outside a compact set,
then there is a continuation map 
\begin{equation} \label{e:Phi}
\Phi_{HG}^a \colon \HF^a (L_q,L_q';H) \to \HF^a (L_q,L_q';G) ,
\end{equation}
and if one can find a path $H_s$ of Hamiltonians from $H$ to~$G$ such that $a$ 
is a regular value of~$\ca_{H_s}$ for all $s$, then $\Phi_{HG}^a$ is an isomorphism.
The groups $\HF^a (L_q,L_q';H)$ can be equipped with a grading that for a geodesic Hamiltonian~$G$ 
is given by the Morse index of a geodesic arc.

Denote by $\Omega_{q,q'}Q$ the space of smooth paths $x \colon [0,1] \to Q$ from $q$ to~$q'$, 
and for $a \in \RR$ write $\Omega^a_{q,q'}Q$
for its subset of paths whose energy 
$\ce (x) = \frac 12 \int_0^1 \|x(t) \|^2 \,dt$ does not exceed~$a$.
The construction in~\S \ref{s:disc} can be adapted to the situation at hand to show that
\begin{equation} \label{iso:based}
\HF^a_* (L_q,L_q';G) \,\cong\, \H_* (\Omega^a_{q,q'}Q) .
\end{equation}

\subsection{Lower bounds for the topological entropy of Reeb flows}
\label{s:entropy}

The topological entropy of a self-map of a metrisable compact space
is a numerical invariant that gives a measure for 
the complexity of the orbit structure of the map.
For instance, isometries, periodic maps and gradient flows have vanishing topological entropy, 
and expanding maps and geodesic flows of negative curvature metrics have positive topological entropy.
A classical problem is to show that a {\it generic}\/ map in a given class of mappings (e.g.\ geodesic flows) 
has positive entropy. A different problem is to see that {\it all}\/ maps 
(of a certain type, or in a given isotopy class) on {\it certain}\/ manifolds have positive topological entropy.
Here we consider this second problem for Reeb flows on closed contact manifolds.
Since these flows model systems without friction, one can expect that they
often have positive topological entropy. 

There are many equivalent definitions of topological entropy~\cite[\S 3.1]{HaKa95}.
For a $C^\infty$ smooth diffeomorphism~$\gf$ of a compact $d$-dimensional manifold~$M$,
the topological entropy $h_{\top}(\gf)$ can be seen as the maximal exponential volume growth of
submanifolds under iterates of~$\gf$:  
We say that a sequence of positive numbers $a_n$ has exponential growth if its exponential growth rate 
$$
\gamma (a_n) \,=\, \liminf_{n \to \infty} \frac 1n \log a_n 
$$
is positive. 
Let $\cb$ be the set of all smoothly embedded closed balls~$B \subset M$, of any dimension $k \in \{0, \dots, d\}$,
and define the exponential volume growth of~$\gf$ by
$$
\gg_{\vol} (\gf) \,=\,  \sup_{B \in \cb} \gg \bigl( \Vol (\gf^n (B)) \bigr) .
$$
Here, the volume is taken with respect to the measure induced on $\gf^n(B)$ by a Riemannian 
metric~$g$ on~$M$, and clearly $\gg_{\vol} (\gf)$ does not depend on~$g$.
Yomdin~\cite{Yom87} ($\geqslant$) and Newhouse~\cite{New88} ($\leqslant$) proved that 
\begin{equation} \label{e:Yomdin}
h_{\top}(\gf) \,=\, \gg_{\vol} (\gf) .
\end{equation}

Recall that a finitely generated group has exponential growth if for one (and hence any)
set~$S$ of generators, the number of different words of length~$\leqslant n$
that can be written by letters in $S \cup S^{-1}$ has exponential growth.
Examples are the fundamental groups of closed manifolds that admit a Riemannian metric of negative 
curvature.
Further, we denote the based loop space of a closed manifold by~$\Omega Q$, 
and if $Q$ has finite fundamental group 
we say that the homology of~$\Omega Q$ has exponential growth 
if the sequence
\begin{equation} \label{e:expOmega}
\sum_{k=0}^n \dim \H_k (\Omega Q; \ZZ_2)
\end{equation}
has exponential growth. An example is the complex projective plane blown-up in at least two points. 
A ``generic'' closed manifold satisfies one of these two conditions.
The topological entropy of a flow is defined as the topological entropy of its time-1-map.

\begin{theorem} \label{t:htop}
Let $Q$ be a closed manifold which satisfies one of the following two conditions:
\begin{enumerate}[\rm (i)]
\item the fundamental group of $Q$ has exponential growth;
\item the fundamental group of $Q$ is finite and the homology of~$\Omega Q$ 
has exponential growth. 
\end{enumerate}
Then every Reeb flow on the spherisation~$S^*Q$ has positive topological entropy. 
\end{theorem}

\m 
The idea of the proof is to show that any (and hence one) Legendrian sphere $S^*_qQ$
grows exponentially under any Reeb flow.
As in Example~\ref{ex:spherisation} we can identify a Reeb flow~$\phi_\ga^t$ on~$S^*Q$
with the restriction of a Hamiltonian flow~$\phi_H^t$ on~$T^*Q$ 
to a fiberwise starshaped hypersurface~$\Sigma \subset T^*Q$,
where $H \colon T^*Q \to [0,\infty)$ is a Hamiltonian function with 
$H^{-1}(1)=\Sigma$ that is fiberwise homogeneous of degree two.
We thus wish to show that any fiber-sphere~$\Sigma_q = \Sigma \cap T_q^*Q$ grows exponentially 
under the Hamiltonian flow~$\phi_H^t$.
Since $H$ is homogeneous of degree two, it is not hard to see that
this growth agrees with the growth of the disc $D_q \subset T_q^*Q$ 
bounded by~$\Sigma_q$,
$$
\gg (\Vol(\phi_H^n(\Sigma_q))) \,=\, \gg (\Vol(\phi_H^n(D_q))) .
$$
In view of the identity~\eqref{e:Yomdin} it thus suffices to show that
\begin{equation*} 
\gg \bigr( \Vol (\phi_H^n(D_q)) \bigr) >0.
\end{equation*}

Fix $n$, and let $U \subset Q$ be the open and dense set of points~$q'$ for which 
$\phi_H^n (D_q)$ intersects $D_{q'}$ transversally.
Our goal is to show that there are positive constants $a,b$ that do not depend on~$n$
such that  
\begin{equation*}
\# \bigl( \phi_H^n (D_q) \cap D_{q'} \bigr) \,\geqslant\, e^{an-b} \qquad 
\mbox{for all\, $q' \in U$.}
\end{equation*}
Since $U$ has positive measure, this would clearly imply that 
$\gg \bigl( \Vol (\phi_H^n(D_q)) \bigr) \geqslant a >0$,
cf.\ Figure~\ref{fig:cyl}.

Since $\phi_{nH}^t = \phi_H^{nt}$, the diffeomorphism $\phi_H^n$ is generated in time~1
by~$nH$,
and since $H$ is homogenous of degree two,
a time 1 flow line~$x$ of $\phi_{nH}^t$ lies in $\{ H \leqslant 1\}$ 
if and only if $\ca_{nH}(x) \leqslant n$.
For $q' \in U$, the points $\phi_H^n (D_q) \cap D_{q'}$ that we wish to count are thus exactly
the generators of the Floer chain groups $\CF^n (L_q , L_{q'}; nH)$.

If the fundamental group of $Q$ has exponential growth, then 
clearly the sequence
$$
\dim \H_0 (\Omega^{n^2} Q)
$$ 
has exponential growth.
And if the fundamental group is finite and the sequence~\eqref{e:expOmega} built by truncating the index 
has exponential growth, 
then also the sequence 
$$
\sum_{k \geqslant 0} \dim \H_k (\Omega^{n^2} Q)
$$ 
built by truncating the energy has exponential growth
by a theorem of Gromov~\cite{Gro78}. 
We therefore find positive constants $a,b$ such that
$$
\dim \H (\Omega_{q,q'}^{n^2}Q) \,\geqslant\, e^{an-b} \qquad \mbox{for all $n \in \NN$ and all $q,q' \in Q$.}
$$ 

Now fix a Riemannian metric and assume for a moment that $H=G$ is its geodesic Hamiltonian.
Then using the isomorphism~\eqref{iso:based} we can estimate
\begin{eqnarray*}
\# \bigl( \phi_G^n (D_q) \cap D_{q'} \bigr) &=& \dim \CF^n (L_q ,L_{q'}; nG) \\
&\geqslant& \dim \HF^n (L_q, L_{q'}; nG) \\
&=& \dim \H (\Omega_{q,q'}^{n^2}Q) \\
&\geqslant& e^{an-b} \quad 
\end{eqnarray*}
for all $q' \in U$, as we wished to show.

If $H$ is not a geodesic Hamiltonian we use a sandwiching argument: 
Choose positive constants $c_- < c_+$ such that
$$
c_-G \,\leqslant\, H \,\leqslant\, c_+ G 
$$
which is possible since the levels of $H$ are fiberwise starshaped.
Then set $G_+ := c_+G$ and choose interpolating functions $G_-$ and $\widetilde H$ 
such that
\begin{equation} \label{e:GH}
G_- \,\leqslant\, \widetilde H \,\leqslant\, G_+ ,
\end{equation}
$G_- = c_-G$ and $\widetilde H = H$ on $\{ H \leqslant 1\}$,
and $G_- = \widetilde H = G_+$ far away from $\{ H \leqslant 1\}$.
Then there is the commutative diagram of Floer homologies 
$$  
\xymatrix{ 
&\HF^{n/\gs}(nG_+) \ar[dr]^{\iota^n}& \\
\HF^{n}(nG_-) \ar[ur]^{\cong}
\ar[rr]
\ar[dr] && \HF^{n}(nG_+) \\
& \HF^{n}(n\widetilde H) \ar[ur]&}
$$
in which we omitted $L_q,L_{q'}$ from the notation and set $\gs = c_+/c_-$.
The homomorphism $\iota^n$ is defined by inclusion of complexes, 
and all other maps are continuation maps as in~\eqref{e:Phi}, which exist thanks to~\eqref{e:GH}.
Since the diagram commutes, 
$$
\dim \HF^{n} (n\widetilde H) \,\geqslant\, 
\rank \iota^n \colon \HF^{n/\gs}(nG_+) \to \HF^{n}(nG_+) .
$$
We already know that $\HF^{n/\gs}(nG_+)$ grows exponentially, 
but in fact $\rank \iota^n$ grows exponentially, since under our assumptions
not only $\dim \H (\Omega_{q,q'}^{n^2}Q)$ but also
$$
\rank \iota^n \colon \H (\Omega_{q,q'}^{n^2}Q) \to \H (\Omega_{q,q'}Q)
$$ 
grows exponentially
and since the isomorphisms~\eqref{iso:based} commute with the maps induced by the inclusions 
$\CF^{a}(L_q,L_{q'};G) \subset \CF^{a'}(L_q,L_{q'};G)$ and $\Omega_{q,q'}^aQ \subset \Omega_{q,q'}^{a'}Q$.
\proofend

Usually, a Floer homology is computed by using its invariance and showing
that for a special Hamiltonian ($C^2$-small and time independent, or geodesic)
the Floer homology can be identified with the homology of a space. 
Examples for this classical procedure are the isomorphisms~\eqref{iso:HF}, \eqref{iso:HL}
and~\eqref{iso}, \eqref{iso:based}.

If one cannot find an isomorphism to the homology of a space, 
one can still hope to understand the Floer homology of one Hamiltonian~$H$. 
With some luck one finds a~$H$ for which one can determine all its Hamiltonian orbits, 
and with even more luck their indices have all the same parity, so that the boundary operator
vanishes and the Floer homology is freely generated by the Hamiltonian orbits of~$H$.
This strategy was used by Alves~\cite{Alv16a, Alv16b} to extend Theorem~\ref{t:htop} in dimension~3 
far beyond spherisations. 
To save words, we say that a closed contact manifold~$(Y,\xi)$ has {\it positive entropy}\/ 
if for every contact form~$\ga$ for~$\xi$ the time-1-map~$\phi_\ga$ of the Reeb flow
of~$\ga$ has positive topological entropy.

\begin{theorem}
There is an infinite family of closed $3$-manifolds 
that carry infinitely many different contact structures of positive entropy.
\end{theorem}

The complexity of the Reeb flows in these examples still comes (as in Theorem~\ref{t:htop})
from the rich topology of~$Y$. 
But very recently, Alves and Meiwes~\cite{AlMe17}
constructed many examples of contact manifolds of dimension~$\geqslant 5$ 
that have very simple topology but still positive entropy.
The most spectacular examples are odd-dimensional spheres with their standard smooth structure.

\begin{theorem} \label{t:AlMe}
The spheres $S^{2n-1}$ with $n \geqslant 4$ and $S^2 \times S^3$ 
admit contact structures of positive entropy.  
\end{theorem}

The contact structures in this theorem are not even exotic: 
They are contact structures on the boundary of Stein domains.
But they are very exotic from a contact dynamical point of view.
Indeed, the usual Reeb flow of their usual contact structure 
(the Hopf flow on $S^{2n-1}$ and the geodesic flow of the round metric on~$S^3$)
are periodic and thus have vanishing topological entropy.

The proof of Theorem~\ref{t:AlMe} uses a more radical form of invariance.
For spheres, 
it starts from the co-disc bundle $D^*Q = \left\{ \frac 12 \|p\|^2 \leqslant 1 \right\}$ 
over a closed Riemannian manifold that is a homology sphere and has fundamental group of exponential growth. 
By~\eqref{iso:based} two Lagrangian fibers~$D_q$ and~$D_{q'}$ 
have rich Floer homology.
Applying to $D^*Q$ suitable surgery operations (plumbing and subcritical handle attachment) 
in the complement of~$D_q$
one arrives at a Liouville domain~$(W,\go)$ that still contains~$D_q$,
and whose contact boundary~$(Y,\xi)$ is a simply connected homology sphere 
and hence is homeomorphic to~$S^{2n-1}$ by the h-cobordism theorem.
The surgeries can be arranged such that $Y$ is even diffeomorphic to~$S^{2n-1}$.
Further, the Lagrangian Floer homology of~$D_q$ and a nearby fiber~$D_{q'}$ in~$W$ 
(that is still defined) is at least as rich as for the original~$D^*Q$.
Indeed,
a plumbing or handle attachment applied to an intermediate Liouville domain~$W_1$
away from $D_q$ and~$D_{q'}$ produces a Liouville embedding $W_1 \subset W_2$,
and there is again a Viterbo transfer map from the Floer homology 
of $(W_1,D_q,D_{q'})$ to the one of $(W_2,D_q,D_{q'})$,
that turns out to be an isomorphism in the case of a subcritical handle attachment 
(by a version of Cieliebak's theorem for symplectic homology from~\cite{cie02})
and injective in the case of plumbing
(by a theorem of Abouzaid and Smith~\cite{AbSm12}).

\begin{figure}[h]  
 \begin{center}
  \psfrag{Sq}{$S_q$}
  \psfrag{f}{$\phi_\ga^n(S_q)$}
  \leavevmode\includegraphics{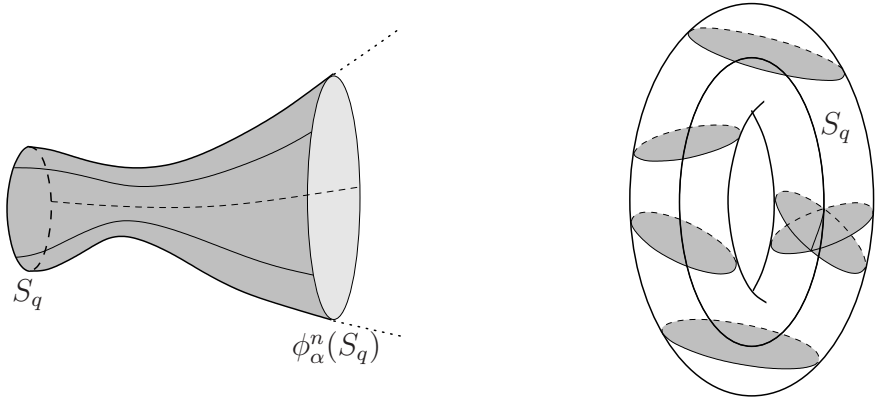}
 \end{center}
 \caption{The set $Z_\ga^n$ intersecting the neighbourhood $B^n \times S_q$ of $S_q$} 
 \label{fig.marcelo}
\end{figure}

To finish the proof one localizes the argument in the proof of Theorem~\ref{t:htop}:
While the Liouville domain $(W,\omega)$ is not a global fibration of Lagrangian discs anymore, 
it still contains a set of the form $B^n \times D_q$ in which each fiber over the
ball~$B^n$ is a disc~$D_{q'}$.
Let $S_{q'}$ be the boundary sphere of~$D_{q'}$.
Taking as set $U$ those $q'$ in~$B^n$ for which the ``cylinder''
$$
Z_\ga^n = \bigcup_{0 \leqslant t \leqslant n} \phi_\ga^n (S_q)
$$ 
intersects $S_{q'}$ transversally,
and by a sandwiching argument as at the end of the proof of Theorem~\ref{t:htop},
it follows that for $q' \in U$ the number of Reeb chords from~$S_q$ to~$S_{q'}$ 
still grows exponentially with time. 
Hence the volume of $Z_\ga^n$ 
and thus also the volume of $\phi_\ga^n (S_q)$
grows exponentially, see Figure~\ref{fig.marcelo}.

%

\subsection*{Historical notes and bibliography}
Theorem~\ref{t:lag} on Lagrangian intersections in~$T^*Q$
can be proved in an elementary way 
(Morse theory for generalized generating functions, see~\cite{LaSi85}), 
but Floer's proof by Lagrangian Floer homology leads to far reaching generalisations.
Floer created Lagrangian Floer homology (under the assumption that $\pi_2 (M,L)$ vanishes)
at the same time as Hamiltonian Floer homology, see~\cite{Flo88:unreg, Flo88:Lag, Flo89:Witten}.
His Theorem~\ref{t:lagclosed} can be proved whenever the Floer homology of~$L \subset (M,\go)$
can be defined. A thorough study of the obstructions to this is given by
Fukaya--Oh--Ohta--Ono in their books~\cite{FOOO}.

Seidel used Lagrangian Floer homology in~\cite{Se99} to provide examples of symplectic four-manifolds
that contain pairs of Lagrangian spheres that are smoothly isotopic,
but non-isotopic through Lagrangian embeddings.
In the same paper he constructed symplectomorphisms that are smoothly, but not symplectically, 
isotopic.
For generalisations of these results see~\cite{Ke14, Se00}.

The isomorphism~\eqref{iso:based} was proved in~\cite{as06}.
It was applied to proving Theorem~\ref{t:htop} by Macarini and the second author in~\cite{MaSch11},
based on earlier applications of Lagrangian Floer homology to volume growth in~\cite{FrSch05, FrSch06}.
Theorem~\ref{t:htop} generalizes work of Dinaburg, Gromov, Paternain and Petean on geodesic flows,
see~\cite{Pat.book}, and in turn can be generalized to those contactomorphisms~on spherisations
which can be reached by a contact isotopy that is everywhere positively transverse to the contact 
structure, see~\cite{Dah16}.
The proper setting of the classical results on the topological entropy of geodesic flows is thus 
contact topology.

The generalisation of the Floer homology $\HF (L_q,L_{q'};H)$ of two cotangent fibers 
in~$T^*Q$, that is used in the proof of Theorem~\ref{t:AlMe},
is the Lagrangian Floer homology $\HF (L,L';H)$
defined for two exact Lagrangian submanifolds~$L,L'$ of a Liouville domain $(W,d\lambda)$
that near the boundary~$\pp W$ are invariant under the (local) flow of the Liouville vector field,
and for a Hamiltonian~$H$ that on the symplectisation part is a function of~$\rho$ 
rapidly increasing with~$\rho$.
Taking direct limits over such functions as in the construction of the symplectic homology~$\SH (W)$ described in~\S \ref{s:SH}
one obtains the so-called {\it wrapped Floer homology}\/ $\HW (L,L')$.
The foundations of this Lagrangian version of the symplectic homology~$\SH (W)$ 
were laid by Abouzaid and Seidel in~\cite{AbSe10}.


\section{Floer homologies for contact manifolds}
\label{s:contact}

Until now the only (perturbed) $J$-holomorphic curves that we have used
are cylinders or strips, that connect periodic orbits or Hamiltonian chords. 
Given the power of the resulting Floer homologies, one is tempted to use
more general $J$-holomorphic curves to study problems in symplectic geometry and dynamics:
curves with possibly positive genus that may be asymptotic to more than one orbit. 
A general framework for such a theory in symplectic cobordisms between contact manifolds
is the {\it symplectic field theory}\/ (SFT) outlined by
Eliashberg, Hofer and Givental.
In this section we look at the special situation where the symplectic manifold
is the symplectisation of a contact manifold,
and describe two special versions of SFT:
Contact homology, and the embedded contact homology due to Hutchings, 
that is defined for 3-dimensional contact manifolds.

\subsection{Contact homology} \label{s:ch}
              
We consider the symplectisation of a contact manifold~$(Y,\ga)$,
but instead of looking at $Y \times (0,\infty)$ with symplectic form $d(r \ga)$
we make a convenient coordinate change and look at
$Y \times \RR$ with symplectic form $d (e^s \ga)$. 
Thus $(Y,\ga)$ now appears as $Y \times \{0\}$ in its symplectisation.
The almost complex structures are chosen somewhat more specific than in the
construction of symplectic homology: 
We require that $J$ leaves $\xi = \ker \ga$
invariant and is compatible with the symplectic form $d\ga$ on~$\xi$,
and maps $\frac{\pp}{\pp s}$ to the Reeb vector field~$R_\ga$
(these conditions are equivalent to~\eqref{e:rJ}), 
and in addition require that $J$ is invariant under shift along the $s$-coordinate.

Every closed orbit of the Reeb flow~$\gf_\ga^t$ is degenerate, since 
the differential of $\gf_\ga^t$ maps the Reeb vector field~~$R_\ga$ to itself.
A contact form $\ga$ is called non-degenerate if for every Reeb orbit~$x(t)$  
with $x(0) = x(T)$ for some $T>0$ all eigenvalues of 
$d\gf_{\ga}^T \colon \xi_{x(0)} \to \xi_{x(T)}$ are different from~$1$.
Any contact form $\ga$ for $\xi$ can be perturbed to a non-degenerate contact form.
Let $\cp_\ga$ be the set of unparametrized periodic Reeb orbits of~$\ga$:
Orbits with the same trace and the same multiplicity, 
but with different starting point, are taken only once, 
but orbits with the same trace but different multiplicities are different elements of~$\cp_\ga$.
We from now on assume $\ga$ is non-degenerate. Then $\cp_\ga$ is a countable set.

\begin{figure}[h]   
 \begin{center}
  \psfrag{x}{$x$}
  \psfrag{hy}{$y$}
  \psfrag{z}{$z$}
  \leavevmode\includegraphics{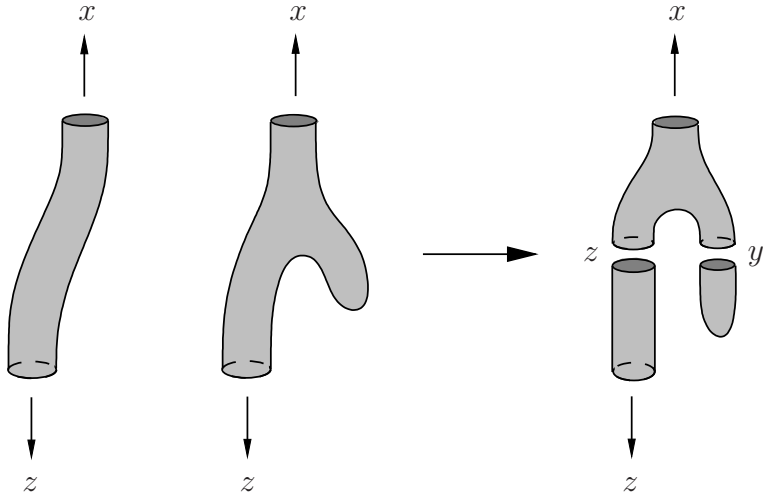}
 \end{center}
 \caption{A holomorphic plane bubbles off}  \label{fig.plane}
\end{figure}

The first attempt to build a Floer homology for $(Y,\ga)$ is by taking as generators the elements of~$\cp_\ga$
and by defining the boundary operator~$\pp x$ by counting the 1-dimensional components of the 
$J$-holomorphic cylinders in $Y \times \RR$ asymptotic to~$x$ and to another orbit~$y$. 
More precisely, we look $J$-holomorphic maps 
$u = (u_Y, u_{\RR}) \colon \RR \times \TT \to Y \times \RR$
such that 
\begin{equation} \label{e:floerch}
\lim_{s \to \mp \infty} u_{\RR}(s,t) = \pm \infty
\quad \mbox{ and } \quad
\lim_{s \to -\infty} u_Y(s,t) = x (T_xt), \quad 
\lim_{s \to +\infty} u_Y(s,t) = y (T_yt),
\end{equation}
where $T_x$ and $T_y$ are the periods of $x$ and~$y$.
Note that contrary to the Floer equation, there is no Hamiltonian term in the equation for the cylinders~$u$.
This attempt in general does not work, since now
the boundary of the space of such $J$-homomorphic cylinders 
can contain configurations different from broken cylinders:
Straight cylinders $z \times \RR$ over a Reeb orbit~$z$
are solutions of~\eqref{e:floerch},
and it can happen that a 1-parametric family of cylinders
from~$x$ to~$z$ breaks along~$z$ and another Reeb orbit~$y$,
at which a $J$-holomorphic plane bubbles off, see Figure~\ref{fig.plane}.
As a result, $\pp^2 =0$ does not hold in general.

This ``accident'' can be excluded by requiring, for instance, that $\cp_\ga$ contains no contractible
periodic orbits.
Decompose $\cp_\ga = \coprod \cp_\ga^h$ according to the free homotopy classes of maps $\TT \to Y$.
The cylinders solving~\eqref{e:floerch} respect this decomposition.
Fix a homotopy class~$h$.
To elements in~$\cp_\ga^h$ we associate a Conley--Zehnder index in $\ZZ_2$ as follows:
Fix a loop $\gamma$ in class~$h$ and a unitary trivialisation of~$\xi$ along~$\gamma$.
For $x \in \cp_\ga^h$ choose an annulus~$A_x$ with boundary $\gamma$ and~$x$,
and extend the trivialisation of~$\xi$ over~$A_x$.
Looking at the linearized flow $d \phi_{\ga}^t \colon \xi_{x(0)} \to \xi_{x (t)}$, 
$t \in [0,T_x]$ in this trivialisation of~$\xi$ over~$x$ we 
obtain $\mu_{\sCZ}(x;A_x) \in \ZZ$.
If we choose another annulus from $\gg$ to~$x$, we obtain an integer with the same parity,
thus $\mu_{\sCZ}(x) \in \ZZ_2$ is defined.

Following the steps in the construction of Hamiltonian Floer homology in Section~\ref{s:Ham}
we can now define the {\it cylindrical contact homology}\/ $\HC^{\cyl}_*(Y,\ga;h)$ with $* \in \ZZ_2$,
provided that $\ga$ is non-degenerate and has no contractible Reeb orbits and that $h$ is a primitive class.
Transversality can be achieved as for Hamiltonian Floer homology, since none of the $J$-holomorphic cylinder 
from~$x$ to~$y$ is multiply covered,
because $x$ and~$y$ are embedded orbits, $h$ being primitive.
This homology does not depend on the choice of a non-degenerate contact form without contractible orbits, 
and if $\gf \colon (Y,\xi_1) \to (Y,\xi_2)$ is a contactomorphism (namely a diffeomorphism 
$\gf \colon Y_1 \to Y_2$ with $d\gf (\xi_1) = \xi_2$), 
then the homologies $\HC^{\cyl}_*(Y,\ga;h)$ and  $\HC^{\cyl}_*(Y,\gf_*\ga;\gf_*h)$ are isomorphic.
This theory suffices to prove, for instance, the following results.

\begin{theorem} \label{t:CZ}
For $k \in \NN$ consider the contact structure $\xi_k$ on the 3-torus $\TT^3$
given by the kernel of the contact form
$$
\ga_k \,=\, \cos (2\pi k \theta_1) \,d \theta_2 + \sin (2\pi k \theta_1) \,d \theta_3,
\qquad (\theta_1, \theta_2, \theta_3) \in \TT^3 .
$$
\begin{itemize}
\item[\rm (i)]
For $k \neq \ell$
the contact structures $\xi_k$ and $\xi_\ell$ are not diffeomorphic. 

\s \ni
\item[\rm (ii)]
For every $k$
the fundamental group of the space of contact structures on~$\TT^3$ that are isotopic to~$\xi_k$
contains an infinite cyclic subgroup. 
\end{itemize}
\end{theorem}

For the proof of (i) one can take any class $h$ represented by a simple non-contractible curve in $\{ \theta_1 = 0\}$,
and one finds that for a suitable perturbation $\hat{\ga}_k$ of~$\ga_k$
the dimension of both $\ZZ_2$-vector spaces $\HC^{\cyl}_{\even}(Y,\hat{\ga}_k;h)$ 
and $\HC^{\cyl}_{\odd}(Y,\hat{\ga}_k;h)$ is~$k$.

\s
Our next goal is to distinguish contact structures on spheres.
Then every closed Reeb orbit is contractible, and so 
the cylindrical contact homology is not defined in general.

To proceed, we first notice that if $(Y,\ga)$ is a sphere, 
or more generally a contact manifold with vanishing homotopy groups $\pi_1(Y)$ and~$\pi_2(Y)$, 
we can associate to $x \in \cp_\ga$ a Conley--Zehnder index $\mu_{\sCZ}(x) \in \ZZ$
by choosing a unitary trivialisation of~$\xi$ along a disc in~$Y$ that bounds~$x$ and 
by looking at the linearized flow $d \phi_{\ga}^t \colon \xi_{x(0)} \to \xi_{x (t)}$, 
$t \in [0,T_x]$, in this trivialisation.

Another difference to the previous theory is that now we must orient the moduli spaces of 
$J$-holomorphic cylinders of dimensions~1 and~2.
Indeed, since all Reeb orbits are contractible, the $J$-holomorphic cylinders we look at
may now be multiply covered, and hence there are various ways to glue broken 
cylinders. 
To arrive at $\partial^2 =0$, multiply covered cylinders 
must therefore be counted with a certain rational weight, 
whence we must work with rational coefficients.
For this, the moduli spaces of $J$-holomorphic cylinders of dimensions~1 and~2 must be oriented;
and for this, all orbits~$x$ which are an even multiple of a simple orbit~$\underline x$
and for which $\mu_{\sCZ}(x) -\mu_{\sCZ}(\underline x)$ is odd must be excluded.
All other orbits are called good.

We thus now start with all good orbits $\cp_\ga^{\good}$.
A closer inspection shows that the holomorphic planes in the above accident can obstruct 
$\pp^2=0$ only if they converge to an orbit of index~$1$.
Further, the resulting homology $\HC_*(Y,\ga)$ is independent of the non-degenerate choice
of~$\ga$ if there are no Reeb orbits of index 0 or $-1$. 
For spheres we one can thus define the contact homology $\HC_*(Y,\xi)$
provided that there exists a non-degenerate contact form for~$\xi$ 
all of whose Reeb orbits have Conley--Zehnder index different from $-1,0,1$.
This turns out to be enough to prove the following:

\begin{theorem} \label{t:ust}
For every $m \geqslant 1$ the sphere $S^{4m+1}$ with its usual smooth structure
carries infinitely many non-diffeomorphic contact structures.
\end{theorem}

\ni 
{\it Outline of the proof.}
To fix the ideas we take $m=1$.
For natural numbers $a_0, \dots, a_3 \geqslant 2$ define 
the Brieskorn manifold $\Sigma (\aa)$ as the intersection of the
(in the origin singular) hypersurface 
$$
z_0^{a_0} + z_1^{a_1} + z_2^{a_2} + z_3^{a_3} \,=\, 0
$$
in $\CC^4$ with the unit sphere $S^7 \subset \CC^4$.
This is a smooth manifold, and the real form
\begin{equation} \label{e:aaa}
\ga_{\aa} \,=\,\frac{i}8 \sum_{j=0}^3 a_j (z_j d \overline{z_j} - \overline{z_j} dz_j) 
\end{equation}
restricts to a contact form on $\Sigma (\aa)$.
From now on we take $a_0 = p \equiv \pm 1 \pmod 8$ and $a_1 = a_2=a_3 =2$.
Then $\Sigma (\aa)$ is diffeomorphic to the standard~$S^5$.
We wish to prove that the contact structures $\xi_p$ defined by $\ga_p := \ga_{(p,2,2,2)}$
are mutually non-diffeomorphic.

The flow of the Reeb vector field 
$R_{\ga_p} = 2i \bigl( \frac{2}{p} z_0, z_1, z_2,z_3 \bigr)$ is periodic, 
and so $\ga_p$ is degenerate.
However, for $\gve \in (0,1)$ the function $f_\gve \colon \CC^4 \to \RR$
defined by
$$
f_\gve (z) \,=\, |z|^2 + 2\gve \,\mbox{Im} (z_2 \overline{z_3})
$$
is positive on~$\Sigma (\aa)$, and for irrational $\gve$ the contact form
$f_\gve \,\ga_p$ is non-degenerate, with only three simple closed orbits.
The Conley--Zehnder indices of these orbits and their iterates
turn out to be even and at least~$2$, 
and taking $\gve$ small one finds that the number of these orbits in degree~$2k$ for 
$2k = 2, 4, 6, \dots, p-1, p+1$ is
$$
1,\; 2,\; 2, \; \dots, \; 2, \; 1 .
$$
Hence the homologies $\HC_*^{\cyl} (\Sigma (\aa), f_\gve \,\ga_p)$ are defined, 
are an invariant of~$\xi_p$, 
and are all different.
\proofend

\begin{figure}[ht] 
 \begin{center}
   \psfrag{q1}{$q_1$}
   \psfrag{q2}{$q_2$}
   \psfrag{q3}{$q_3$}
   \psfrag{p}{$p$}
   \psfrag{u}{$u$}
   \psfrag{x}{$x$}
   \psfrag{y1}{$y_1$}
   \psfrag{y2}{$y_2$}
   \psfrag{y3}{$y_3$}
   \psfrag{c}{$\cong$}
   \leavevmode\epsfbox{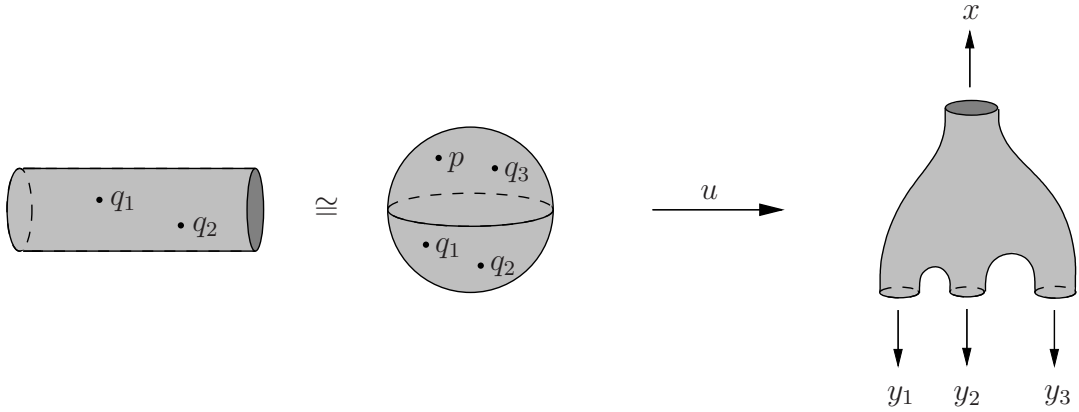}
 \end{center}
 \caption{Curves used in the construction of $\HC_*$} 
 \label{fig.ch}
\end{figure}

In the previous example we were a bit lucky with the indices. 
To obtain a more robust homology for contact manifolds,
we return to the accident described in Figure~\ref{fig.plane}.
It suggests how to proceed in a more general case:
Given $x,y$ in $\cp_{\ga}^{\good}$ we take $k$ punctures $q_1, \dots, q_k$ on the cylinder
$\RR \times \TT$ and look at holomorphic maps 
$u \colon (\RR \times \TT) \setminus \{ q_1, \dots, q_k \} \to Y \times \RR$ 
that are positively asymptotic to~$x$ as $s \to -\infty$, 
negatively asymptotic to $y$ as $s \to + \infty$,
and negatively asymptotic to some orbit $y_j \in \cp_\ga$ near~$q_j$.
And we do this for all integers~$k \geqslant 0$.
In this way, the problematic configurations are included in the picture from the outset.
Counting such punctured cylinders leads to the so-called linearized contact homology
mentioned in the notes to Section~\ref{s:symphom}.
Or, giving up the special role of~$y$, 
we fix $k \geqslant 1$, points $q, p_1, \dots, p_k$ on the sphere $S^2$, 
and a complex structure~$j$ on the $k+1$ fold punctured sphere 
$\dot S_k := S^2 \setminus \{q, p_1, \dots, p_k \}$,
and look at $(j,J)$-holomorphic maps 
$u \colon \dot S_k \to Y \times \RR$ such that, 
in holomorphic polar coordinates 
$(s,t) \in (0,\gve) \times \TT$ around the punctures, 
$$
\lim_{s \to 0} u_{\RR}(s,t) = +\infty \,\mbox{ at~$q$} 
\quad \mbox{ and } \quad
\lim_{s \to 0} u_{\RR}(s,t) = -\infty \,\mbox{ at~$p_j$} 
$$
and
$$
\lim_{s \to 0} u_Y(s,t) = x (-T_xt) \,\mbox{ at~$q$} 
\quad \mbox{ and } \quad
\lim_{s \to 0} u_Y(s,t) = y_j (T_jt) \,\mbox{ at~$p_j$} 
$$
for some $y_j \in \cp_\ga^{\good}$ of period~$T_j$, see Figure~\ref{fig.ch}.
For $x, y_1, \dots, y_k \in \cp_\ga^{\good}$ define $\cm (x; y_1, \dots, y_k)$
as the space of maps as above, where we allow all almost complex structures
on $\dot S_k$ but mod out by biholomorphisms fixing the punctures.

%
%

Set $\mu (x) := \mu_{\sCZ} (x)+n-3$, where $\dim Y = 2n-1$,
and consider the unital free super-commutative graded algebra $\ca_\ga$ over~$\QQ$
generated by the orbits in~$\cp_\ga^{\good}$.
Here the unit is a formal element in degree zero, 
the grading of the generators is given by~$\mu$, 
and super-commutative means that $xy = (-1)^{\mu(x) \mu(y)} yx$.
In other words, $\ca_\ga$ is a polynomial algebra on generators of even degree
and an exterior algebra on generators of odd degree.
The differential 
$$
\pp x \,=\, \sum_{k \geqslant 1} \sum_{y_1, \dots, y_k} \nu (x;y_1, \dots, y_k) \, y_1 \cdots y_k
$$
on generators of this algebra is defined by a certain weighted count $\nu (x;y_1, \dots, y_k) \in \QQ$
of the 1-dimensional components of $\cm (x;y_1, \dots, y_k)$, 
and is extended to products by the graded Leibnitz rule.
Then $\pp^2 =0$, and the {\it contact homology}\/ $\HC_* (Y,\ga)$ is independent of 
the generic choice of~$J$ and the choice of a non-degenerate~$\ga$, 
and is an invariant of the contact structure~$\xi$.

At the chain level, $(\ca_a, \pp)$ is a differential graded algebra, 
and the contact homology is an algebra, not just a vector space. 
It is thus not surprising that one obtains a theory with deep algebraic structures
if one takes into account all $J$-holomorphic curves, of finitely many positive and finitely many
negative punctures, that possibly have genus, cf.\ Figure~\ref{fig:curvesgG}.
The resulting homology is called the full SFT of~$(Y,\xi)$.

\subsection{Embedded contact homology} \label{s:ech}

Embedded contact homology is a version of SFT for three dimensional contact manifolds,
whose power comes from its relation to Seiberg--Witten theory.
We sketch its construction, and give a few of its many applications to symplectic geometry and dynamics in low dimensions.

\m \ni
{\bf Outline of the construction.}\
Embedded contact homology can be constructed for any closed 3-manifold with a contact form,
but for simplicity we here only look at the sphere~$S^3$ with a contact form~$\alpha$
whose kernel is the standard contact structure~$\xi_0$ from Example~\ref{ex:star}.
In other words, we look at the boundary of a starshaped domain $U \subset \RR^4$ with contact form~$\lambda_0 |_{\partial U}$.
We also assume that $\alpha$ is non-degenerate. 
As generators of the chain complex $\ECC_*(S^3,\alpha)$ one takes 
orbit~{\it sets}, namely finite collections $\Xi = \{(x_i, m_i)\}$, $m_i \in \NN$, 
where the $x_i$ are distinct simple closed Reeb orbits on~$S^3$, and $m_i$ indicates the multiplicity of~$x_i$.
One also requires that the $(x_i,m_i)$ are ``good'', 
which this time means that $m_i=1$ in case that the eigenvalues of the linearized return map 
$d\gf_{\alpha}^{T_i} \colon \xi_{x_i(0)} \to \xi_{x_i(T_i)}$ are real.

\begin{figure}[ht] 
 \begin{center}
   \psfrag{G}{$\Xi$}
   \psfrag{Gi}{$\Xi_j$}
   \psfrag{1}{{\footnotesize$1 \times$}}
   \psfrag{2}{\footnotesize{$3 \times$}}
   \psfrag{s-}{\footnotesize$s \to -\infty$}
   \psfrag{s+}{\footnotesize$s \to +\infty$}
   \leavevmode\epsfbox{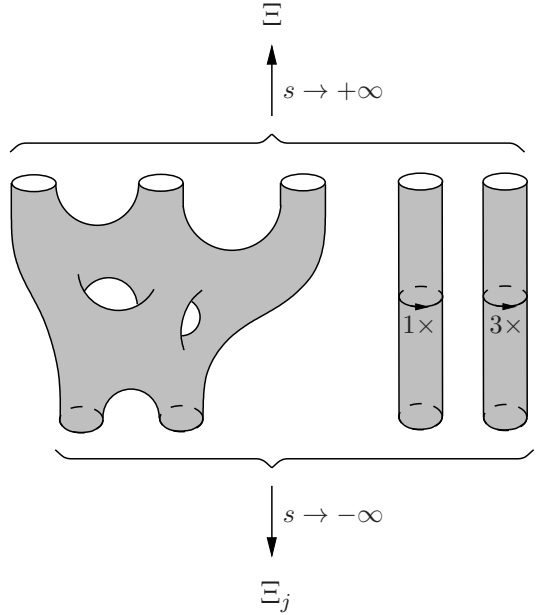}
 \end{center}
 \caption{Curves used in the construction of $\ECH_*$} 
 \label{fig:curvesgG}
\end{figure}

Take an almost complex structure on the symplectisation $S^3 \times \RR$ as in Section~\ref{s:ch}.
The boundary operator is now defined by counting possibly disconnected $J$-holomorphic curves in~$S^3 \times \RR$ 
with ends asymptotic to orbit sets $\Xi$ and~$\Xi_j$. 
The coefficient $\nu (\Xi, \Xi_j)$ of $\pp \1 \Xi = \sum_j \nu (\Xi,\Xi_j)\, \Xi_j$ counts mod~2
the 1-dimensional components of $\RR$-families of curves 
as indicated in Figure~\ref{fig:curvesgG}:
There is exactly one ``main curve'', that is embedded (whence the name),
may have genus and is not an $s$-invariant cylinder, 
and besides that there can be finitely many $s$-invariant cylinders that may be multiply covered.
The index of a generator~$\Xi$ is also different to the one in other Floer homologies
and much subtler a story.

%
%

The resulting homology is called the embedded contact homology $\ECH_*(S^3,\alpha)$. 
It turns out to be independent of~$\alpha$,
\begin{equation} \label{e:ECH}
\ECH_* (S^3,\alpha) \,=\, 
\left\{\begin{array} {lcl}
\ZZ_2  & &  \mbox{if $* \geqslant 0$ is even}, \\ [0.0em] 
0     & &   \mbox{otherwise} .
\end{array}\right.
\end{equation}

There are several reasons why ECH can be constructed and computed, 
while the full SFT is not yet fully established and much harder to compute.
One reason is that in ECH the symplectic manifold is four-dimensional.
In this dimension, $J$-holomorphic curves can be better understood than in higher dimensions, 
for instance because of ``positivity of intersections'',
which says that in a symplectic four-manifold two different $J$-holomorphic curves 
intersect positively at every intersection point.
Further, transversality in ECH is easier to achieve than in SFT, because the main curve is embedded.  

But the most important difference between ECH and the Floer homologies we have seen earlier is 
its tight relation to a gauge-theoretic homology. 
In fact, most choices in the construction of ECH are made in such a way that this relation exists, 
and this relation is key for establishing most properties of ECH:
The gauge-theoretic Seiberg--Witten equations can be written down in both dimension~three and four, 
hence for $S^3$ and~$S^3 \times \RR$.
The Seiberg--Witten--Floer homology $\SWF (S^3,\xi_0)$ is a Floer homology in which the generators are solutions 
of the Seiberg--Witten equations on~$S^3$ and the connecting orbits (defining the boundary operator)
are solutions of the Seiberg--Witten equations on~$S^3 \times \RR$.
Taubes proved in~\cite{Ta07} that $\SWF_* (S^3,\xi_0)$ is isomorphic to $\ECH_* (S^3,\alpha)$.
This implies that $\ECH_*(S^3,\alpha)$ does not depend on the choice of the almost complex structure used 
in its construction,
nor on the choice of the contact form~$\alpha$ for~$\xi_0$, 
and computing $\SWF_* (S^3,\xi_0)$ one finds~\eqref{e:ECH}.
The basic idea of Taubes' proof is as follows:
A solution to the SW-equation is a pair $(A, \psi)$, where $A$ is a connection 
and $\psi$ is a section of a spinor bundle over~$S^3$.
Deforming the SW-equation along $d\alpha$ one obtains a sequence of solutions 
whose spinor components are close to zero only
on a set that approximates more and more a collection of periodic Reeb orbits.
From this one constructs a map of chain complexes that descends to 
an isomorphism in homology.

\b \ni
{\bf ECH capacities.}
Let $U$ be a bounded domain in~$\RR^4$ that is starshaped with respect to the origin and
has smooth boundary~$\pp U$.
We also assume that the restriction of $\lambda_0$ to~$\pp U$ is non-degenerate.
The homology $\ECH_* (\pp U) = \ECH_*(\pp U,\lambda_0 |_{\partial U})$ 
comes with a filtration by action, where the action of a generator $\Xi = \{ (x_i,m_i) \}$
is defined by $\ca (\Xi) = \sum_i m_i \int_{x_i} \lambda_0$.
The inclusion $\ECC_*^a(\pp U) \to \ECC_*(\pp U)$ of the complex generated by those~$\Xi$
with $\ca (\Xi) \leqslant a$ descends to maps $\iota^a_*$ in homology. 
Let $\gg_k$ be the generator of $\ECH_{2k}(\pp U) \cong \ZZ_2$
and similar to~\eqref{def.cga} define
\begin{equation} \label{def.ck}
c_k(U) \,=\, \inf \left\{ a \mid \gg_k \in \im \iota^a_{2k} \right\}.
\end{equation}
In words, the $k$-th ECH capacity of~$U$ is the smallest $a$ such that the generator~$\gg_k$ of $\ECH_{2k}(\pp U)$
can be represented by a sum of orbit sets $\left\{ (x_i,m_i) \right\}$
with $\sum_i m_i \, \int_{x_i}\lambda_0 \leqslant a$.
For arbitrary starshaped domains with possibly non-smooth boundaries (such as polydiscs), 
these capacities are defined by approximation.

The numbers $c_k (U)$ are monotone increasing, 
$$
0 = c_0 (U) \leqslant c_1(U) \leqslant c_2(U) \leqslant \dots
$$
Moreover, they are monotone in the sense that $c_k(U) \leqslant c_k(V)$
if there exists a symplectic embedding $\gf \colon U \to V$,
since then there are natural maps $\ECH_*^a (\pp V) \to \ECH_*^a (\gf (\pp U))$
with properties similar to the transfer maps in symplectic homology.

\begin{example} \label{ex:ellECH}
{\rm
Let $U$ be an ellipsoid $\E (a,b)$ with $a/b$ irrational.
Then the only embedded closed characteristics on~$\pp U$ are the oriented circles $x_a = \{z_2=0\}$ and $x_b = \{z_1=0\}$.
A generator of the ECH chain complex has the form $(x_a,m) \cup (x_b,n)$ with $m,n \in \NN_{\geqslant 0}$.
Its action is $am+bn$.
It turns out that every generator has even index, hence the differential~$\pp$ vanishes identically.
Further, in each degree there is exactly one generator
(this yields another proof of~\eqref{e:ECH}, since $\ECH_*(\pp U)$ does not depend on the starshaped domain~$U$)
and the index is monotone with respect to the action. 
It follows that 
\begin{equation} \label{e:cN}
c_k (\E(a,b)) \,=\, N_k(a,b) ,
\end{equation}
where $(N_k(a,b))_{k \geqslant 0}$ is the sequence of numbers formed by arranging all the linear combinations 
$ma + nb$ with $m,n \geqslant 0$ in nondecreasing order (with repetitions). 
For rational ellipsoids, the identity~\eqref{e:cN} now follows from the monotonicity of the capacities and by taking inner
and outer approximations by irrational ellipsoids.
}
\end{example}

\ni
{\bf Application 1: Symplectic embedding obstructions for domains in~$\RR^4$}

\s 
\ni
Since $c_1(\B^4(a)) = a$ and $c_1(\E(b_1,b_2)) = b_1$ whenever $b_1 \leqslant b_2$,
the monotonicity of~$c_1$ implies another proof of the Nonsqueezing Theorem~\ref{t:Gromov}.
But ECH capacities have much finer applications to symplectic embedding problems.
Consider for instance the question of when $\E(a_1,a_2)$ symplectically embeds into~$\E(b_1,b_2)$.
ECH capacities provide a complete set of invariants for this problem:

\begin{theorem} \label{t:ellsharp}
There exists a symplectic embedding
$\E (a_1,a_2) \to \E (b_1,b_2)$ if and only if $N_k (a_1,a_2) \leqslant N_k (b_1,b_2)$ for all $k \geqslant 1$.
\end{theorem}

\b \ni
{\bf Application 2: Hamiltonian and contact closing lemmas}

\s \ni
For this application we need the extension of ECH capacities to Liouville domains~$(W,\omega)$.
Like the numbers $c_k(U)$ defined in~\eqref{def.ck}, the numbers $c_k(W)$ belong to the
multispectrum $\multispec (W)$, namely to the set of finite sums of actions of
closed Reeb orbits on~$\pp W$.
It is remarkable that from the capacities~$c_k(W)$ one can recover the volume of the domain~$W$.
Indeed, it has been proved in~\cite{CHR15} that
\begin{equation} \label{e:ckvol}
\lim_{k \to \infty} \frac{(c_k(W))^2}{k} \,=\, 4  \Vol (W) \,=\, 2 \int_W \omega \wedge \omega .
\end{equation}

\begin{theorem} 
[\bf $C^\infty$ Closing Lemma for geodesic flows on surfaces]  \label{t:Irie1}
Let $\Sigma$ be a closed surface with Riemannian metric~$g$.
Then for any open set $U \subset \Sigma$ and any $\gve >0$ there exists a 
smooth function $f \colon \Sigma \to \RR$ with $f-1$ supported in~$U$ and $d_{C^\infty} (f,1) < \gve$ 
such that for the conformal Riemannian metric $f g$ there exists a closed geodesic passing through~$U$. 
\end{theorem}

In the theorem, the distance is defined by
$$
d_{C^\infty} (f,g) \,=\, \sum_{k=0}^\infty 2^{-k}  \frac{\|f-g\|_{C^k}}{1+\|f-g\|_{C^k}} .
$$

\m \ni
{\it Idea of the proof.}
Consider the co-disc bundle 
$D^* (\Sigma,g) = \left\{ (q,p) \in T^*\Sigma \mid \|p\|_g \leqslant 1 \right\}$
with the usual symplectic form $d \lambda_{\can}$.
The multi-spectrum of this Liouville domain is the set of sums of lengths of closed geodesics, 
and by formula~\eqref{e:ckvol}
\begin{equation} \label{e:lim2vol}
\lim_{k \to \infty} \frac{\bigl( c_k(D^*(\Sigma,g)) \bigr)^2}{k} \,=\, 4 \pi \Vol (\Sigma, g) .
\end{equation} 
Now choose a smooth function $h \colon \Sigma \to \RR_{\geqslant 0}$ that does not vanish identically
but is supported in $U$, and such that $d_{C^\infty}(h,0) \leqslant \gve$.
For the metrics $g_t := (1+th)g$,
the function $t \mapsto \Vol (\Sigma, g_t)$ is then strictly increasing.
We claim that there exists $t \in [0,1]$ such that for~$g_t$
there exists a closed geodesic passing through~$U$.
If not, the metrics $g_t$ all have the same closed geodesics and hence the same multispectrum $\sigma \subset \RR$. 
Since ECH capacities depend continuously on~$t$ and since $\sigma$ has zero Lebesgue measure, 
$c_k(D^* (\Sigma,g_t))$ does not depend on~$t$, for each~$k$. 
Hence~\eqref{e:lim2vol} shows that $\Vol (\Sigma,g_t)$ does not depend on~$t$, which is not true. 
\proofend

Theorem \ref{t:Irie1} implies that for a $C^\infty$-generic subset of the space of Riemannian metrics on a closed surface~$\Sigma$, 
the set of closed geodesics is dense in~$\Sigma$.

A similar closing result holds for Reeb flows on a closed contact 3-manifold.
Combining this result with a suspension construction, one obtains

\begin{theorem} \label{t:closing}
Let $\Sigma$ be a closed oriented surface with an area form~$\go$,
and let $\gf$ be a Hamiltonian diffeomorphism of~$\Sigma$.
Then for any non-empty open set $U \subset \Sigma$ there exists a sequence of Hamiltonian diffeomorphisms
$(\gf_j)$ that converges to~$\gf$ in the $C^\infty$ topology 
such that every $\gf_j$ has a periodic point in~$U$.
\end{theorem}

If follows that for a closed surface~$\Sigma$ with area form~$\go$,
the group $\Ham (\Sigma,\go)$ contains a $C^\infty$-generic subset of diffeomorphisms with dense set of periodic points.

\begin{corollary} 
[\bf $C^\infty$ Hamiltonian closing lemma for surfaces]
\label{cor:closing}

For every point $p \in \Sigma$ and every $\gf \in \Ham (\Sigma,\go)$ there exists a sequence $(\Phi_j) \subset \Ham (\Sigma,\go)$
that converges to $\gf$ in the $C^\infty$ topology and such that $p$ is a periodic point of every $\Phi_j$.
\end{corollary}

\proof
By Theorem~\ref{t:closing} we find a sequence $(\gf_j) \subset \Ham (\Sigma,\go)$ with periodic points $p_j$
such that $\gf_j \to \gf$ in the $C^\infty$ topology and $p_j \to p$ in~$\Sigma$. 
It is easy to construct Hamiltonian diffeomorphisms $\tau_j$ with $\tau_j(p) =p_j$ and such that
$\tau_j \to \id$ and $\tau_j^{-1} \to \id$ in~$C^\infty$.
Then $\tau_j^{-1} \circ \gf_j \circ \tau_j \to \gf$,
and if $\gf_j^{\ell_j}(p_j) = p_j$, then $\Phi_j^{\ell_j}(p) =p$.
\proofend

For the $C^1$ topology, Corollary~\ref{cor:closing} is the famous Hamiltonian $C^1$ closing lemma of Pugh and Robinson~\cite{PuRo83}, 
that holds in arbitrary dimension.
Corollary~\ref{cor:closing} answers the Hamiltonian case of Problem~10 in Smale's list~\cite{Sma98} 
for surfaces and the $C^\infty$ topology.

\m \ni
{\bf Historical notes and bibliography}

\s \ni
{\bf Contact homology} was first sketched in~\cite{El98},
and SFT was outlined by Eliashberg, Givental and Hofer in~\cite{ElGiHo00}.
More details can be found in Wendl's lecture notes~\cite{We16}.
Surveys on contact homologies and SFT are~\cite{bo09, el07}.
The compactness theorem necessary to establish symplectic field theory is proven in~\cite{boelhowyze03}.
Transversality for contact homology is established in~\cite{Par15},
but for general~SFT relies on the polyfold theory from~\cite{Ho04, HWZ17}.

Theorem~\ref{t:CZ} (i) was first proved by Giroux~\cite{gir94} and Kanda~\cite{kan97}
by methods specific to dimension~three.
For a detailed proof by cylindrical contact homology we refer to Lecture~10 in~\cite{We16}.
Cylindrical contact homology can distinguish many contact structures also in higher dimensions.
Several examples can be found in~\cite{bou02}.
Among them is Ustilovsky's Theorem~\ref{t:ust} from ~\cite{ust99},
that was the first application of contact homology.
A Morse--Bott version of contact homology is given in~\cite{bou02}.
This has the advantage that in many interesting situations one does not need to find a non-degenerate contact form,
but can directly work with a symmetric and particularly nice contact form like~$\ga_{\aa}$ in~\eqref{e:aaa},
see~\cite{bou02, vko08} for many examples.

Theorem~\ref{t:CZ} (ii) on the fundamental group of the components of the space of contact structures 
on the 3-torus was first proved by Geiges and Gonzalo~\cite{gego04}.
But again, cylindrical contact homology can be used to prove such results also in higher dimensions. 
For instance, Bourgeois showed in~\cite{bou06} that 
for any closed $2n$-dimensional manifold~$Q$ with $\H_1(Q;\ZZ) \neq 0$ that admits a Riemannian metric of negative curvature,
the homotopy group~$\pi_{2n-1}$ of the space of contact structures 
on~$S^*Q$ in the component of~$\xi_{\can}$ contains an infinite cyclic subgroup.

Other applications of cylindrical contact homology are 
a contact version of the nonsqueezing theorem~\cite{elkipo06, fra16}, 
and proofs of cases of the Weinstein conjecture and of versions of the Conley conjecture 
for Reeb flows, see the survey~\cite{GiGu15}. 
Rational SFT (namely SFT where all $J$-curves are punctured spheres) 
is used in \cite{ElGiHo00} and~\cite{cimo14}
to find obstructions to the existence of Lagrangian submanifolds in 
symplectic manifolds that contain many holomorphic spheres (like complex projective space).
The full SFT package has also found applications to questions in contact topology, see~\cite{lawe11}.

In the same way as Lagrangian Floer homology is a relative version of Hamiltonian Floer homology, 
there are contact homologies for Legendrian submanifolds, 
see for instance~\cite{eketsu07}.
Chekanov and Eliashberg used them to distinguish Legendrian knots with the same 
``classical'' invariants~\cite{ch02}. 
They also give rise to very fine invariants of smooth knots in~$\RR^3$, see~\cite{CELN16, ekngsh16, ng05}.

\s \ni
{\bf Embedded contact homology} is due to Hutchings.
His lectures~\cite{Hu14} also contain a 
comparison of ECH and~SFT.
Theorem~\ref{t:ellsharp} was first proved by McDuff~\cite{Mc11}.
Her proof does not explicitely use ECH capacities, but is inspired
by their properties, and the proof using ECH capacities is shorter~\cite{Hu11b}. 
For many other applications of ECH capacities to symplectic embedding problems
see the original~\cite{Hu11} and Section~12 of~\cite{Sch17}.
The two surveys~\cite{CHLS, Sch17} describe motivations and applications of symplectic embedding 
problems to questions outside of symplectic geometry.
Theorem~\ref{t:Irie1} is due to Irie~\cite{Ir15} and Theorem~\ref{t:closing} with its Corollary~\ref{cor:closing}
to Asaoka and Irie~\cite{AsIr16}.

The embedded contact homology $\ECH_* (Y,\alpha)$ of a closed 3-dimensional contact manifold~$(Y,\alpha)$ is constructed along the same lines as
in the special case of a contact form on~$S^3$ outlined above,
and Taubes proved in~\cite{Ta07} that $\ECH_* (Y,\alpha)$ is isomorphic to 
the Seiberg--Witten--Floer homology $\SWF_* (Y,\xi)$, where $\xi = \ker \alpha$.
Kronheimer--Mrowka~\cite{KrMr07} showed that $\SWF_* (Y,\xi)$ is non-zero in 
infinitely many degrees.
Since $\ECH_* (Y,\alpha)$ is generated by closed Reeb orbits and the empty set, 
there must exist at least one closed Reeb orbit 
(otherwise, $\ECH_* (Y,\alpha)$ would vanish in all degrees $\neq 0$). 
This proves the general Weinstein conjecture in dimension~3.
For a survey on Taubes' proof see~\cite{Hu10}.

The previous result can be improved by using ECH, see~\cite{CGHu16}:
{\it Every Reeb flow on a closed contact 3-manifold has at least two simple periodic orbits.}
Take, for instance, a contact form $\alpha$ on~$(S^3,\xi_0)$.
By the proof of the Weinstein conjecture there exists a simple 
periodic Reeb orbit~$x$. 
Let $T = \ca (x)$ be its period, and 
let $U \subset \RR^4$ be the starshaped domain with $(\pp U, \lambda_0|_{\pp U}) \cong (S^3, \alpha)$.
If $x$ is the only simple periodic Reeb orbit, 
then the definition~\eqref{def.ck} shows that
$c_k(U) = n_k T$ for a strictly increasing sequence of natural numbers~$n_k$.
In particular, $n_k \geqslant k$, and so
$$
\lim_{k \to \infty} \frac{(c_k^2(U))^2}{k} \,\geqslant\, 
\lim_{k \to \infty} \frac{k^2 T^2}{k} \,=\, \infty,
$$
in contradiction to~\eqref{e:ckvol}.
Taking the boundary of an irrational ellipsoid $\E (a,b)$ we see that the above result result is sometimes optimal.
Is it true that every Reeb flow on a closed contact 3-manifold has
either two or infinitely many simple periodic orbits, 
and that the first possibility
can happen only on $S^3$ and on lens spaces?
Partial affirmative answers to this question are given in~\cite{CGHuPo17, HoWyZe03, HuTa09}.

Recall from Section~\ref{s:basicsII} that Arnol'd conjectured 
that every Legendrian knot in~$\RR^3$ with the standard 
contact structure $\xi_{\st} = \ker (dz+xdy)$ has a Reeb chord, 
for every contact form for~$\xi_{\st}$.
Since $\xi_{\st}$ on~$\RR^3$ is contactomorphic to the contact structure 
$\xi_0 = \ker \lambda_0$ on~$S^3$ without one point, 
this conjecture follows from the analogous conjecture for~$(S^3,\xi_0)$,
which was proved by Mohnke~\cite{Moh01}.
Hutchings and Taubes~\cite{HuTa11,HuTa13} proved the Arnol'd chord conjecture 
for all closed contact 3-manifolds:
{\it Every Legendrian knot in a closed contact 3-manifold $(Y,\xi)$ has a Reeb chord, 
for any contact form for~$\xi$.}
%
The proof uses again properties of the Floer homologies 
$\SWF_* (Y,\xi) \cong \ECH_* (Y,\alpha)$.


\section{Applications to topology} \label{s:top}

The classical applications of Floer homology are to Hamiltonian dynamics and symplectic geometry. 
But Floer theory is rather a concept than a collection of theorems.
It is thus not surprising that Floer homology found applications in many other fields. 
We have already seen a topological application of Lagrangian Floer homology in cotangent bundles  (Corollary~\ref{c:lens}), 
and as mentioned in the notes to the previous section 
Legendrian contact homology gives rise to powerful knot invariants. 
In this section we describe a few topological applications of other Floer homologies.

\s
Already Floer~\cite{Flo88:3man} associated to every 3-dimensional oriented integral homology sphere~$Y$
its ``instanton Floer homology'' $I_*(Y)$ by looking at the Chern--Simons functional 
on the space of $\SU (2)$ connections on~$Y$.
The critical points are connections with vanishing curvature (modulo gauge), 
and the gradient flow lines are solutions of the anti-self-dual Yang--Mills equation on $Y \times \RR$.
These groups have been later defined for all closed 3-manifolds.
They are important also in the study of 4-manifolds, since the relative Donaldson invariants
of 4-manifolds with boundary take values in instanton Floer homology.
Since Floer's pioneering work,
many other Floer homologies have been associated to 3-manifolds.
Each of these Floer homologies has a gauge theoretic and a symplectic version, 
that are (at least conjecturally) isomorphic. 
For instance, take a closed 3-manifold~$Y$ and choose a Heegard decomposition $Y = U_0 \cup_{\Sigma} U_1$,
where $\Sigma \subset Y$ is a closed surface and the~$U_j$ are handle bodies. 
The space of $\SU (2)$ flat connections $M(\Sigma)$ is a symplectic manifold, and those flat connections
on~$\Sigma$ that extend to~$U_j$ form a Lagrangian submanifold~$L_j$ in~$M(\Sigma)$.
According to the Atiyah--Floer conjecture, the instanton homology~$I_*(Y)$ 
is isomorphic to the Lagrangian Floer homology $\HF_*(M(\Sigma),L_0,L_1)$, 
and this has been proved for mapping tori in~\cite{DoSa94}.
Similarly, the gauge-theoretic Seiberg--Witten--Floer homology $\SWF (Y,\xi)$ is isomorphic to the embedded contact homology 
$\ECH (Y,\alpha)$, cf.\ \S \ref{s:ech}.
The isomorphisms between the gauge-theoretic and the symplectic versions make these Floer homologies 
very powerful, since insights on one side can be transferred to the other side; 
they have many applications to low dimensional topology.
And recently, Manolescu~\cite{Ma16} gave a spectacular application to higher dimensions, 
disproving the triangulation conjecture.
His survey~\cite{Ma15} gives an excellent description of the Floer homologies for 3-manifolds
and of their topological applications.
We therefore simply state a few  results, 
that give an idea which kind of topological applications Floer homologies have.

\m
Let $\kk \subset S^3$ be a knot.
Remove a tubular neighbourhood $\nu(\kk) \cong S^1 \times D^2$ from $S^3$.
A longitude is a simple closed curve $\lambda$ on the 2-torus $\pp \nu (\kk)$ 
that is contractible in~$S^3 \setminus \nu (\kk)$, 
and a meridian is a simple closed curve $\mu \subset \pp \nu (\kk)$ 
that is contractible in~$\nu (\kk)$.
Both curves are unique up to isotopy in $\pp \nu (\kk)$. 
Now take relatively prime integers $p,q$, and glue back 
the solid torus $S^1 \times D^2$ along its boundary to~$\pp \nu (\kk)$
in such a way that the meridian $\{1\} \times S^1 \subset T^2$ 
is taken to a simple closed curve in the homology
class $p [\mu] + q [\lambda] \in H_1(\pp \nu (\kk))$.
The resulting oriented 3-manifold 
$$
S^3_{p/q}(\kk) \,=\, S^3 \setminus (\nu (\kk)) \cup_{\pp \nu (\kk)} (S^1 \times D^2) 
$$
is called the Dehn surgery of~$S^3$ on~$\kk$ with slope~$p/q$.

By a classical result of Lickorish and Wallace, every closed connected orientable 3-manifold can be obtained
by Dehn surgery on a collection of disjoint knots in~$S^3$. 
An important problem in 3-manifold topology is thus to understand which manifold can 
be obtained by Dehn surgery of~$S^3$ on which knot with which slope.
Let $\kk_0$ be the unknot. 

\begin{theorem} \label{t:P}
{\rm \bf (Property P for knots~\cite{KrMr04})}
If $\kk \neq \kk_0$, then $\pi_1 (S^3_{p/q}(\kk)) \neq 0$ for all $p,q$.
\end{theorem}

This theorem shows that at least by Dehn surgery on a single knot one cannot produce
a counterexample to the 3-dimensional Poincar\'e conjecture.
This conclusion is now superseded by Perelman's proof of the 3-dimensional Poincar\'e conjecture,
which also implies Theorem~\ref{t:P} 
since Dehn surgery on a non-trivial knot can never yield~$S^3$.

\begin{theorem} \label{t:KMO}
{\rm (\cite{KMOS07})}
Let $\kk \subset S^3$ be a knot. 
If there is an orientation preserving diffeomorphism $S^3_{p/q} (\kk) \to S^3_{p/q} (\kk_0)$
for some $p/q \in \QQ$, then $\kk = \kk_0$.
\end{theorem}

It is well known that $S^3_{p/q} (\kk_0) = S^1 \times S^2$ if $p=0$
and that $S^3_{p/q} (\kk_0)$ is the Lens space $L(p,q)$ if $p>0$. 
Since the lens space $L(2,q)$ is $\RP^3$ for all odd $q$, Theorem~\ref{t:KMO} thus implies 
that $\RP^3$ cannot be obtained by Dehn surgery on a non-trivial knot.

The next result shows that an oriented closed 3-manifold can be obtained by at most two different Dehn surgeries on the same knot. 

\begin{theorem}
{\rm (\cite{Wu11, NiWu15})}
Let $\kk \subset S^3$ be a knot. 
If there is an orientation preserving diffeomorphism 
$S^3_{p/q} (\kk) \to S^3_{p'/q'} (\kk)$
for $p/q \neq p'/q'$, then $p/q = -p'/q'$ and $q^2 \equiv -1 \pmod p$.
\end{theorem}

\m
Recall that an $n$-dimensional compact topological manifold is a compact Hausdorff space  
that is locally homeomorphic to~$\RR^n$.
A {\it triangulation}\/ of a compact topological manifold~$X$ is a homeomorphism $X \to K$
to a finite simplicial complex~$K$.
Kneser asked in 1924 whether every compact topological manifold admits a triangulation.
This is so in dimension $n \leqslant 3$ and also for compact smooth manifolds of arbitrary dimension.
On the other hand, Casson showed that there are 4-dimensional compact topological manifolds 
without triangulation.
This is also the case in higher dimensions:

\begin{theorem}
{\rm (\cite{Ma15})}
For every $n \geqslant 5$ there exist compact $n$-dimensional topological manifolds
that do not admit a triangulation. 
\end{theorem}

It was known that the triangulation problem in dimension $\geqslant 5$ is equivalent
to a problem in dimensions~3 and~4, and it was this problem that Manolescu
solved by means of a Floer theory.

\end{document}